\documentclass[11pt]{article}

\usepackage{subcaption}
\usepackage{amsmath, amsfonts, amsthm, amssymb}
\usepackage{mathabx, comment}
\usepackage{authblk, color, bm, graphicx}
\usepackage[margin = 1.5cm, font = small, labelfont = bf, labelsep = period]{caption}
\usepackage{xspace}
\usepackage[utf8]{inputenc}
\usepackage{multirow}
\usepackage{setspace}
\usepackage[round]{natbib}
\usepackage{mathtools}
\usepackage{tabularx}
\usepackage{hyperref}
\usepackage{enumitem}

\usepackage{algorithmicx, algorithm,algpseudocode}

\usepackage[dvipsnames]{xcolor}
\usepackage[framemethod=TikZ]{mdframed}

\usepackage{fullpage}[2cm]
\newtheorem{defi}{Definition}
\newtheorem{thm}{Theorem}
\newtheorem{prop}{Proposition}
\newtheorem{lem}{Lemma}
\newtheorem{obs}{Observation}
\newtheorem{coro}{Corollary}
\newtheorem{ex}{Example}
\newtheorem{rem}{Remark}

\newcommand{\eg}{\textit{e.g.}}

\title{An MILP-Based Solution Scheme for Factored \\ and Robust Factored Markov Decision Processes}

\author[1]{Huikang Liu}
\author[2]{Wolfram Wiesemann}
\author[3]{Man-Chung Yue}

\affil[1]{\small \textit{Shanghai Jiao Tong University, China,} \texttt{hkl1u@sjtu.edu.cn}}
\affil[2]{\small \textit{Imperial College London, United Kingdom,} \texttt{ww@imperial.ac.uk}}
\affil[3]{\small \textit{The University of Hong Kong, Hong Kong,} \texttt{mcyue@hku.hk}}

\begin{document}
	
\maketitle

\begin{abstract}
Factored Markov decision processes (MDPs) are a prominent paradigm within the artificial intelligence community for modeling and solving large-scale MDPs whose rewards and dynamics decompose into smaller, loosely interacting components. Through the use of dynamic Bayesian networks and context-specific independence, factored MDPs can achieve an exponential reduction in the state space of an MDP and thus scale to problem sizes that are beyond the reach of classical MDP algorithms. However, factored MDPs are typically solved using custom-designed algorithms that can require meticulous implementations and considerable fine-tuning. In this paper, we propose a mathematical programming approach to solving factored MDPs. In contrast to existing solution schemes, our approach leverages off-the-shelf solvers, which allows for a streamlined implementation and maintenance; it effectively capitalizes on the factored structure present in both state and action spaces; and it readily extends to the largely unexplored class of robust factored MDPs, whose transition kernels are only known to reside in a pre-specified ambiguity set. Our numerical experiments demonstrate the potential of our approach.

\textbf{Keywords:} Factored Markov decision processes; robust Markov decision processes.
\end{abstract}

\section{Introduction}

Dating back to the seminal work of \cite{Bellman}, Markov decision processes (MDPs) have emerged as one of the major paradigms to model, analyze and solve dynamic decision problems affected by uncertainty \citep{Puterman, Bertsekas}. Nevertheless, traditional MDPs are hampered by the \emph{curse of dimensionality} since their characterization of a problem's dynamics and reward structure tends to increase exponentially in the length of the problem's description. In recent decades, various research communities have devised strategies to circumvent the unfavorable scalability of traditional MDPs through approximation, decomposition or a blend of both principles.

A large number of physical, biological and social systems possess a nearly decomposable, hierarchic structure \citep{S62:architecture_complexity}, and the corresponding MDP representations---while large in scale---naturally embody aggregations of smaller components that evolve more or less independently. This structure is leveraged by \emph{factored MDPs} (FMDPs), a methodology popularized predominantly within the artificial intelligence community. Instead of adopting the flat state space representation of traditional MDPs that enumerates all states sequentially (such as $s = 1, 2, \ldots$), FMDPs characterize each state of an MDP via value assignments to state variables (such as ``\texttt{velocity} = \texttt{high}, \texttt{acceleration} = \texttt{low}''), and these assignments evolve and contribute to the system's rewards largely independently. This enables a compact representation of the stochastic evolution of FMDPs through dynamic Bayesian networks, which describe the one-step changes in the state variables as functions of a small number of ancestral state variables, as opposed to the entire state. FMDPs have found applications in, among others, resource allocation \citep{GPS02:app:resource_allocation, DD06:app:resource_allocation}, recommender systems \citep{TB14:app:recommender_systems}, power plant \citep{RIS04:app:power_plant, RSS09:app:power_plant} and micro grid management \citep{AC14:app:micro_grid, AC15:app:micro_grid}, manufacturing \citep{BDG00:app:manufacturing_robots, DSdB11:app:manufacturing_traffic}, infrastructure design \citep{HCZ17:app:infrastructure}, traffic management \citep{DSdB11:app:manufacturing_traffic} and robot task planning \citep{BDG00:app:manufacturing_robots}. Factored representations are also employed in reinforcement learning \citep{SH04:rl:factored, DSW06:rl:learning, OVR14:rl:near_optimal} and partially observable MDPs \citep{WPY05:pomdp:dialogue}.

Perhaps due to its artificial intelligence heritage, the solution of FMDPs has largely relied upon hand-designed algorithms that iteratively manipulate decision trees, algebraic decision diagrams and/or rule systems representing the dynamics of the system as well as the value function \citep{Boutellier, Hoey, ZP99:csi}. While there is no doubt that the resulting algorithms have significantly enhanced our comprehension of FMDPs, thereby facilitating numerous successful practical applications, these algorithms come with their own set of challenges. Most notably, they require considerable care during (re-)implementation, they tend not to exploit recent advancements in computational hardware (such as the parallelizability provided by contemporary multi-core architectures), and they demand\ considerable technical and domain-specific knowledge for fine-tuning their hyperparameters.

In this paper, we put forth a new framework that leverages mathematical programming to model and solve FMDPs. Instead of requiring the state evolution and reward structure to exhibit a sparse dependence on the current state, our model imposes a sparse dependence on a mixed-integer linear programming (MILP) representable feature vector derived from the current state and action. Compared to the existing literature on FMDPs, our MILP feature-based representation enjoys at least three advantages: \emph{(i)} it enables us to model richer classes of dependencies since the feature vector may exhibit a complex dependence on the entire state without impacting the scalability of our approach; \emph{(ii)} it allows us to achieve an exponential reduction not only in the state space but also in the action space; and \emph{(iii)} it unifies and significantly extends the mainstream approaches for capturing context-specific independence (which we elaborate upon later). To solve the resulting FMDP formulation, we combine a well-known approximate linear programming (LP) formulation for traditional MDPs, which replaces the exponentially large value function with a linear combination of basis functions, with a novel cutting plane algorithm that identifies violated constraints through the solution of an MILP (rather than sequentially probing exponentially many constraints). Our method naturally extends to factored representations of \emph{robust MDPs}, where the transition kernel governing the stochastic state evolution is only known to reside in a pre-specified ambiguity set, and where a policy is sought that performs best in view of the most adversarial transition kernel contained in the ambiguity set. We consider two representations for the ambiguity set: a non-factored representation that contains all transition kernels whose marginal projections along the individual state variables lie in pre-specified sets, and a factored representation that additionally imposes independence between the transitions of the state variables. We investigate the numerical performance of our algorithms in two case studies, which also compare them against alternative approaches from the literature.

The main contributions of this paper may be summarized as follows.
\begin{enumerate}
    \item[\emph{(i)}] We propose a new framework for FMDPs that is centered around MILP-representable features. Our methodology allows us to exploit factored structure in larger classes of MDPs, and it enables us to achieve an exponential reduction in both the state and the action spaces.
    \item[\emph{(ii)}] We show how an FMDP in our framework can be solved via an MILP-based cutting plane approach. In contrast to previous work, our approach relies on off-the-shelf solvers and avoids the need for hand-designed solution schemes.
    \item[\emph{(iii)}] We extend our framework to robust FMDPs, which are largely unexplored to date. We introduce two classes of ambiguity sets, we show how they relate to each other, and we demonstrate how robust FMDPs can be solved with our cutting plane approach.
\end{enumerate}
More broadly, we see this work as a step towards the ambitious goal of integrating the insights of two influential yet distinct views on dynamic decision-making under uncertainty: the artificial intelligence approach, which has been pivotal in introducing, examining and leveraging structure in MDPs, and the mathematical programming approach, which has been exceptionally successful in solving MDPs. The source code of our framework, together with all experimental results, is made available open-source as part of this work to facilitate reuse as well as computational comparisons.\footnote{Anonymous GitHub repository: \href{https://anonymous.4open.science/r/Factored-MDPs-A0F0/}{https://anonymous.4open.science/r/Factored-MDPs-A0F0/}. \label{fn:github}}

We are not the first ones to study FMDPs through the lens of mathematical programming. In a series of papers, \cite{GKP01:max_norm, Guestrin1, GVK02:context_specific, GKPV03:efficient_solution} and \cite{SP01:direct_value} devise compact LP formulations for an approximate value iteration and an approximate policy iteration as well as an approximate linear programming approach for FMDPs. The authors conclude that their approximate linear programming approach is superior as it exhibits better scaling while at the same time requiring less stringent assumptions about the structure of the problem. In contrast to our cutting plane scheme, their approach relies on non-serial dynamic programming \citep{BVB72:nonserial} to equivalently replace an exponential number of constraints with a smaller, hopefully polynomial, number. Computing the smallest number of equivalent constraints is known to be NP-hard, and there are instances for which no equivalent reformulation with a small number of constraints exists. In practice, problem-specific heuristics are employed to generate small equivalent representations without performance guarantees. The non-serial dynamic programming approach differs from our method in multiple aspects: \emph{(i)} it does not seem to accommodate MILP-based feature representations, which provide us with significant additional modelling flexibility, \emph{(ii)} it cannot leverage factored structure in the action space, and \emph{(iii)} its extension to robust FMDPs requires the solution of non-convex multilinear optimization problems for which currently neither global nor local optimality guarantees are available. We compare our cutting plane scheme with the non-serial dynamic programming approach in our numerical experiments.

\cite{DSdB11:app:manufacturing_traffic,DdBCS11:using_mp_to_solve_fmdp,DdBDS16:real_time_DP} develop several solution schemes for robust FMDPs. In contrast to the non-robust setting, however, these approaches require the solution of highly nonlinear optimization models as subproblems. The authors solve those problems heuristically using local optimization techniques, and they do not discuss any convergence guarantees. In summary, our comprehension of robust FMDPs is not as advanced as our knowledge of robust non-factored MDPs, for which a relatively complete theory exists (see, \emph{e.g.}, \citealp{iyengar2005robust,nilim2005robust,wiesemann2013robust} as well as the literature review of \citealp{HPW22:robust_phi_div}).

%- since those approaches are based on those from the previous paragraphs, they inherit their disadvantages as well.
%- Comparison to existing literature on solving robust FMDPs. 
%- ambiguity sets similar to our factored ambiguity sets
%\cite{DSdB11:app:manufacturing_traffic}:
%- require an additional rectangularity assumption on the ambiguity set (rectangularity across subsets of target states)
%- extend algebraic decision diagrams and the associated algorithms to robust setting
%- \cite{DdBDS16:real_time_DP} use asynchronous algorithms for the subclass of stochastic shortest path problems.
%\cite{DdBCS11:using_mp_to_solve_fmdp}:
%- they the approximate linear programming formulation from the previous paragraph but optimize simultaneously over the value function approximation and the transition kernel
%- they reduce the number of constraints using non-serial dynamic programming techniques from the previous paragraph.
%the approaches have in common that they require the solution of highly nonlinear optimization models as subproblems, and they solve them heuristically using local optimization techniques. No convergence guarantees are provided.
%since those approaches are based on those from the previous paragraphs, they inherit their disadvantages as well.

Beyond FMDPs, other frameworks have been proposed to alleviate the curse of dimensionality of MDPs. Approximate dynamic programming, for example, offers a very general methodology to control large-scale MDPs by approximating the value function through a linear combination of basis functions \citep{Bertsekas_Tsitsiklis, P11:adp}. In contrast to FMDPs, approximate dynamic programming does not rely on any specific problem structure. However, this generality implies that the associated solution approaches either require stronger assumptions (such as access to the stationary distribution under the unknown optimal policy) or provide weaker guarantees (such as asymptotic rather than finite convergence). Weakly coupled dynamic programming, on the other hand, studies systems that decompose into largely independently evolving constituent MDPs that are coupled by a small number of resource constraints \citep{Hawkins, Adelman_Mersereau, Bertsimas_Misic}. Weakly coupled dynamic programming is less expressive than FMDPs as it requires the constituent MDPs to evolve completely independently once a policy has been fixed. A notable exception is the recent work of \citet{BZ22:dps_with_signals}, which we compare ourselves against in the numerical experiments.

The remainder of the paper unfolds as follows. Section~\ref{sec:FMDP} defines FMDPs as well as sparse representations via feature vectors, and it introduces our MILP-based cutting plane solution scheme. Sections~\ref{sec:robust_FMDP} and~\ref{sec:robust_FMDP:factored} extend our methodology to non-factored and factored robust FMDPs, respectively. We report the results of numerical experiments in Section~\ref{sec:num_experiments}. Supplementary material and all proofs are provided in an electronic companion to maintain the conciseness of the main paper.

\textbf{Notation.} We define $\mathbb{B} = \{ 0, 1 \}$ and let $[n] = \{ 1, \ldots, n \}$ for any $n \in \mathbb{N}$. We denote by $\mathbf{1}$ the vector of all ones and by $\mathrm{e}_i$ the $i$-th canonical basis vector, respectively; in both cases, the dimension will be clear from the context. For a finite set $\mathcal{X} = [n]$, we denote by $\Delta (\mathcal{X})$ the set of all probability distributions over $\mathcal{X}$, which we interpret either as functions $\Delta = \{ (p : \mathcal{X} \rightarrow [0, 1]) \, : \, \sum_{x \in \mathcal{X}} p(x) = 1 \}$ or as vectors $\Delta = \{ p \in [0, 1]^n \, : \, \mathbf{1}^\top p = 1 \}$, depending on the context. We define $\mathbf{1} [\mathcal{E}]$ as the indicator function that takes the value $1$ ($0$) if the condition $\mathcal{E}$ is (not) satisfied.

%  and set $\mathbb{B} (\mathcal{X}) = \mathbb{B}^{| \mathcal{X} |}$ for any finite index set $\mathcal{X}$

\section{Factored Markov Decision Processes}\label{sec:FMDP}

A (non-factored) \emph{Markov decision process} (MDP) is defined by the tuple $(\mathcal{S}, \mathcal{A}, q, p, r, \gamma)$, where $\mathcal{S}$ and $\mathcal{A}$ represent the finite state and action spaces, $q \in \Delta (\mathcal{S})$ describes the distribution of the initial state, the transition probabilities $p : \mathcal{S} \times \mathcal{A} \rightarrow \Delta (\mathcal{S})$ characterize the stochastic one-step evolution of the system, $r : \mathcal{S} \times \mathcal{A} \rightarrow \mathbb{R}$ is the reward function, and $\gamma \in (0, 1)$ is the discount factor. MDPs do not capture any intrinsic structure that an application domain may possess, and they thus often scale exponentially in a natural representation of the problem.

\begin{ex}[Predictive Maintenance]\label{ex:machine_rep:initial}
    Consider a maintenance problem where multiple pieces of equipment are operated in sequence, and where the condition of upstream equipment impacts the deterioration rate of downstream equipment. This is the case, for example, \emph{(i)} in assembly lines where a poor calibration of the initial cutting step results in subtle defects that can impact the subsequent welding and soldering steps, \emph{(ii)} in the oil and gas industry, where poorly maintained drilling machines result in sediments and other impurities that lead to a faster wear and tear of downstream equipment such as pipelines, pumps, and processing units, and \emph{(iii)} in power plants, where poorly maintained boilers can generate steam at an incorrect pressure or temperature that causes excessive wear and tear at downstream equipment such as turbines.

    Assume that $10$ machines are operated in sequence, each of which is in one of $10$ states ranging from $1$ (malfunction) to $10$ (optimal performance). In each period, exactly $2$ machines should undergo maintenance. The problem can be modeled as a non-factored MDP with $10^{10}$ states and ${10 \choose 2} = 45$ actions, \mbox{which requires the consideration of $10^{10} \cdot 45 \cdot 10^{10} = 4.5 \cdot 10^{21}$ transitions probabilities.}
\end{ex}

Oftentimes, we can achieve an exponential reduction in the representation of a problem by adopting a factored representation of the state space \citep{BDG00:app:manufacturing_robots, GKPV03:efficient_solution}. 

\begin{defi}[Factored MDP]\label{def:fmdp}
    A \emph{factored  MDP} (FMDP) is defined by the tuple $(\mathcal{S}, \mathcal{A}, q, p, r, \gamma)$, where the finite state space $\mathcal{S} = \bigtimes_{n \in [N]} \mathcal{S}_n$ constitutes a product of binary sub-state spaces $\mathcal{S}_n \subseteq \mathbb{B}^{S_n}$, the finite action space is $\mathcal{A} \subseteq \mathbb{B}^A$, the initial distribution is $q = \{ q_n \}_{n \in [N]}$ with $q_n \in \Delta (\mathcal{S}_n)$, the transition kernel is $p = \{ p_n \}_{n \in [N]}$ with $p_n : \mathcal{S} \times \mathcal{A} \rightarrow \Delta (\mathcal{S}_n)$, the reward function is $r = \{ r_j \}_{j \in [J]}$ with $r_j : \mathcal{S} \times \mathcal{A} \rightarrow \mathbb{R}$, and the discount factor is $\gamma \in (0, 1)$.
\end{defi}

The restriction to binary sub-state and action spaces in Definition~\ref{def:fmdp} is standard in the FMDP literature, and it is taken in view of the solution methods that we develop later. That said, all of our theory readily extends to general finite sub-state and action spaces.

A deterministic and stationary (\emph{i.e.}, time-invariant) \emph{policy} $\pi : \mathcal{S} \rightarrow \mathcal{A}$ for an FMDP assigns an action to each state. Under the policy $\pi$, the FMDP evolves as follows. The initial state $\tilde{s}^0$ is random and satisfies $\tilde{s}^0 = s$, where $s = (s_1, \ldots, s_N) \in \mathcal{S}$, with probability $\prod_{n \in [N]} q_n (s_n)$. If the random state $\tilde{s}^t$ of the FMDP at time $t \in \mathbb{N}_0$ satisfies $\tilde{s}^t = s$, $s \in \mathcal{S}$, then the action $\pi (s)$ is taken and the system transitions to the new state  $\tilde{s}^{t+1} = s'$, $s' = (s'_1, \ldots, s'_N) \in \mathcal{S}$, with probability $\prod_{n \in [N]} p_n (s'_n \, | \, s, \pi (s))$. In other words, both the initial and the transition probabilities of an FMDP multiply across sub-states. For the transition probabilities, this implies that conditional on the current state and the selected action, each sub-state evolves independently. We are interested in policies $\pi^\star$ that maximize the \emph{expected total reward} $\mathbb{E} \left[ \sum_{t \in \mathbb{N}_0} \gamma^t \cdot r (\tilde{s}^t, \pi (\tilde{s}^t)) \, | \, \tilde{s}^0 \sim q \right]$, where the stochastic state evolution $\{ \tilde{s}^t \}_{t \in \mathbb{N}_0}$ is governed by $q$, $p$ and $\pi$ as described above, and where the one-step rewards $r (s, a)$ constitute sums of the reward components $r_j (s, a)$, $j \in [J]$. (The purpose of the reward components will become clear shortly.) It is well-known that there always exists a deterministic and stationary policy that maximizes the expected total reward among all (possibly non-stationary and/or randomized) policies, see, \eg,~\citet[Theorem 6.2.7]{Puterman}.

\begin{ex}[Predictive Maintenance, Cont'd]\label{ex:machine_rep:factored}
    Assume that in Example~\ref{ex:machine_rep:initial}, different machines deteriorate independently of another, conditional on the current state and action. We can model this problem as an FMDP with $10$ sub-state spaces, one for each machine's condition, and $45$ actions, thus requiring the consideration of $10 \cdot 10^{10} \cdot 45 \cdot 10 = 4.5 \cdot 10^{13}$ transition probabilities.
\end{ex}

The representation of an FMDP can be compressed further when the rewards and transition probabilities only depend on parts of the current state. To formalize this idea, we define the \emph{scope} of a function $f : \mathcal{X} \rightarrow \mathbb{R}$ with $\mathcal{X} \subseteq \mathbb{B}^I$ as $\mathfrak{s} [f] = [I] \setminus \{ i \in [I] \, : \, x_{-i} = x'_{-i} \; \Rightarrow \; f(x) = f(x') \;\; \forall x, x' \in \mathcal{X} \}$. Thus, the scope of $f$ comprises the indices of those components of $f$'s input that can impact $f$'s value. We informally say that $f$ has a \emph{low scope} if $| \mathfrak{s} [f] | \ll I$. More formally, a sequence of functions $\{ f^j \}_{j \in \mathbb{N}}$, $f^j : \mathcal{X}^j \rightarrow \mathbb{R}$ with $\mathcal{X}^j \subseteq \mathbb{B}^{I_j}$ and $I_j$ monotonically increasing, has a low scope if $| \mathfrak{s} [f^j] |$ is bounded from above by a logarithm of $I_j$ as $j$ grows large. We also define $\mathfrak{S} [f] = \{ \sum_{i \in [I]} \mathbf{1} [i \in \mathfrak{s} [f] ] \cdot x_i \cdot \mathrm{e}_i \, : \, x \in \mathcal{X} \}$, where $\mathrm{e}_i \in \mathbb{R}^I$, as the set of all sub-vectors with components $i \in \mathfrak{s} [f]$, padded with zeros to form a valid input to $f$. (If $\mathcal{X}$ does not contain some of these padded vectors due to the presence of the zeros, then the domain of $f$ can be extended in the obvious way to accommodate for those inputs.) Note that $| \mathfrak{S} [f] | \ll | \mathcal{X} |$ whenever $f$ has a low scope. FMDPs are particularly powerful when the reward components $r_j$ \mbox{and the transition probabilities $p_n$ exhibit a low scope.}

\begin{ex}[Predictive Maintenance, Cont'd]\label{ex:machine_rep:cpds}
    Assume that in Example~\ref{ex:machine_rep:factored}, the transition probabilities of each machine only depend on the current state of that machine, whether or not that machine is being repaired, as well as the current state of its immediate predecessor in the sequence  (if existent). The resulting FMDP requires the consideration of of $10 \cdot 2 \cdot 10 = 200$ transition probabilities for the first machine as well as $(10 \cdot 10) \cdot 2 \cdot 10 = 2{,}000$ transition probabilities for each subsequent machine, thus reducing the storage requirements to $18{,}200$ transition probabilities.
\end{ex}

Informally, the transition probabilities in Example~\ref{ex:machine_rep:cpds} exhibit a low scope since each $p_n$ depends only on the current state of the machine $n$ as well as at most one predecessor machine, as opposed to the current state of all machines. In view of our formal definition, consider a family of maintenance problems where the number $N$ of machines varies. In that setting, each of the transition probabilities $p_n$ would still only depend on the current state of at most two machines, that is, its scope would remain constant when the size of the problem grows.

So far, we have focused on exploiting a low scope in the transition probabilities. Analogously, we will exploit a low scope in the reward components $r_j (s, a)$ by assuming that each $r_j : \mathcal{S} \times \mathcal{A} \rightarrow \mathbb{R}$ exhibits a low scope in $s \in \mathcal{S}$. In the context of our maintenance problem, for example, it is natural to assume that the rewards in each period constitute a sum of machine-specific rewards $r_j (s, a)$, $j \in [N]$, that are determined by each machine's state $s_j$.

In addition to a low-scope structure, FMDPs can exploit context-specific independence. Section~\ref{sec:FMDP:feature_representation} describes this concept and introduces an MILP-based feature representation that unifies and extends existing approaches for context-specific independence. Afterwards, Section~\ref{sec:FMDP:cutting_plane} develops and analyzes a cutting plane scheme to solve FMDPs with MILP-based feature representations.

\subsection{MILP-Based Feature Representation}\label{sec:FMDP:feature_representation}

So far, our parsimonious representation of MDPs exploits \emph{(i)} a \emph{factored structure}, under which the one-step transitions are characterized by the product $\prod_{n \in [N]} p_n (s'_n \, | \, s, a)$ of individual transition probabilities $p_n$ for each sub-state $n \in [N]$, as well as \emph{(ii)} a \emph{low-scope structure}, under which the transition probabilities $p_n$ and reward components $r_j$ exhibit a low scope with respect to the current state $s \in \mathcal{S}$ of the system. Additionally, many application domains possess \emph{context-specific independence}, under which each individual transition probability $p_n$ and reward component $r_j$ only depend on a small subset of its scope in any particular context \mbox{(\emph{i.e.}, the current state and selected action).}

\begin{ex}[Predictive Maintenance, Cont'd]\label{ex:machine_rep:csi}
    Assume that in Example~\ref{ex:machine_rep:cpds}, the transition probabilities of each machine only depend on the current state of that machine, whether or not that machine is being repaired, as well as whether or not its immediate predecessor in the sequence is in a poor condition (state $1$, $2$ or $3$). We can model this problem as an FMDP with context-specific independence that requires the storage of $10 \cdot 2 \cdot 10 = 200$ transition probabilities for the first machine as well as $(2 \cdot 10) \cdot 2 \cdot 10 = 400$ transition probabilities for each subsequent machine, thus reducing the consideration to $3{,}800$ transition probabilities. Here, the reduction from $2{,}000$ to $400$ transition probabilities for the machines $2, \ldots, 10$ is achieved by distinguishing between the two contexts where the predecessor machine is in  poor condition or not, rather than storing explicitly the transition probabilities for each of the $10$ possible states of the predecessor machine.
\end{ex}

\begin{figure}[!p]
    \includegraphics[width = \textwidth]{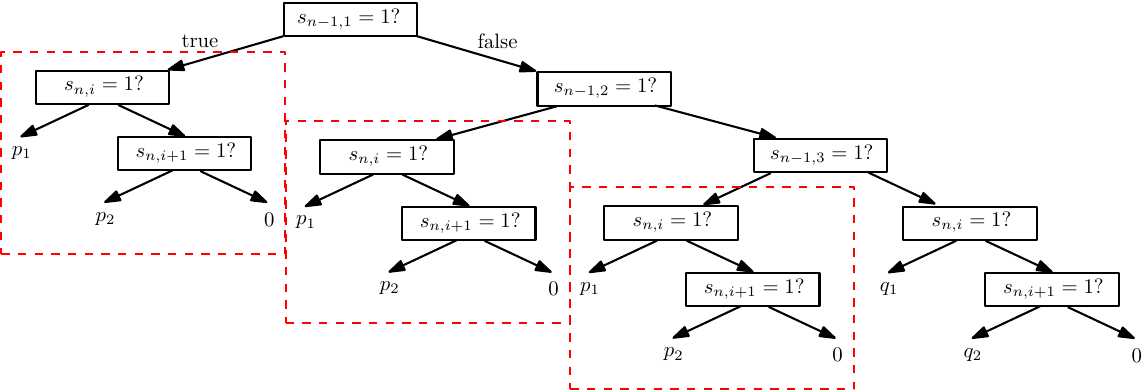}
    \begin{center} \textbf{(a)} Tree representation \end{center}
    \begin{center} \begin{tabular}{l}
        $f^1 (s_{n-1,1}, \ldots, s_{n, i + 1}) \; = \; s_{n-1, 1} \cdot g^1 (s_{n,i}, s_{n,i+1}) + (1 - s_{n-1, 1}) \cdot f^2 (s_{n-1,2}, \ldots, s_{n, i + 1})$ \\
        $f^2 (s_{n-1,2}, \ldots, s_{n, i + 1}) \; = \; s_{n-1, 2} \cdot g^1 (s_{n,i}, s_{n,i+1}) + (1 - s_{n-1, 2}) \cdot f^3 (s_{n-1,3}, \ldots, s_{n, i + 1})$ \\
        $f^3 (s_{n-1,3}, \ldots, s_{n, i + 1}) \; = \; s_{n-1, 3} \cdot g^1 (s_{n,i}, s_{n,i+1}) + (1 - s_{n-1, 3}) \cdot h^1 (s_{n,i}, s_{n,i+1})$ \\
        $g^1 (s_{n,i}, s_{n,i+1}) \; = \; s_{n,i} \cdot p_1 + (1 - s_{n,i}) \cdot g^2 (s_{n,i+1})$ \\
        $g^2 (s_{n,i+1}) \; = \; s_{n,i+1} \cdot p_2 + (1 - s_{n,i+1}) \cdot 0$ \\
        $h^1 (s_{n,i}, s_{n,i+1}) \; = \; s_{n,i} \cdot q_1 + (1 - s_{n,i}) \cdot h^2 (s_{n,i+1})$ \\
        $h^2 (s_{n,i+1}) \; = \; s_{n,i+1} \cdot q_2 + (1 - s_{n,i+1}) \cdot 0$
    \end{tabular} \end{center}
    \begin{center} \textbf{(b)} ADD representation \end{center}
    \begin{center} \begin{tabular}{ll}
        $s_{n-1,1} = 1 \; \wedge \; s_{n,i} = 1 \;\; \Rightarrow \;\; p_1$; & $s_{n-1,1} = 1 \; \wedge \; s_{n,i+1} = 1 \;\; \Rightarrow \;\; p_2$ \\
        $s_{n-1,2} = 1 \; \wedge \; s_{n,i} = 1 \;\; \Rightarrow \;\; p_1$; & $s_{n-1,2} = 1 \; \wedge \; s_{n,i+1} = 1 \;\; \Rightarrow \;\; p_2$ \\
        $s_{n-1,3} = 1 \; \wedge \; s_{n,i} = 1 \;\; \Rightarrow \;\; p_1$; & $s_{n-1,3} = 1 \; \wedge \; s_{n,i+1} = 1 \;\; \Rightarrow \;\; p_2$ \\
        \multicolumn{2}{l}{$s_{n-1,1} = 0 \; \wedge \; s_{n-1,2} = 0 \; \wedge \; s_{n-1,3} = 0 \; \wedge \; s_{n,i} = 1 \;\; \mspace{18mu} \Rightarrow \;\; q_1$} \\
        \multicolumn{2}{l}{$s_{n-1,1} = 0 \; \wedge \; s_{n-1,2} = 0 \; \wedge \; s_{n-1,3} = 0 \; \wedge \; s_{n,i+1} = 1 \;\; \Rightarrow \;\; q_2$} \\
        $s_{n,i} = 0 \; \wedge \; s_{n,i+1} = 0 \;\; \Rightarrow \;\; 0$
    \end{tabular} \end{center}
    \begin{center} \textbf{(c)} Rule-based representation \end{center}
    $\mspace{-25mu}$
    \begin{tabular}{llll}
        $\phi_1 (s, a) \; \geq \; s_{n-1,j} \;\; \forall j \in [3]$;
        & $\phi_1 (s, a) \leq s_{n-1,1} + s_{n-1,2} + s_{n-1,3}$;
        & $\phi_2 (s, a) \; = \; s_{n,i}$;
        & $\phi_3 (s, a) \; = \; s_{n, i+1}$
    \end{tabular}
    \begin{center} \textbf{(d)} MILP-based feature representation \end{center}
    \caption{Different representations of the transition probabilities $p_n (s_{n,i} = 1 \, | \, s, a)$, $n > 1$ and $a_n = 0$, for the FMDP from Example~\ref{ex:machine_rep:csi_reps}. In part~(b), we assume that $p_n (s_{n,i} = 1 \, | \, s, a) = f^1 (s)$. \label{fig:csi}}
\end{figure}

%\clearpage

Traditionally, FMDPs capture context-specific independence via decision trees \citep{Boutellier}, algebraic decision diagrams \citep{Hoey} or rule systems \citep{ZP99:csi}.

A \emph{decision tree} represents a function $f : \mathcal{S} \rightarrow \mathbb{R}$ as a tree that branches upon the values of some of the components of $s \in \mathcal{S}$ and that contains the function values $f (s)$ as leaf nodes. We can exploit context-specific independence in the transition probabilities $p_n (s_n' \, | \, s, a)$ by constructing a separate tree for each sub-state $n \in [N]$, each possible value $s_n' \in \mathcal{S}_n$ of that sub-state and each action $a \in \mathcal{A}$. Likewise, we can exploit context-specific independence in the rewards $r (s, a)$ by constructing a separate tree for each reward component $r_j$, $j \in [J]$, and each action $a \in \mathcal{A}$.

An \emph{algebraic decision diagram} (ADD) $f : \mathbb{B}^I \rightarrow \mathbb{R}$ of binary vectors $x = (x_1, \ldots, x_I)$ is defined recursively as follows: An empty function $f : \emptyset \rightarrow \mathbb{R}$ satisfying $f() = c$, for any $c \in \mathbb{R}$, is an ADD. If two functions $f', f'' : \mathbb{B}^{I - i} \rightarrow \mathbb{R}$ are ADDs of binary vectors $(x_{i+1}, \ldots, x_I)$, then the function $f : \mathbb{B}^{I - i + 1} \rightarrow \mathbb{R}$ satisfying $f (x_i, \ldots, x_I) = x_i \cdot f' (x_{i+1}, \ldots, x_I) + (1 - x_i) \cdot f'' (x_{i+1}, \ldots, x_I)$ is an ADD of binary vectors $(x_i, \ldots, x_I)$. For each sub-state $n \in [N]$ and each possible value $s_n' \in \mathcal{S}_n$ of that sub-state, the associated one-step transition probabilities $p_n (s_n' | s, a)$ can be represented as a collection of ADDs of the state $s$, one for each action $a \in \mathcal{A}$. Likewise,  we can exploit context-specific independence in the rewards $r (s, a)$ by constructing a separate ADD for each reward component $r_j$, $j \in [J]$, and each action $a \in \mathcal{A}$. To facilitate an efficient representation, the order in which the components of $s$ are inspected \mbox{can vary across different transition probabilities and reward components.}

Finally, a \emph{rule system} represents a function $f : \mathcal{S} \rightarrow \mathbb{R}$ via a number of mutually exclusive and jointly exhaustive conjunctions that involve (negations of) some of the components of the state $s$, together with the associated function values $f(s)$ if the respective conjunctions are satisfied. We can exploit context-specific independence in the transition probabilities $p_n (s_n' \, | \, s, a)$ by constructing a separate rule system for each sub-state $n \in [N]$, each possible value $s_n' \in \mathcal{S}_n$ of that sub-state, and each action $a \in \mathcal{A}$. Likewise, we can exploit context-specific independence in the rewards $r (s, a)$ by constructing a separate rule system for each reward component $r_j$, $j \in [J]$, and each action $a \in \mathcal{A}$.

\begin{ex}[Predictive Maintenance, Cont'd]\label{ex:machine_rep:csi_reps}
    Consider again the maintenance problem from Example~\ref{ex:machine_rep:csi}. Representing the problem as an FMDP, we can define the state space as $\mathcal{S} = \bigtimes_{n \in [10]} \mathcal{S}_n$, where each sub-state space $\mathcal{S}_n = \{ s_n \in \mathbb{B}^{10} \, : \, \sum_{i \in [10]} s_{n,i} = 1 \}$ adopts a unary description of the machine's state under which $s_{n,i} = 1$ if and only if machine $n$ is in state $i \in [10]$. Likewise, we can define the action space as $\mathcal{A} = \{ a \in \mathbb{B}^{10} \, : \, \sum_{n \in [10]} a_n = 2 \}$ such that machine $n$ is repaired if and only if $a_n = 1$. Fix any machine $n > 1$ and action $a \in \mathcal{A}$ with $a_n = 0$ and assume that $p_n (s'_{n, i} = 1 \, | \, s, a) = p_1$ if $s_{n-1,j} = 1$ for some $j \in [3]$ and $s_{n,i} = 1$; $= p_2$ if $s_{n-1,j} = 1$ for some $j \in [3]$ and $s_{n,i+1} = 1$; $= q_1$ if $s_{n-1,j} = 0$ for all $j \in [3]$ and $s_{n,i} = 1$; $= q_2$ if $s_{n-1,j} = 1$ for all $j \in [3]$ and $s_{n,i+1} = 1$. Here, $p_i$ and $q_i$, $i = 1, 2$, are constants from the application domain. Figure~\ref{fig:csi}~(a)--(c) provide representations of the transition probabilities $p_n (s'_{n, i} = 1 \, | \, s, a)$ as trees, ADDs and rule systems, respectively.
\end{ex}

Instead of adopting trees, ADDs or rule systems, we propose to exploit context-specific independence via MILP-based feature representations. More precisely, we express the dependence of the reward components $r_j$ and transition probabilities $p_n$ on both $s$ and $a$ via a \emph{feature map} $\phi : \mathcal{S} \times \mathcal{A} \rightarrow \mathbb{B}^F$ that maps state-action pairs $(s, a)$ to binary feature vectors $\phi (s, a) \in \mathbb{B}^F$:
\begin{subequations}\label{eq:phi_p_and_r}
    \begin{equation}\label{eq:p_and_r}
        p_n (\cdot \, | \, s, a) \; = \; \psi_n( \phi (s, a))
        \quad \text{and} \quad
        r_j (s, a) \; = \; \rho_j (\phi (s, a))
        \qquad \forall (s, a) \in \mathcal{S} \times \mathcal{A},
    \end{equation}
    where $\psi_n: \mathbb{B}^F \rightarrow \Delta(\mathcal{S}_n) $, $n \in [N]$ and $\rho_j : \mathbb{B}^F \rightarrow \mathbb{R}$, $j \in [J]$, can be any functions as long as they exhibit low scopes with respect to their inputs $\varphi \in \{ \phi (s, a) \, : \, (s, a) \in \mathcal{S} \times \mathcal{A} \}$. We require the feature map $\phi$ to have a mixed-integer linear representation of the form
    \begin{equation}
        \phi (s, a) = \varphi
        \quad \Longleftrightarrow \quad
        (\varphi, \zeta; s, a) \in \mathcal{F}
        \;\; \text{for some} \;\; \zeta \in \mathbb{R}^{F_l} \times \mathbb{B}^{F_b},
    \end{equation}
\end{subequations}
where $\mathcal{F}$ is a polyhedron described by $F_c$ linear constraints and $\zeta$ is a vector of $F_l$ continuous and $F_b$ binary auxiliary variables. Note that the representation via feature maps does not restrict generality since $\phi$ can be chosen to be the identity map. In practice, however, we are interested in feature maps under which $\psi_n$ and $\rho_j$ exhibit a low scope. For the sake of simplicity and to support intuition, we henceforth commit a slight abuse of notation and write
\begin{equation*}
        p_n (\cdot \, | \, s, a) \; = \; p_n(\cdot\,|\, \phi (s, a))
        \quad \text{and} \quad
        r_j (s, a) \; = \; r_j (\phi (s, a))
        \qquad \forall (s, a) \in \mathcal{S} \times \mathcal{A}
\end{equation*}
instead of $\psi_n$ and $\rho_j$ in \eqref{eq:p_and_r}.

\begin{ex}[Predictive Maintenance, Cont'd]\label{ex:machine_rep:feature_maps}
    Using the state and action spaces of Example~\ref{ex:machine_rep:csi}, we can represent the same transition probabilities $p_n (s'_{n, i} = 1 \, | \, s, a)$ with the MILP-based features $\phi_1 (s, a) = \mathbf{1} [ s_{n-1,j} = 1 \text{ for some } j \in [3] ]$, $\phi_2 (s, a) = \mathbf{1} [ s_{n,i} = 1]$ and $\phi_3 (s, a) = \mathbf{1} [ s_{n,i + 1} = 1 ]$ by setting $p_n (s'_{n, i} = 1 \, | \, s, a) = p_1$ if $\phi_1 (s, a) = 1$ and $\phi_2 (s, a) = 1$; $= p_2$ if $\phi_1 (s, a) = 1$ and $\phi_3 (s, a) = 1$; $= q_1$ if $\phi_1 (s, a) = 0$ and $\phi_2 (s, a) = 1$; $= q_2$ if $\phi_1 (s, a) = 0$ and $\phi_3 (s, a) = 1$; $= 0$ otherwise. The MILP representations of $\phi_1$, $\phi_2$ and $\phi_3$ are presented in Figure~\ref{fig:csi}~(d).
\end{ex}

Compared to trees, ADDs and rule systems, MILP-based feature representations offer various benefits: \emph{(i)} while perhaps subjective, we would argue that they are easier to define and maintain, \emph{(ii)} they can exploit parsimony in both the state and the action space, and \emph{(iii)} they are at least as efficient and sometimes substantially more efficient in characterizing context-specific independence. The second advantage is illustrated in the following example.

\begin{ex}[Predictive Maintenance, Cont'd]\label{ex:machine_rep:actions}
    Consider a variant of our maintenance problem that comprises $N$ machines in total, out of which $m$ should be repaired in each period. In this case, the action space satisfies $| \mathcal{A} | = {N \choose m}$. Since the representations of trees, ADDs and rule systems grow linearly in the number of actions, they grow exponentially in a natural problem description: for $N = 100$ machines and $m = 20$ repairs per period, for example, we have approximately $5.4 \cdot 10^{20}$ actions. In contrast, for each transition probability $p_n (s'_{n,i} = 1 \, | s, a)$ the MILP-based feature representation of Example~\ref{ex:machine_rep:feature_maps} can readily accommodate this setting through the inclusion of a single additional feature $\phi_4 (s, a) = \mathbf{1} [ a_n = 1]$ for each machine $n \in [N]$.
\end{ex}

We next discuss the third advantage of MILP-based feature representations, namely that they are at least as efficient and sometimes significantly more efficient than trees, ADDs and rule systems in characterizing context-specific independence. To see this, consider a variant of our maintenance problem where the transition probabilities of each machine $n > 1$ depend on whether or not at least $50\%$ of its predecessor machines are in states $8$, $9$ or $10$. In this case, any representation of the transition probabilities via trees or rule systems would necessarily grow exponentially in the problem description since every subset of at least half of each machine's predecessors needs to be covered either by a path from the root node to a leaf (in the case of trees) or a probability rule (in the case of rule systems). In contrast, ADDs as well as MILP-based feature representations can handle this variant efficiently, the latter through the inclusion of $n - 1$ features $\phi_n (s, a) = \mathbf{1} [ \sum_{i \in [n]} (s_{i, 8} + s_{i, 9} + s_{i, 10}) \, \geq \, n/2 ]$, $n \in [N - 1]$. More formally, we can establish that \emph{any} of the three representations from the literature is dominated by our MILP-based feature representation in the following sense.

\begin{thm}[Complexity of Representing Context-Specific Independence]\label{thm:complexity_representing_csi}
    ~
    \begin{enumerate}
        \item[\emph{(i)}] For any description of a function $f : \mathbb{B}^I \rightarrow \mathbb{R}$ as a tree, ADD or rule system, there is an equivalent description as an MILP-based feature representation whose size is polynomially bounded in the original description.
        \item[\emph{(ii)}] There are descriptions of functions $f : \mathbb{B}^I \rightarrow \mathbb{R}$ as MILP-based feature representations for which any equivalent description as a tree, ADD or rule system requires a size that is exponential in the original description.
    \end{enumerate}
\end{thm}

Despite their versatility, even MILP-based feature representations can scale unfavorably in the worst case. This is shown by the next result, which provides an asymptotically tight characterization of the complexity of MILP-based feature representations in the context of our FMDPs.

\begin{thm}[Feature Complexity]\label{thm:generic_and_nasty_features}
    Any feature map $\phi : \mathcal{S} \times \mathcal{A} \rightarrow \mathbb{B}^F$ affords a mixed-integer linear representation with at most $F_c = 2F \cdot |\mathcal{S} \times \mathcal{A}|$ constraints and no auxiliary variables, that is, $F_l = F_b = 0$. Furthermore, unless the two complexity classes P and NP coincide, there are features $\phi$ that require $\Omega (|\mathcal{S} \times \mathcal{A}|)$ many continuous auxiliary variables and/or constraints or $\Omega (\log |\mathcal{S} \times \mathcal{A}|)$ binary auxiliary variables.
\end{thm}

Theorem~\ref{thm:generic_and_nasty_features} shows that any feature can be represented with a number of constraints that is polynomially bounded by $| \mathcal{S} \times \mathcal{A} |$ (which itself typically scales exponentially in a natural description of the problem), and that this upper bound (applied to the constraints and/or continuous auxiliary variables) is tight in the worst case. The introduction of binary auxiliary variables can help to reduce this bound \mbox{at the expense of additional computational effort in the subsequent solution schemes.}

While MILP-based feature representations may scale exponentially in the worst case, compact representations often suffice for powerful features that record complex dependencies on $s$ and $a$.

\begin{prop}[Feature Maps]\label{prop:example_feature_mapping}
    ~
    \begin{enumerate}
        \item \emph{Sub-state selection (single).} The feature $\phi (s, a) = \mathbf{1}[s_n \in \mathcal{S}'_n]$, where $\mathcal{S}'_n \subseteq \mathcal{S}_n$, can be represented using $F_l = |\mathcal{S}'_n|$ and $F_b = 0$ variables as well as $F_c = 1 + (1 + S_n) |\mathcal{S}'_n|$ constraints.
        \item \emph{Sub-state selection (universal).} The feature $\phi (s, a) = \mathbf{1}[s_n \in \mathcal{S}'_n \;\; \forall n \in \mathcal{N}]$, where $\mathcal{N}\subseteq [N]$ and $\mathcal{S}'_n \subseteq \mathcal{S}_n$, $n \in \mathcal{N}$, can be represented using $F_l = \sum_{n \in \mathcal{N}} | \mathcal{S}'_n |$ and $F_b = 0$ variables as well as $F_c = 1 + | \mathcal{N} | + \sum_{n \in \mathcal{N}} (1 + S_n) \cdot | \mathcal{S}'_n |$ constraints.
        \item \emph{Sub-state selection (existential).} The feature $\phi (s, a) = \mathbf{1}[\exists n \in \mathcal{N}\, : \, s_n \in \mathcal{S}'_n]$, where $\mathcal{N}\subseteq [N]$ and $\mathcal{S}'_n \subseteq \mathcal{S}_n$, $n \in \mathcal{N}$, can be represented using $F_l = \sum_{n \in \mathcal{N}} | \mathcal{S}_n \setminus \mathcal{S}'_n |$ and $F_b = 0$ variables as well as $F_c = 1 + | \mathcal{N} | + \sum_{n \in \mathcal{N}} (1 + S_n) \cdot | \mathcal{S}_n \setminus \mathcal{S}'_n |$ constraints.
        \item \emph{Sub-state selection (lower bound).} The feature $\phi (s, a) = \mathbf{1}[s_n \in \mathcal{S}'_n \text{ for at least } \nu \text{ different } n \in \mathcal{N}]$, where $\mathcal{N}\subseteq [N]$ and $\mathcal{S}'_n \subseteq \mathcal{S}_n$, $n \in \mathcal{N}$, can be represented using $F_l = \sum_{n \in \mathcal{N}} | \mathcal{S}'_n |$ and $F_b = 1$ variables as well as $F_c = 2 + \sum_{n \in \mathcal{N}} (1 + S_n) \cdot | \mathcal{S}'_n |$ constraints.
        \item \emph{Sub-state selection (upper bound).} The feature $\phi (s, a) = \mathbf{1}[s_n \in \mathcal{S}'_n \text{ for at most } \nu \text{ different } n \in \mathcal{N}]$, where $\mathcal{N}\subseteq [N]$ and $\mathcal{S}'_n \subseteq \mathcal{S}_n$, $n \in \mathcal{N}$, can be represented using $F_l = \sum_{n \in \mathcal{N}} | \mathcal{S}'_n |$ and $F_b = 1$ variables as well as $F_c = 2 + \sum_{n \in \mathcal{N}} (1 + S_n) \cdot | \mathcal{S}'_n |$ constraints.
    \end{enumerate}
\end{prop}

Proposition~\ref{prop:example_feature_mapping} can be readily extended to analogous features for the selected action $a \in \mathcal{A}$. For the sake of brevity, we omit the (largely similar) statement. We close this subsection with a formal definition of the FMDPs that we will study in the remainder of this paper.

\begin{defi}[Featured-Based Factored MDP]\label{def:fb-fmdp}
    A \emph{feature-based factored MDP} is defined by the tuple $(\mathcal{S}, \mathcal{A}, \phi, q, p, r, \gamma)$, where $\mathcal{S}$, $\mathcal{A}$, $q$ and $\gamma$ coincide with their counterparts from Definition~\ref{def:fmdp}, and $\phi$, $p$ and $r$ are defined as in equation~\eqref{eq:phi_p_and_r}.
\end{defi}

We will use the abbreviation FMDP to refer to both FMDPs (\emph{cf.}~Definition~\ref{def:fmdp}) and feature-based FMDPs; the context will make it clear which definition we are referring to.

\subsection{Cutting Plane Solution Scheme}\label{sec:FMDP:cutting_plane}

Exploiting the factored and low-scope structure as well as the context-specific independence that is often present in practical applications allows us to derive compact representations of MDPs. We next develop an iterative algorithm that computes approximately optimal policies for MDPs described by those representations.

Optimal policies for MDPs are typically computed via value iteration, (modified) policy iteration or linear programming (LP), see \S 6.3--6.5 and \S 6.9 of \cite{Puterman}. In this paper, we use the LP approach to determine approximately optimal policies for FMDPs. Our choice is motivated by the observation of \cite{GKPV03:efficient_solution} that the LP approach tends to outperform both value and policy iteration for FMDPs. To this end, consider the following LP, which injects the factored representation of FMDPs into the standard LP model for MDPs:
\begin{equation}\label{opt:fmdp:exact_lp}
    \begin{array}{l@{\quad}l@{\qquad}l}
        \displaystyle \mathop{\text{minimize}}_{v} & \displaystyle \sum_{s \in \mathcal{S}} \left[ \prod_{n \in [N]} q_n (s_n) \right] v (s) \\
        \displaystyle \text{subject to} & \displaystyle v (s) \geq \sum_{j \in [J]} r_j (\phi(s, a)) + \gamma \sum_{s' \in \mathcal{S}} \left[ \prod_{n \in [N]} p_n (s'_n \, | \, \phi(s, a)) \right] v (s') & \displaystyle \forall (s, a) \in \mathcal{S} \times \mathcal{A} \\
        & \displaystyle v : \mathcal{S} \rightarrow \mathbb{R}
     \end{array}
\end{equation}
The decision variable $v$ in this problem represents the value function of the FMDP; at optimality, $v (s)$ records the expected total reward if the system starts in state $s \in \mathcal{S}$ and is operated under an optimal policy. Note that despite the factored structure of FMDPs, the optimal value function does not possess a factored representation in general. The constraints ensure that the value function $v$ weakly overestimates the expected total reward associated with any state $s \in \mathcal{S}$ under any possible policy, while the objective function minimizes the expected total reward under the initial distribution $q$. Whenever $q_n (s_n) > 0$ for all $s \in \mathcal{S}_n$ and $n \in [N]$, problem~\eqref{opt:fmdp:exact_lp} has a unique optimal solution that coincides with the fixed point of the Bellman operator for the FMDP \citep[Theorem~6.2.2 and Section~6.9.1]{Puterman}. Irrespective of whether problem~\eqref{opt:fmdp:exact_lp} has a unique solution, any optimal solution $v^\star$ to the problem allows us to recover an optimal policy $\pi^\star$ for the FMDP via the \emph{greedy policy}
\begin{equation*}
    \pi^\star (s) \in \mathop{\arg \max}_{a \in \mathcal{A}} \left\{ \sum_{j \in [J]} r_j (\phi(s, a)) + \gamma \sum_{s' \in \mathcal{S}} \left[ \prod_{n \in [N]} p_n (s'_n \, | \, \phi(s, a)) \right] v^\star (s') \right\}
    \qquad \forall s \in \mathcal{S},
\end{equation*}
see \citet[Theorem~6.2.7]{Puterman}. Any optimal policy constructed this way attains an expected total reward that coincides with the optimal objective value of problem~\eqref{opt:fmdp:exact_lp}.

Problem~\eqref{opt:fmdp:exact_lp} appears to be difficult to solve as it comprises exponentially many decision variables and constraints as well as sums over exponentially many terms. In fact, we next show that in contrast to ordinary MDPs, which can be solved in polynomial time, problem~\eqref{opt:fmdp:exact_lp} is EXPTIME-hard due to its ability to represent complex problems in a compact way.

\begin{thm}\label{thm:complexity_factored_mdps}
    Computing the optimal expected total reward of an FMDP is EXPTIME-hard even when the feature map $\phi$ is the identity and all transition probabilities and reward components have scopes that are logarithmic in $N$.
\end{thm}

In the following, we simplify problem~\eqref{opt:fmdp:exact_lp} in two steps. We first replace the exponentially-sized value function $v$ with a linearly weighted value function approximation, which will also allow us to eliminate the exponentially-sized sums in~\eqref{opt:fmdp:exact_lp}. After that, we introduce the constraints of problem~\eqref{opt:fmdp:exact_lp} iteratively through a cutting plane scheme.

%\begin{appr}[Value Function Approximation]\label{appr:fmdp:factored_value_function}
%    We approximate the value function $v : \mathcal{S} \rightarrow \mathbb{R}$ in problem~\eqref{opt:fmdp:exact_lp} by the function
%    \begin{equation*}
%        v' (s) = \sum_{k \in [K]} w_k \cdot \nu_k (s),
%    \end{equation*}
%    where $\nu_k : \mathcal{S} \rightarrow \mathbb{R}$ are pre-selected basis functions and $w_k \in \mathbb{R}$ become the new decision variables in problem~\eqref{opt:fmdp:exact_lp}, $k \in [K]$. 
%\end{appr}

In view of the first step, we approximate the value function $v : \mathcal{S} \rightarrow \mathbb{R}$ in problem~\eqref{opt:fmdp:exact_lp} by
\begin{equation*}
    v' (s) = \sum_{k \in [K]} w_k \cdot \nu_k (s),
\end{equation*}
where $\nu_k : \mathcal{S} \rightarrow \mathbb{R}$ are pre-selected linearly independent basis functions and $w_k \in \mathbb{R}$ become the new decision variables in problem~\eqref{opt:fmdp:exact_lp}, $k \in [K]$. We assume that the \emph{block-wise scope} of each basis function $\nu_k$, $\overline{\mathfrak{s}} [\nu_k] = [N] \setminus \{ n \in [N] \, : \, s_{-n} = s'_{-n} \; \Rightarrow \; \nu_k(s) = \nu_k(s) \;\; \forall s, s' \in \mathcal{S} \}$, $k \in [K]$, is low. In contrast to our earlier definition of the scope $\mathfrak{s}$, which considers the individual binary function arguments, the block-wise scope $\overline{\mathfrak{s}}$ considers the sub-states of the input (which consist of $S_1, S_2, \ldots, S_N$ binary components each) as atomic units. In analogy to $\mathfrak{S}$, we also define $\overline{\mathfrak{S}} [\nu_k]$ as the set of all sub-vectors with sub-states $n \in \overline{\mathfrak{s}} [\nu_k]$, padded with zeros to form a valid input to $\nu_k$.

Value function approximations are very common in approximate dynamic programming, see, \emph{e.g.}, \cite{SS85:adp}, \cite{dFVR03:adp} and \S 4.7.2 of \cite{P11:adp}. In contrast to the LP approach employed in that literature, which requires knowledge of the unknown optimal policy to sample constraints efficiently and which results in probabilistic guarantees, we can exploit the factored representation of FMDPs to design a cutting plane solution scheme that only relies on available information about the system and that results in deterministic guarantees.

\begin{prop}\label{prop:value_function_reformulation}
%    Under Approximation~\ref{appr:fmdp:factored_value_function}, problem~\eqref{opt:fmdp:exact_lp} is equivalent to the LP
    Under our value function approximation, problem~\eqref{opt:fmdp:exact_lp} is equivalent to the LP
    \begin{equation}\label{opt:fmdp:basis_fcts_lp_simplified}
        \mspace{-15mu}
        \begin{array}{l@{\quad}l}
            \displaystyle \mathop{\text{\emph{minimize}}}_{w} & \displaystyle \sum_{k \in [K]} w_k  \sum_{s \in \overline{\mathfrak{S}} [\nu_k]}  \nu_k (s) \prod_{n \in \overline{\mathfrak{s}} [\nu_k]} q_n (s_n) \\
            \displaystyle \text{\emph{subject to}} & \displaystyle \sum_{k \in [K]} w_k \cdot \nu_k (s) \geq \sum_{j \in [J]} r_j (\phi (s, a)) + \gamma \sum_{k \in [K]} w_k  \sum_{s' \in \overline{\mathfrak{S}} [\nu_k]} \nu_k (s')  \prod_{n \in \overline{\mathfrak{s}} [\nu_k]} p_n (s'_n \, | \, \phi (s, a)) \\
            & \displaystyle \mspace{500mu} \forall (s, a) \in \mathcal{S} \times \mathcal{A} \\
            & \displaystyle w \in \mathbb{R}^K.
        \end{array}
    \end{equation}
    Moreover, if there is a linear combination $w^\mathbf{1} \in \mathbb{R}^K$ of the basis functions that produces the all-one vector $\mathbf{1} \in \mathbb{R}^{| \mathcal{S} |}$, then problem~\eqref{opt:fmdp:basis_fcts_lp_simplified} attains its (finite) optimal value.
    %Moreover, if $q (s) > 0$ for all $s \in \mathcal{S}$ and the vectors $(\cdots [\nu_k (s)]_{s \in \mathcal{S}} \cdots)^\top \in \mathbb{R}^{| \mathcal{S} |}$, $k \in [K]$, are linearly independent and contain the all-one vector $\mathbf{1} \in \mathbb{R}^{| \mathcal{S} |}$ in their span, then problem~\eqref{opt:fmdp:basis_fcts_lp_simplified} attains its (finite) optimal value.
\end{prop}

Contrary to problem~\eqref{opt:fmdp:exact_lp}, the optimal solution to problem~\eqref{opt:fmdp:basis_fcts_lp_simplified} may not be unique. Any optimal solution $w^\star$ to~\eqref{opt:fmdp:basis_fcts_lp_simplified} allows us to recover an \emph{approximately} optimal policy $\pi^\star$ for the FMDP via
\begin{equation}\label{opt:approx_optimal_policy} 
    \pi^\star (s) \in \mathop{\arg \max}_{a \in \mathcal{A}} \left\{ \sum_{j \in [J]} r_j (\phi (s, a)) + \gamma \sum_{k \in [K]} w^\star_k \sum_{s' \in \overline{\mathfrak{S}} [\nu_k]} \nu_k (s') \prod_{n \in \overline{\mathfrak{s}} [\nu_k]} p_n (s'_n \, | \, \phi (s, a)) \right\}
\end{equation}
for all $s \in \mathcal{S}$. Since the greedy policy $\pi^\star$ is constructed from an approximation of the optimal value function, it will no longer coincide with the optimal policy in general. We can estimate the expected total reward generated by $\pi^\star$ via simulation; due to the suboptimality of $\pi^\star$, this results in a lower bound on the expected total reward generated by the unknown optimal policy. On the other hand, the optimal value of problem~\eqref{opt:fmdp:basis_fcts_lp_simplified} bounds the expected total reward generated by the unknown optimal policy from above since the value function approximation reduces the feasible region in problem~\eqref{opt:fmdp:basis_fcts_lp_simplified}, and thus the optimal value of problem~\eqref{opt:fmdp:basis_fcts_lp_simplified} exceeds that of problem~\eqref{opt:fmdp:exact_lp}. In summary, the solution of problem~\eqref{opt:fmdp:basis_fcts_lp_simplified} provides an implementable greedy policy as well as a bound on the suboptimality of that policy relative to the unknown optimal policy.

Determining $\pi^\star (s)$ according to equation~\eqref{opt:approx_optimal_policy} appears to be difficult as it involves the consideration of exponentially many actions $a \in \mathcal{A}$. We will see below in Corollary~\ref{cor:determine_best_action}, however, that $\pi^\star (s)$ can be evaluated through the solution of an MILP of moderate size. This allows for a compact implicit representation of $\pi^\star$ in terms of the weight vector $w \in \mathbb{R}^K$, as opposed to the exponentially large value functions $v$ in problem~\eqref{opt:fmdp:exact_lp}.

\begin{algorithm}[tb]
\caption{Cutting Plane Scheme to Solve Problem~\eqref{opt:fmdp:basis_fcts_lp_simplified}} \label{alg:CP}
\begin{algorithmic}[1]
\State \textbf{Input:} Initial constraint set $\mathcal{C} \subset \mathcal{S} \times \mathcal{A}$ and an optimality tolerance $\epsilon \geq 0$.
\State Set constraint violation to $+\infty$.
\While{constraint violation exceeds $\epsilon$}
\State \textbf{Master Problem:} Solve the variant of problem~\eqref{opt:fmdp:basis_fcts_lp_simplified} that only contains the constraints pertaining to $(s, a) \in \mathcal{C}$ and record the optimal weight vector $w^\star$.
\State \textbf{Subproblem:} Determine the constraint violation and the associated constraint
    \begin{equation*}
        \mspace{-28mu}
        (s^\star, a^\star) \in \mathop{\arg \max}_{(s, a) \in \mathcal{S} \times \mathcal{A}} \Bigg\{
        \sum_{j \in [J]} r_j (\phi (s, a)) + \gamma \sum_{k \in [K]} w^\star_k \sum_{s' \in \overline{\mathfrak{S}} [\nu_k]} \nu_k (s') \prod_{n \in \overline{\mathfrak{s}} [\nu_k]} p_n (s'_n \, | \, \phi (s, a)) - \sum_{k \in [K]} w^\star_k \cdot \nu_k (s)
        \Bigg\}.
    \end{equation*}
\State If the constraint violation exceeds $\epsilon$, add $(s^\star, a^\star)$ to $\mathcal{C}$.
\EndWhile
\State \textbf{Output:} An approximately optimal solution $w^\epsilon = w^\star$ to problem~\eqref{opt:fmdp:basis_fcts_lp_simplified}.
\end{algorithmic}
\end{algorithm}

For suitable value function approximations, problem~\eqref{opt:fmdp:basis_fcts_lp_simplified} has a small number of decision variables while the number of constraints remains large. To combat this issue, Algorithm~\ref{alg:CP} employs a cutting plane scheme to iteratively introduce only those constraints that are relevant for the solution of the problem. The algorithm assumes that the initial constraint set $\mathcal{C} \subset \mathcal{S} \times \mathcal{A}$ is large enough so that the restricted master problem is not unbounded, which can be readily verified.

The key characteristics of Algorithm~\ref{alg:CP} are summarized in the next result.

\begin{thm}\label{thm:Algo_1}
    Algorithm~\ref{alg:CP} terminates in finite time. Let $F (w)$ and $F^\star$ denote the objective function value of $w \in \mathbb{R}^K$ in problem~\eqref{opt:fmdp:basis_fcts_lp_simplified} as well as the optimal value of~\eqref{opt:fmdp:basis_fcts_lp_simplified}, respectively. Let $w^\epsilon$ be the output of Algorithm~\ref{alg:CP}, and assume that there is a linear combination $w^\mathbf{1} \in \mathbb{R}^K$ of the basis functions that produces the all-one vector $\mathbf{1} \in \mathbb{R}^{| \mathcal{S} |}$. Then $\hat{w} = w^\epsilon + \frac{\epsilon}{1-\gamma} \cdot w^\mathbf{1}$ is feasible in~\eqref{opt:fmdp:basis_fcts_lp_simplified} and satisfies
    \begin{equation*}
        F^\star - \kappa_1 \cdot \epsilon
        \;\; \leq \;\;
        F(w^\epsilon)
        \;\; \leq \;\;
        F^\star
        \;\; \leq \;\;
        F(\hat{w})
        \;\; \leq \;\;
        F^\star + \kappa_2 \cdot \epsilon,
    \end{equation*}
    where $\kappa_1, \kappa_2 > 0$ are independent of $\epsilon$.
\end{thm}

Theorem~\ref{thm:Algo_1} implies that Algorithm~\ref{alg:CP} optimally solves problem~\eqref{opt:fmdp:basis_fcts_lp_simplified} when $\epsilon = 0$. The master problem of Algorithm~\ref{alg:CP} is an LP that can be solved efficiently. Without further structure, the subproblem of Algorithm~\ref{alg:CP} would require a complete enumeration of all state-action pairs $(s, a) \in [\mathcal{S} \times \mathcal{A}] \setminus \mathcal{C}$. The structure inherent in FMDPs allows us to determine $(s^\star, a^\star)$ much more efficiently.

\begin{thm}\label{thm:solution_subproblem}
    A maximally violated constraint $(s^\star, a^\star) \in \mathcal{S} \times \mathcal{A}$ for fixed $w^\star \in \mathbb{R}^K$ in problem~\eqref{opt:fmdp:basis_fcts_lp_simplified} can be extracted from an optimal solution $(s^\star, a^\star, \varphi^\star, \zeta^\star, \eta^\star, \xi^\star, \beta^\star)$ to the MILP
    \begin{equation*}
        \mspace{-75mu}
        \begin{array}{l@{\quad}l@{\quad}l}
            \displaystyle \mathop{\text{\emph{maximize}}}_{\substack{s, \, a, \, \varphi, \, \zeta \\ \eta, \, \xi, \, \beta}} & \multicolumn{2}{l}{\mspace{-8mu} \displaystyle \sum_{j \in [J]} \sum_{f \in \mathfrak{S} [r_j] } r_j(f) \cdot \eta_{jf} + \gamma \sum_{k \in [K]} w^\star_k \sum_{f \in \mathfrak{S} [\overline{\nu}_k] } \overline{\nu}_k(f) \cdot \xi_{kf} - \sum_{k \in [K]} w^\star_k \sum_{s' \in \mathfrak{S} [\nu_k] } \nu_k(s') \cdot \beta_{ks'}} \\
            \text{\emph{subject to}} & \displaystyle \eta_{jf} \in [0,1], \;\; 1 + \sum_{i \in \mathfrak{s} [r_j] } \frac{\varphi_i - f_i}{2 f_i - 1} \leq \eta_{jf} \leq (2 f_l - 1) \varphi_l + 1 - f_l & \displaystyle \forall j \in [J], \; \forall f \in \mathfrak{S} [r_j], \; \forall l \in \mathfrak{s} [r_j] \\
            & \displaystyle \xi_{kf} \in [0, 1], \;\; 1 + \sum_{i \in \mathfrak{s} [\overline{\nu}_k]} \frac{\varphi_i - f_i}{2 f_i - 1} \leq \xi_{kf} \leq (2 f_l - 1) \varphi_l + 1 - f_l & \displaystyle \forall k \in [K], \; \forall f \in \mathfrak{S} [\overline{\nu}_k], \; \forall l \in \mathfrak{s} [\overline{\nu}_k] \\
            & \displaystyle \beta_{ks'} \in [0,1], \;\; 1 + \sum_{i \in \mathfrak{s} [\nu_k]} \frac{[s]_i - [s']_i}{2 [s']_i -1} \leq \beta_{ks'} \leq (2 [s']_l - 1) [s]_l + 1 - [s']_l & \displaystyle \forall k \in [K], \; \forall s' \in \mathfrak{S} [\nu_k], \; \forall l \in \mathfrak{s} [\nu_k] \\
            & \displaystyle (s, a) \in \mathcal{S} \times \mathcal{A}, \;\; (\varphi, \zeta; s, a) \in \mathcal{F}, \;\; \zeta \in \mathbb{R}^{F_l} \times \mathbb{B}^{F_b},
        \end{array}
    \end{equation*}
    where the functions $\overline{\nu}_k: \mathbb{B}^F \rightarrow \mathbb{R}$, $k \in [K]$, are defined as
    \begin{equation*}
        \overline{\nu}_k(\phi (s, a)) = \sum_{s' \in \overline{\mathfrak{S}} [\nu_k]} \nu_k (s') \cdot \prod_{n \in \overline{\mathfrak{s}} [\nu_k]} p_n (s'_n \, | \, \phi (s, a))
    \end{equation*}
    and $[s]_l$ refers to the $l$-th (binary) component of the binary vector $s$, $l = 1, \ldots, \sum_{n \in [N]} S_n$.
\end{thm}

One may wonder whether it is necessary to solve NP-hard MILPs in Algorithm~\ref{alg:CP} or whether maximally violated constraints can be identified in polynomial time. To this end, the next result shows that more efficient solution techniques than the one of Theorem~\ref{thm:solution_subproblem} are unlikely to exist.

\begin{prop}\label{prop:complexity_fmdp_subproblem}
  Identifying a maximally violated constraint in Algorithm~\ref{alg:CP} is strongly NP-complete even if all transition probabilities and reward components have scope $1$.
\end{prop}

We can use a straightforward variation of the MILP in Theorem~\ref{thm:solution_subproblem} to recover the greedy policy $\pi^\star$ characterized by equation~\eqref{opt:approx_optimal_policy}.

\begin{coro}\label{cor:determine_best_action}
    For fixed $s \in \mathcal{S}$ and $w^\star \in \mathbb{R}^K$, let $(a^\star, \varphi^\star, \zeta^\star, \eta^\star, \xi^\star)$ be optimal in the MILP
    \begin{equation*}%\label{eq:subproblem_milp}
        \mspace{-40mu}
        \begin{array}{l@{\quad}l@{\quad}l}
            \displaystyle \mathop{\text{\emph{maximize}}}_{a, \, \varphi, \, \zeta, \, \eta, \, \xi} & \multicolumn{2}{l}{\mspace{-8mu} \displaystyle \sum_{j \in [J]} \sum_{f \in \mathfrak{S} [r_j] } r_j(f) \cdot \eta_{jf} + \gamma \sum_{k \in [K]} w^\star_k \sum_{f \in \mathfrak{S} [\overline{\nu}_k] } \overline{\nu}_k(f) \cdot \xi_{kf}} \\
            \text{\emph{subject to}} & \displaystyle \eta_{jf} \in [0,1], \;\; 1 + \sum_{i \in \mathfrak{s} [r_j] } \frac{\varphi_i - f_i}{2 f_i - 1} \leq \eta_{jf} \leq (2 f_l - 1) \varphi_l + 1 - f_l & \displaystyle \forall j \in [J], \; \forall f \in \mathfrak{S} [r_j], \; \forall l \in \mathfrak{s} [r_j] \\
            & \displaystyle \xi_{kf} \in [0, 1], \;\; 1 + \sum_{i \in \mathfrak{s} [\overline{\nu}_k]} \frac{\varphi_i - f_i}{2 f_i - 1} \leq \xi_{kf} \leq (2 f_l - 1) \varphi_l + 1 - f_l & \displaystyle \forall k \in [K], \; \forall f \in \mathfrak{S} [\overline{\nu}_k], \; \forall l \in \mathfrak{s} [\overline{\nu}_k] \\
            & \displaystyle a \in \mathcal{A}, \;\; (\varphi, \zeta; s, a) \in \mathcal{F}, \;\; \zeta \in \mathbb{R}^{F_l} \times \mathbb{B}^{F_b},
        \end{array}
    \end{equation*}
    where $\overline{\nu}_k$ is defined in Theorem~\ref{thm:solution_subproblem}. Then $a^\star$ is a maximizer for state $s$ in equation~\eqref{opt:approx_optimal_policy}.
\end{coro}

We close this section with some refinements of Algorithm~\ref{alg:CP} that, while not improving upon its theoretical worst-case complexity, help to speed up the solution of typical problem instances.

\begin{rem}[Practical Considerations]\label{rem:practical_considerations}
    Instead of adding a single most violated constraint in each iteration of Algorithm~\ref{alg:CP}, one can add any constraint (possibly multiple ones) whose violation exceeds a given threshold. Such constraints can be identified from a prematurely terminated solution of the MILP from Theorem~\ref{thm:solution_subproblem} or any other heuristic. Moreover, since subsequent instances of the master problem in Algorithm~\ref{alg:CP} differ only through the inclusion of additional constraints, one can warm-start the LP solver from a previously optimal solution. 
\end{rem}

\section{Robust Factored Markov Decision Processes}\label{sec:robust_FMDP}

In practice, the transition kernel $p$ of an FMDP typically differs from the true one-step dynamics of a system for at least two reasons. Firstly, any model constitutes a simplification of reality, and as such its underlying assumptions may not accurately capture second-order effects of the actual system. Secondly, the transition probabilities of an FMDP are often estimated from past observations of the system, and they are thus affected by estimation errors. Following the rich literature on robust MDPs, this section assumes that the true transition kernel is  unknown, and that it is only known to reside in an ambiguity set $\mathcal{P} \subseteq \{ p : \mathbb{B}^F \rightarrow \Delta (\mathcal{S}) \}$ of transition kernels,
\begin{gather*}
    \mathcal{P} = \left\{
        \left[ p : \mathbb{B}^F \rightarrow \Delta (\mathcal{S}) \right] \, : \,
        p (\cdot \, | \, \phi(s, a)) \in \mathcal{P} (\phi(s, a)) \;\; \forall (s, a) \in \mathcal{S} \times \mathcal{A}
    \right\} \\
    \text{with} \quad
    \mathcal{P} (\phi(s, a)) = \left\{
        p (\cdot \, | \, \phi(s, a)) \in \Delta (\mathcal{S}) \, : \, \mathrm{Marg}_n (p (\cdot \, | \, \phi(s, a))) \in \mathcal{P}^n (\phi (s, a)) \;\; \forall n \in [N]
    \right\}
\end{gather*}
for some low-scope set-valued mappings $\mathcal{P}^n : \mathbb{B}^F \rightrightarrows \Delta (\mathcal{S}_n)$. Our ambiguity set $\mathcal{P}$ is an instance of $(s, a)$-rectangular ambiguity sets that characterize the uncertainty underlying each element $p (\cdot | \phi(s, a))$ of the transition kernel separately for every state-action pair $(s, a) \in \mathcal{S} \times \mathcal{A}$. This rectangularity guarantees that the Bellman optimality principle continues to hold for our robust FMDPs \citep{iyengar2005robust, nilim2005robust, wiesemann2013robust}. Our ambiguity set stipulates that the marginals of $p (\cdot | \phi(s, a))$ with regards to the different sub-states $n \in [N]$ live in the marginal ambiguity sets $\mathcal{P}^n (\phi(s, a)) \subseteq \Delta (\mathcal{S}_n)$. To this end, we define the marginalization map $\mathrm{Marg}_{\mathcal{N}} : \Delta (\mathcal{S}) \rightarrow \Delta (\prod_{n \in \mathcal{N}} \mathcal{S}_n)$ for a subset $\mathcal{N} \subseteq [N]$ of sub-states as
\begin{equation*}
    \big[ \mathrm{Marg}_{\mathcal{N}} (p) \big] (\sigma) = \sum_{s \in \mathcal{S}} \mathbf{1} \Big[ \bigwedge_{n \in \mathcal{N}} (s_n = \sigma_n) \Big] \cdot p (s).
\end{equation*}
For each $\sigma \in \prod_{n \in \mathcal{N}} \mathcal{S}_n$, $\big[ \mathrm{Marg}_n (p) \big] (\sigma)$ returns the probability that $\tilde{s}_n = \sigma_n$ for all $n \in \mathcal{N}$ under the distribution $p \in \Delta(\mathcal{S})$. For notational convenience, we abbreviate $\mathrm{Marg}_{\{ n \}}$ by $\mathrm{Marg}_n$ for $n \in [N]$.

One readily verifies that the ambiguity set $\mathcal{P}$ is non-empty if and only if all marginal ambiguity sets $\mathcal{P}^n (\phi (s, a))$, $n \in [N]$ and $(s, a) \in \mathcal{S} \times \mathcal{A}$, are non-empty. Indeed, $\mathcal{P}$ is trivially empty if any of the marginal sets $\mathcal{P}^n (\phi (s, a))$ is empty, and $\mathcal{P}$ contains the transition kernel $p$ satisfying $p (s' | \phi(s, a)) = \prod_{n \in [N]} p_n (s'_n | \phi(s, a))$ if each marginal set $\mathcal{P}^n (\phi (s, a))$ contains a distribution $p_n (\cdot | \phi(s, a))$. We emphasize, however, that the ambiguity set $\mathcal{P}$ is not factored in general, that is, it will typically also contain transition kernels other than product distributions as long as those distributions are compatible with the marginal ambiguity sets.

Throughout this section, we assume that each mapping $\mathcal{P}^n$, $n \in [N]$, is polyhedral-valued. Natural choices for the ambiguity set are norm-balls of the form
\begin{equation*}
    \mathcal{P}^n (\varphi) = \left\{ p_n \in \Delta(\mathcal{S}_n) \, : \, \left \lVert p_n - \hat{p}_n (\cdot \, | \, \varphi) \right \rVert \leq \epsilon (\varphi) \right\}
\end{equation*}
that are centered at some empirical transition kernel $\hat{p}$ (\emph{e.g.}, the maximum likelihood estimate given historical observations) and where $\lVert \cdot \rVert$ is a polyhedral norm such as the $1$-norm, the $\infty$-norm or the intersection of both.

In the spirit of robust optimization, we seek a deterministic and stationary policy $\pi : \mathcal{S} \rightarrow \mathcal{A}$ that performs best in view of the most adversarial transition kernel from within the ambiguity set $\mathcal{P}$. Figuratively speaking, we imagine that the decision maker plays a Stackelberg leader-follower game against an adversarial nature which observes the decision maker's policy and subsequently chooses a transition kernel from within the ambiguity set that minimizes the expected total reward. More concretely, our robust FMDPs evolve as follows. The selection of the initial state is the same as for non-robust FMDPs. If the random state $\tilde{s}^t$ of the robust FMDP at time $t \in \mathbb{N}_0$ satisfies $\tilde{s}^t = s$, $s \in \mathcal{S}$, then the action $\pi (s)$ is taken and the system transitions to the new state $\tilde{s}^{t+1} = s'$, $s' = (s'_1, \ldots, s'_N) \in \mathcal{S}$, with probability $p (s' \, | \, \phi(s, \pi (s)))$ corresponding to \emph{some} $p \in \mathcal{P} (\phi(s, \pi(s)))$. Recall that the ambiguity set $\mathcal{P}$ is not factored, and hence $p (s' \, | \, \phi(s, \pi (s))) \neq \prod_{n \in [N]} p_n (s'_n \, | \, \phi(s, \pi (s)))$ in general. We are interested in policies $\pi^\star$ that maximize the \emph{worst-case expected total reward} $\inf_{p \in \mathcal{P}} \, \mathbb{E} \left[ \sum_{t \in \mathbb{N}_0} \gamma^t \cdot r (\tilde{s}^t, \pi (\tilde{s}^t)) \, | \, \tilde{s}^0 \sim q \right]$. Since the ambiguity set $\mathcal{P}$ is $(s, a)$-rectangular by construction, it is well-known that there always exists a deterministic and stationary policy that maximizes the worst-case expected total reward among all (possibly non-stationary and/or randomized) policies, \mbox{see, \eg, \citet{iyengar2005robust}, \citet{nilim2005robust} and \citet{wiesemann2013robust}.}

In the following, we discuss how our cutting plane solution scheme from Section~\ref{sec:FMDP:cutting_plane} extends to the robust setting (Section~\ref{sec:robust_FMDP:cutting_plane}), and how the revised master problem and subproblems can be solved efficiently (Sections~\ref{sec:robust_FMDP:master_problem} and~\ref{sec:robust_FMDP:sub_problem}).

\subsection{Column-and-Constraint Generation Scheme}\label{sec:robust_FMDP:cutting_plane}

Recall from Section~\ref{sec:FMDP:cutting_plane} that the LP~\eqref{opt:fmdp:basis_fcts_lp_simplified} presented in Proposition~\ref{prop:value_function_reformulation} solves an FMDP with our value function approximation if the transition kernel $p$ is known exactly. Standard results from the literature on robust MDPs \citep{iyengar2005robust, nilim2005robust, wiesemann2013robust} imply that a robust FMDP with value function approximation is solved by the bilinear program
\begin{equation}\label{opt:rfmdp:basis_fcts_lp_simplified}
    \begin{array}{l@{\quad}l}
        \displaystyle \mathop{\text{minimize}}_{w, p} & \displaystyle \sum_{k \in [K]} w_k \sum_{s \in \overline{\mathfrak{S}} [\nu_k]} \nu_k (s) \prod_{n \in \overline{\mathfrak{s}} [\nu_k]} q_n (s_n) \\
        \displaystyle \text{subject to} & \displaystyle \sum_{k \in [K]} w_k \cdot \nu_k (s) \geq \sum_{j \in [J]} r_j (\phi (s, a)) + \gamma \sum_{k \in [K]} w_k \sum_{s' \in \overline{\mathfrak{S}} [\nu_k]} \nu_k (s') \cdot p_{sa} (s') \\
        & \displaystyle \mspace{395mu} \forall (s, a) \in \mathcal{S} \times \mathcal{A} \\
        & \displaystyle p_{sa} \in \mathcal{P} (\phi(s, a)), \, (s, a) \in \mathcal{S} \times \mathcal{A},
     \end{array}
\end{equation}
where the adversarial nature optimizes over both the basis function weights $w$ and the transition kernel $p$ from within the ambiguity set $\mathcal{P}$. Problem~\eqref{opt:rfmdp:basis_fcts_lp_simplified} suffers from three challenges: \emph{(i)} it contains exponentially many constraints and decision variables $p_{sa}$; \emph{(ii)} it contains bilinear terms involving the basis function weights $w_k$ and transition probabilities $p_{sa} (s')$; and \emph{(iii)} the decision vectors $p_{sa} \in \Delta (\mathcal{S})$ have exponential length. As for \emph{(i)}, we will employ a column-and-constraint generation scheme that iteratively introduces violated constraints and their associated transition probabilities. In view of \emph{(ii)}, we will relax the objective of solving problem~\eqref{opt:rfmdp:basis_fcts_lp_simplified} exactly to computing a stationary point of the problem. As for \emph{(iii)}, we next impose an additional assumption on the structure of our basis functions $\{ \nu_k \}_{k \in [K]}$ that enables us to optimize over low-scope transition probabilities $p_{sa}$.

\begin{defi}[Running Intersection Property]
    The basis functions $\nu_k$, $k \in [K]$, satisfy the \emph{running intersection property} (RIP) if they can be reordered such that for each $k \in [K] \setminus \{ 1 \}$ there is $k' \in [k-1]$ with $\overline{\mathfrak{s}} [\nu_k] \cap \left( \bigcup_{\kappa \in [k-1]} \overline{\mathfrak{s}} [\nu_\kappa] \right) \subseteq \overline{\mathfrak{s}} [\nu_{k'}]$.
\end{defi}

\begin{obs}\label{obs:rip}
    For $(s, a) \in \mathcal{S} \times \mathcal{A}$ and fixed weights $w \in \mathbb{R}^K$, the set of all expressions
    \begin{equation*}
        \left\{ \sum_{k \in [K]} w_k \sum_{s' \in \overline{\mathfrak{S}} [\nu_k]} \nu_k (s') \cdot p_{sa} (s') \, : \, p_{sa} \in \mathcal{P} (\phi (s, a)) \right\}
    \end{equation*}
    is a subset of the set of all expressions
    \begin{equation*}
        \mspace{-25mu}
        \left\{ \sum_{k \in [K]} w_k \sum_{s' \in \overline{\mathfrak{S}} [\nu_k]} \nu_k (s') \cdot p_{sak} (s') \, : \,
        \left[ \begin{array}{l@{\quad}l}
            \displaystyle \mathrm{Marg}_n (p_{sak}) \in \mathcal{P}^n (\phi (s, a)) & \displaystyle \forall k \in [K], \; \forall n \in \overline{\mathfrak{s}} [\nu_k] \\
            \displaystyle \mathrm{Marg}_{\mathcal{N} (k, k')} (p_{sak}) = \mathrm{Marg}_{\mathcal{N} (k, k')} (p_{sak'}) & \displaystyle \forall 1 \leq k < k' \leq K \\
            \displaystyle p_{sak} \in \Delta (\overline{\mathfrak{S}} [\nu_k]), \; k \in [K]
        \end{array} \right]
        \right\},
    \end{equation*}
    where we use the abbreviation $\mathcal{N} (k, k') = \overline{\mathfrak{s}} [\nu_k] \cap \overline{\mathfrak{s}} [\nu_{k'}]$. Moreover, if the basis functions $\nu_k$, $k \in [K]$, satisfy the RIP, then both sets coincide.
\end{obs}

Observation~\ref{obs:rip} allows us to construct a progressive approximation or equivalent reformulation of problem~\eqref{opt:rfmdp:basis_fcts_lp_simplified} where the exponential-size transition probabilities $p_{sa} \in \Delta (\mathcal{S})$ are replaced with polynomial-size probabilities $p_{sak} \in \Delta (\overline{\mathfrak{S}} [\nu_k])$, $k \in [K]$. The quadratic number of marginal consistency conditions in Observation~\ref{obs:rip} can be reduced to a number that is linear in $K$ (\emph{cf.~}Lemma~1 of \citealt{doan2015robustness}). Note that any progressive approximation of nature's problem~\eqref{opt:rfmdp:basis_fcts_lp_simplified} amounts to a conservative approximation from the decision maker's perspective.

\begin{ex}[Sliding Window Basis Functions]\label{ex:low-scope-basis}
    Consider a set of basis functions $\{ \nu_k \}_k$ such that the block-wise scope of each function $\nu_k$ is a subset of the sliding window $\{ k, k + 1, \ldots, k + m \} \cap [N]$, where $m \in \mathbb{N}$. In that case, $\{ \nu_k \}_k$ satisfy the RIP since
    \begin{align*}
        \overline{\mathfrak{s}} [\nu_k] \cap \left( \bigcup_{\kappa \in [k-1]} \overline{\mathfrak{s}} [\nu_\kappa] \right)
        \;\; &\subseteq \;\;
        \{ k, k + 1, \ldots, k + m \} \cap \left( \bigcup_{\kappa \in [k-1]} \{ \kappa, \kappa + 1, \ldots, \kappa + m \} \right) \cap [N] \\
        &= \;\;
        \{ k, k + 1, \ldots, k + m \} \cap [ k + m - 1] \cap [N]\\
        &= \;\;
        \{ k, k + 1, \ldots, k + m - 1 \} \cap [N]
        \;\; \subseteq \;\;
        \overline{\mathfrak{s}} [\nu_{k-1}].
    \end{align*}
\end{ex}

The use of the running intersection property in MDPs is not entirely new. \citet{GG02:distributed_planning} use this property to enforce consistency between the states of local agents in a collaborative multi-agent system where each agent is represented by an MDP. \citet{LP02:learning_markov_games} apply the running intersection property in the context of team Markov games to ensure that the mixed strategies of different agents in a team can be combined to a mixed team strategy. \citet{CCDYWX15:efficient_alp} use the running intersection property to approximate the FMDP constraints generated by non-serial dynamic programming. In contrast, we employ the running intersection property to obtain equivalent polynomial-sized representations of ambiguity sets in FMDPs.

We are now ready to present our column-and-constraint generation scheme in Algorithm~\ref{alg:CP-robust}. The algorithm assumes that the initial constraint set $\mathcal{C} \subset \mathcal{S} \times \mathcal{A}$ is sufficiently large so that the optimal value of the restricted master problem is not unbounded from below. This can be verified efficiently through the solution of an LP if we fix the transition probabilities $p_{sa}$ to any point in the relative interior of $\mathcal{P} (\phi (s, a))$, $(s, a) \in \mathcal{C}$.

The key characteristics of  Algorithm~\ref{alg:CP-robust} are summarized next.

\begin{algorithm}[tb]
\caption{Column-and-Constraint Generation Scheme to Solve Problem~\eqref{opt:rfmdp:basis_fcts_lp_simplified}} \label{alg:CP-robust}
\begin{algorithmic}[1]
\State \textbf{Input:} Initial state-action set $\mathcal{C} \subset \mathcal{S} \times \mathcal{A}$ and an optimality tolerance $\epsilon \geq 0$.
\State Set constraint violation to $+\infty$.
\While{constraint violation exceeds $\epsilon$}
\State \textbf{Master Problem:} Solve, to a stationary point, the variant of problem~\eqref{opt:rfmdp:basis_fcts_lp_simplified} that only contains the constraints and decision variables $p_{sa}$ pertaining to $(s, a) \in \mathcal{C}$, and then record $w^\star$.
\State \textbf{Subproblem:} Determine the constraint violation and the associated triple
    \begin{equation*}
        \mspace{-35mu}
        (s^\star, a^\star, p^\star) \in \mathop{\arg \max}_{(s, a) \in \mathcal{S} \times \mathcal{A}} \Bigg\{ \sum_{j \in [J]} r_j (\phi (s, a)) + \gamma \min_{p_{sa} \in \mathcal{P} (\phi (s, a))} \Bigg[ \sum_{k \in [K]} w^\star_k \sum_{s' \in \overline{\mathfrak{S}} [\nu_k]} \nu_k (s') \cdot p_{sa} (s') \Bigg] - \sum_{k \in [K]} w^\star_k \cdot \nu_k (s) \Bigg\}.
    \end{equation*}
\State If the constraint violation exceeds $\epsilon$, add $(s^\star, a^\star)$ to $\mathcal{C}$.
\EndWhile
\State \textbf{Output:} An approximate stationary point $(w^\epsilon, p^\epsilon) = (w^\star, p^\star)$ to problem~\eqref{opt:rfmdp:basis_fcts_lp_simplified}.
\end{algorithmic}
\end{algorithm}

\begin{thm}\label{thm:lambda-epsilon-kkt-point}
    %Suppose that problem~\eqref{opt:rfmdp:basis_fcts_lp_simplified} has a finite optimal value.
    Algorithm~\ref{alg:CP-robust} terminates in finite time. Let $(w^\epsilon, p^\epsilon)$ be the output of Algorithm~\ref{alg:CP-robust}, and assume that there is a linear combination $w^\mathbf{1} \in \mathbb{R}^K$ of the basis functions that produces the all-one vector. Then $(\hat{w}, \hat{p})$ with $\hat{w} = w^\epsilon + \frac{\epsilon}{1-\gamma} \cdot w^\mathbf{1}$, and $\hat{p}_{sa} = p^\epsilon_{sa}$ for $(s, a) \in \mathcal{C}$ and $\hat{p}_{sa} \in \arg\min \big\{ \sum_{k \in [K]} w^\epsilon_k \sum_{s' \in \overline{\mathfrak{S}} [\nu_k]} \nu_k (s') \cdot p_{sa} (s') \, : \, p_{sa}\in \mathcal{P} (\phi (s, a)) \big\}$  otherwise is feasible in~\eqref{opt:rfmdp:basis_fcts_lp_simplified}. Moreover, if the RFMDP is irreducible for some $p \in \mathcal{P}$, then the Euclidean distance of $(\hat{w}, \hat{p})$ to the stationary points of~\eqref{opt:rfmdp:basis_fcts_lp_simplified} is bounded from above by $\kappa \cdot \epsilon^\tau + \frac{\epsilon}{1-\gamma} \|w^\mathbf{1}\|_2$, where $\kappa>0$ and $\tau\in (0,1)$ are independent of $\epsilon$. 
\end{thm}

The output $(w^\epsilon, p^\epsilon)$ of Algorithm~\ref{alg:CP-robust} overestimates the worst-case expected total reward achievable by the unknown optimal policy for two reasons. Firstly, and similar to our cutting plane scheme from Section~\ref{sec:FMDP:cutting_plane}, Algorithm~\ref{alg:CP-robust} restricts nature's flexibility by imposing a value function approximation. Secondly, and differently to Section~\ref{sec:FMDP:cutting_plane}, Algorithm~\ref{alg:CP-robust} solves problem~\eqref{opt:rfmdp:basis_fcts_lp_simplified} to stationarity as opposed to global optimalilty. Theorem~\ref{thm:lambda-epsilon-kkt-point} guarantees, however, that the output of Algorithm~\ref{alg:CP-robust} can be completed to a feasible solution $(\hat{w}, \hat{p})$ to~\eqref{opt:rfmdp:basis_fcts_lp_simplified} that is an approximate stationary point, and that $\hat{p}$ represents a global worst case for the value function approximation $\sum_{k \in [K]} \hat{w}_k \cdot \nu_k$. Thus, the output of Algorithm~\ref{alg:CP-robust} mirrors that of Algorithm~\ref{alg:CP} in the sense that both algorithms overestimate the (worst-case) expected total reward of the unknown optimal policy, but both algorithms correctly compute the (worst-case) expected total reward for fixed basis function weights $\hat{w}$.

In analogy to Section~\ref{sec:FMDP:cutting_plane}, an approximate stationary point $(\hat{w}, \hat{p})$ to problem~\eqref{opt:rfmdp:basis_fcts_lp_simplified} allows us to recover an approximately optimal policy $\pi^\star$ for the robust FMDP via
\begin{equation}\label{opt:approx_robust_optimal_policy}
    \pi^\star (s) \in \mathop{\arg \max}_{a \in \mathcal{A}} \left\{ \sum_{j \in [J]} r_j (\phi (s, a)) + \gamma \min_{p_{sa} \in \mathcal{P} (\phi (s, a))} \Bigg[ \sum_{k\in[K]} \hat{w}_k \sum_{s' \in \overline{\mathfrak{S}} [\nu_k]} \nu_k (s') \cdot p_{sa} (s') \Bigg] \right\}
\end{equation}
for all $s \in \mathcal{S}$. The robust greedy policy $\pi^\star$ will not coincide with the optimal robust policy in general since it is constructed from an approximation of the optimal value function. We will see in Section~\ref{sec:num_experiments}, however, that the robust greedy policy~\eqref{opt:approx_robust_optimal_policy} can outperform the nominal greedy policy~\eqref{opt:approx_optimal_policy} of Section~\ref{sec:FMDP:cutting_plane} in data-driven experiments where the true transition kernel $p$ governing the dynamics of an FMDP is unknown and needs to be estimated from historical observations. We defer the efficient computation of the max-min problem in~\eqref{opt:approx_robust_optimal_policy} to Section~\ref{sec:robust_FMDP:sub_problem}.

\subsection{Solution of the Master Problem}\label{sec:robust_FMDP:master_problem}

Recall that the master problem~\eqref{opt:rfmdp:basis_fcts_lp_simplified}, even when restricted to a small subset $\mathcal{C} \subseteq \mathcal{S} \times \mathcal{A}$ of constraints and decision variables $p_{sa}$ with low scope, appears computationally challenging due the presence of bilinearities between the basis function weights $w$ and the transition probabilities $p_{sa}$. In fact, the problem bears similarity to robust MDPs, for which currently no polynomial-time solution scheme is known. For this reason, Algorithm~\ref{eq:update} relaxes the ambitious goal of solving~\eqref{opt:rfmdp:basis_fcts_lp_simplified} to global optimality and instead computes a stationary point of~\eqref{opt:rfmdp:basis_fcts_lp_simplified}. 

%In fact, a convex reformulation has been derived under specific assumptions (Marek's paper), but the resulting optimization problem contains exponentially large (small?) coefficients. Likewise, variants of robust value and policy iteration have been derived for traditional (non-factored) robust MDPs, but their runtimes are not polynomial in the discount factor $\gamma$.

\begin{algorithm}[tb]
\caption{Alternating Minimization Scheme for Restricted Master Problem~\eqref{opt:rfmdp:basis_fcts_lp_simplified}} \label{eq:update}
\begin{algorithmic}[1]
\State \textbf{Input:} Initial solution $(w^0, p^0)$, and regularizer sequence $\{ \lambda_t \}_{t \in \mathbb{N}}$.
\For{$t = 1, 2, \dots$}
\State Let $w^t$ be an optimal solution to the LP-representable problem
    \begin{equation*}
        \begin{array}{l@{\quad}l}
            \displaystyle \mathop{\text{minimize}}_{w} & \displaystyle \sum_{k\in [K]} w_k \sum_{s \in \overline{\mathfrak{S}} [\nu_k]} \nu_k (s) \prod_{n \in \overline{\mathfrak{s}} [\nu_k]} q_n (s_n) + \lambda_t \left \lVert w - w^{t-1} \right \rVert_1 \\
            \displaystyle \text{subject to} & \displaystyle \sum_{k \in [K]} w_k \cdot \nu_k (s) \geq \sum_{j \in [J]} r_j (\phi (s, a)) + \gamma \sum_{k \in [K]} w_k \sum_{s' \in \overline{\mathfrak{S}} [\nu_k]} \nu_k (s') \cdot p^{t-1}_{sak} (s') \\
            & \displaystyle \mspace{450mu} \forall (s, a) \in \mathcal{C} \\
            & \displaystyle w \in \mathbb{R}^K.
        \end{array}
    \end{equation*}
\State For each state-action pair $(s, a) \in \mathcal{C}$, let $p^t$ be an optimal solution to the LP
    \begin{equation*}
        \begin{array}{l@{\quad}l@{\qquad}l}
            \displaystyle \mathop{\text{minimize}}_{p} & \displaystyle \sum_{k\in[K]} \sum_{s' \in \overline{\mathfrak{S}} [\nu_k]} w^t_k \cdot \nu_k (s') \cdot p_{sak} (s') \\
            \displaystyle \text{subject to} & \displaystyle \mathrm{Marg}_n (p_{sak}) \in \mathcal{P}^n (\phi (s, a)) & \displaystyle \forall k \in [K], \; \forall n \in \overline{\mathfrak{s}} [\nu_k] \\
            & \displaystyle \mathrm{Marg}_{\mathcal{N} (k, k')} (p_{sak}) = \mathrm{Marg}_{\mathcal{N} (k, k')} (p_{sak'}) & \displaystyle \forall 1 \leq k < k' \leq K \\
            & \displaystyle p_{sak} \in \Delta (\overline{\mathfrak{S}} [\nu_k]), \; k \in [K].
        \end{array}
    \end{equation*}
\EndFor
\State \textbf{Output:} Stationary point $(w^\infty, p^\infty)$ to the restricted master problem~\eqref{opt:rfmdp:basis_fcts_lp_simplified}.
\end{algorithmic}
\end{algorithm}

\begin{thm}\label{thm:kkt-convergence}
    %Suppose that problem~\eqref{opt:rfmdp:basis_fcts_lp_simplified} has a finite optimal value.
    For any non-increasing regularizer sequence $\{ \lambda_t \}_{t \in \mathbb{N}}$ satisfying $\lambda_t > 0$, $t \in \mathbb{N}$, and $\lambda_t \underset{t \rightarrow \infty}{\longrightarrow} 0$, the sequence $\{ (w^t, p^t) \}_{t \in \mathbb{N}}$ produced by Algorithm~\ref{eq:update} is bounded, and any limit point $(w^\infty, p^\infty)$ is a stationary point to the restricted master problem~\eqref{opt:rfmdp:basis_fcts_lp_simplified}.
\end{thm}

Common choices for the regularizer sequence $\{ \lambda_t \}_{t \in \mathbb{N}}$ include $\lambda_t = e^{-t}$, $\lambda_t = 1/t$ and $\lambda_t = 1/t^2$.

\subsection{Solution of the Subproblems}\label{sec:robust_FMDP:sub_problem}

In analogy to Theorem~\ref{thm:solution_subproblem} from Section~\ref{sec:FMDP:cutting_plane}, we can identify a maximally violated constraint in Algorithm~\ref{alg:CP-robust} via the solution of an MILP.

\begin{coro}\label{coro:robust:solution_subproblem}
    A maximally violated constraint $(s^\star, a^\star) \in \mathcal{S} \times \mathcal{A}$ for fixed $w^\star \in \mathbb{R}^K$ in problem~\eqref{opt:rfmdp:basis_fcts_lp_simplified}, together with the corresponding transition probabilities $p_{s^\star a^\star}$, can be extracted from an optimal solution $(s^\star, a^\star, \varphi^\star, \zeta^\star, \eta^\star, \xi^\star, \beta^\star, \chi^\star, \delta^\star)$ to the MILP
    \begin{equation*}%\label{eq:subproblem_milp}
        \mspace{-95mu}
        \begin{array}{l@{\quad}l@{\quad}l}
            \displaystyle \mathop{\text{\emph{maximize}}}_{\substack{s, \, a, \, \varphi, \, \zeta \\ \eta, \, \xi, \, \beta, \, \chi, \, \delta}} & \multicolumn{2}{l}{\mspace{-8mu} \displaystyle \sum_{j \in [J]} \sum_{f \in \mathfrak{S} [r_j]} r_j(f) \cdot \eta_{jf} - \gamma \sum_{k \in [K]} \sum_{n \in \overline{\mathfrak{s}} [\nu_k]} \sum_{f \in \mathfrak{S} [\mathcal{P}^n]} \psi_n (f)^\top \delta_{knf} - \sum_{k \in [K]} w^\star_k \sum_{s' \in \mathfrak{S} [\nu_k]} \nu_k(s') \cdot \beta_{ks'}} \\
            \text{\emph{subject to}}  & \displaystyle w^\star_k \mu_k + \sum_{n \in \overline{\mathfrak{s}} [\nu_k]} \sum_{f \in \mathfrak{S} [\mathcal{P}^n]} \mathrm{M}_{k, n}{}^\top \Psi_n (f)^\top \delta_{knf}  \geq \\
            & \displaystyle \mspace{60mu} \sum_{1 \leq k' < k} \mathrm{M}_{k, \mathcal{N} (k, k')}{}^\top \chi_{k'k} - \sum_{k < k' \leq K} \mathrm{M}_{k, \mathcal{N} (k, k')}{}^\top \chi_{kk'} & \displaystyle \forall k \in [K] \\
            & \displaystyle \eta_{jf} \in [0,1], \;\; 1 + \sum_{i \in \mathfrak{s} [r_j]} \frac{\varphi_i - f_i}{2 f_i - 1} \leq \eta_{jf} \leq (2 f_l - 1) \varphi_l + 1 - f_l & \displaystyle \forall j \in [J], \; \forall f \in \mathfrak{S} [r_j], \; \forall l \in \mathfrak{s}[r_j] \\
            & \displaystyle \xi_{nf} \in [0, 1], \;\; 1 + \sum_{i \in \mathfrak{s} [\mathcal{P}^n]} \frac{\varphi_i - f_i}{2 f_i - 1} \leq \xi_{nf} \leq (2 f_l - 1) \varphi_l + 1 - f_l & \displaystyle \forall n \in [N], \; \forall f \in \mathfrak{S} [\mathcal{P}^n], \; \forall l \in \mathfrak{s} [\mathcal{P}^n] \\
            & \displaystyle \beta_{ks'} \in [0,1], \;\; 1 + \sum_{i \in \mathfrak{s} [\nu_k]} \frac{[s]_i - [s']_i}{2 [s']_i -1} \leq \beta_{ks'} \leq (2 [s']_l - 1) [s]_l + 1 - [s']_l & \displaystyle \forall k \in [K], \; \forall s' \in \mathfrak{S} [\nu_k], \; \forall l \in \mathfrak{s} [\nu_k] \\
            & \multicolumn{2}{l}{\mspace{-8mu} \displaystyle \chi_{kk'} \in \mathbb{R}^{M_{\mathcal{N} (k, k')}}, \, 1 \leq k < k' \leq K, \;\; \delta_{knf} \in \mathbb{R}^{M_{sa}}_+, \, \delta_{knf \ell} \leq \mathrm{M} \cdot \xi_{nf}, \, k \in [K], \, n \in \overline{\mathfrak{s}} [\nu_k], f \in \mathfrak{S} [\mathcal{P}^n] \text{\emph{ and }} \ell \in [M_{sa}]} \\
            & \displaystyle (s, a) \in \mathcal{S} \times \mathcal{A}, \;\; (\varphi, \zeta; s, a) \in \mathcal{F}, \;\; \zeta \in \mathbb{R}^{F_l} \times \mathbb{B}^{F_b},
        \end{array}
    \end{equation*}
    where $\psi_n$, $\Psi_n$, $\mu_k$, $\mathrm{M}$, $\mathrm{M}_{k,n}$, $\mathrm{M}_{k, \mathcal{N}(k, k')}$, $M_{\mathcal{N} (k, k')}$ and $M_{sa}$ are defined in the proof.
\end{coro}

In analogy to Corollary~\ref{cor:determine_best_action}, we recover the greedy policy $\pi^\star$ characterized by equation~\eqref{opt:approx_robust_optimal_policy} by solving a variant of the previous MILP where we fix the state $s \in \mathcal{S}$, disregard the last expression from the objective function and remove the auxiliary variables $\{ \beta_{ks'} \}$ and their associated constraints.

\section{Factored Ambiguity Sets}\label{sec:robust_FMDP:factored}

The ambiguity sets from Section~\ref{sec:robust_FMDP} can be overly conservative since they allow for arbitrary stochastic dependencies between the transitions of different sub-states $\tilde{s}_n$, $n \in [N]$. This implies, for example, that in data-driven settings where the marginal ambiguity sets $\mathcal{P}^n (\phi(s, a))$, $(s,a) \in \mathcal{S} \times \mathcal{A}$, are estimated from data, the ambiguity sets $\mathcal{P} (\phi (s, a))$ will typically not converge to singleton sets even if all marginal ambiguity sets do. We obtain a less conservative class of ambiguity sets if we impose the additional restriction that the sub-state transitions are independent:
\begin{align*}
    \mathcal{P} (\phi (s, a)) = \Big\{ p (\cdot | \phi (s, a)) \in \Delta(\mathcal{S}) \, : \,
    & \Big[ p (s' | \phi (s, a)) = \prod_{n \in [N]} p_n (s'_n | \phi (s, a)) \;\; \forall s' \in \mathcal{S} \Big] \\
    & \mspace{74mu} \text{for some } p_n \in \mathcal{P}^n (\phi (s, a)), \, n \in [N] \Big\}.
\end{align*}
We call the above ambiguity sets \emph{factored}, as opposed to the \emph{non-factored} ambiguity sets from the previous section. We thus distinguish between RFMDPs with factored ambiguity sets (or factored RFMDPs for short) and RFMDPs with non-factored ambiguity sets (or non-factored RFMDPs for short). For the same marginal ambiguity sets $\mathcal{P}^n$, a factored ambiguity set $\mathcal{P}$ is a subset of its non-factored counterpart.

\begin{figure}[tb]
    $\mspace{-110mu}$
    \begin{tabular}{c@{$\mspace{-5mu}$}c@{$\mspace{-5mu}$}c@{$\mspace{-5mu}$}c}
        \includegraphics[height = 3.8cm]{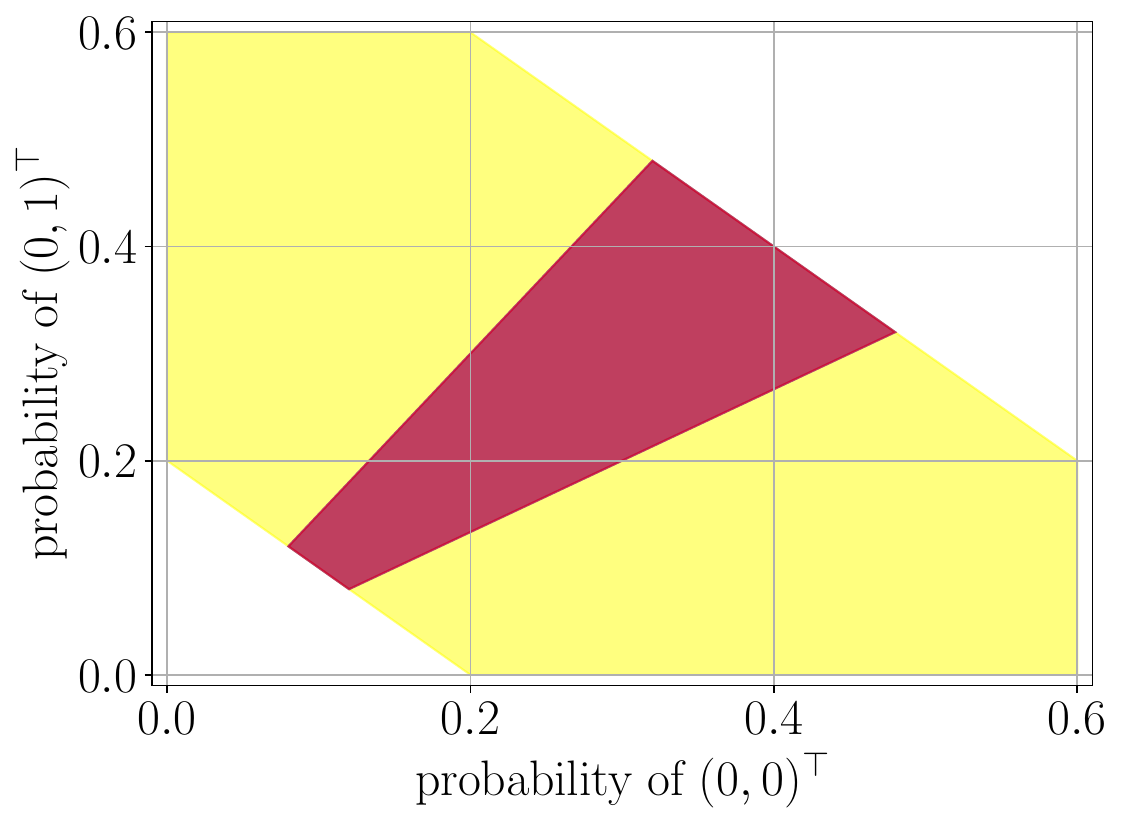} &
        \includegraphics[height = 3.8cm]{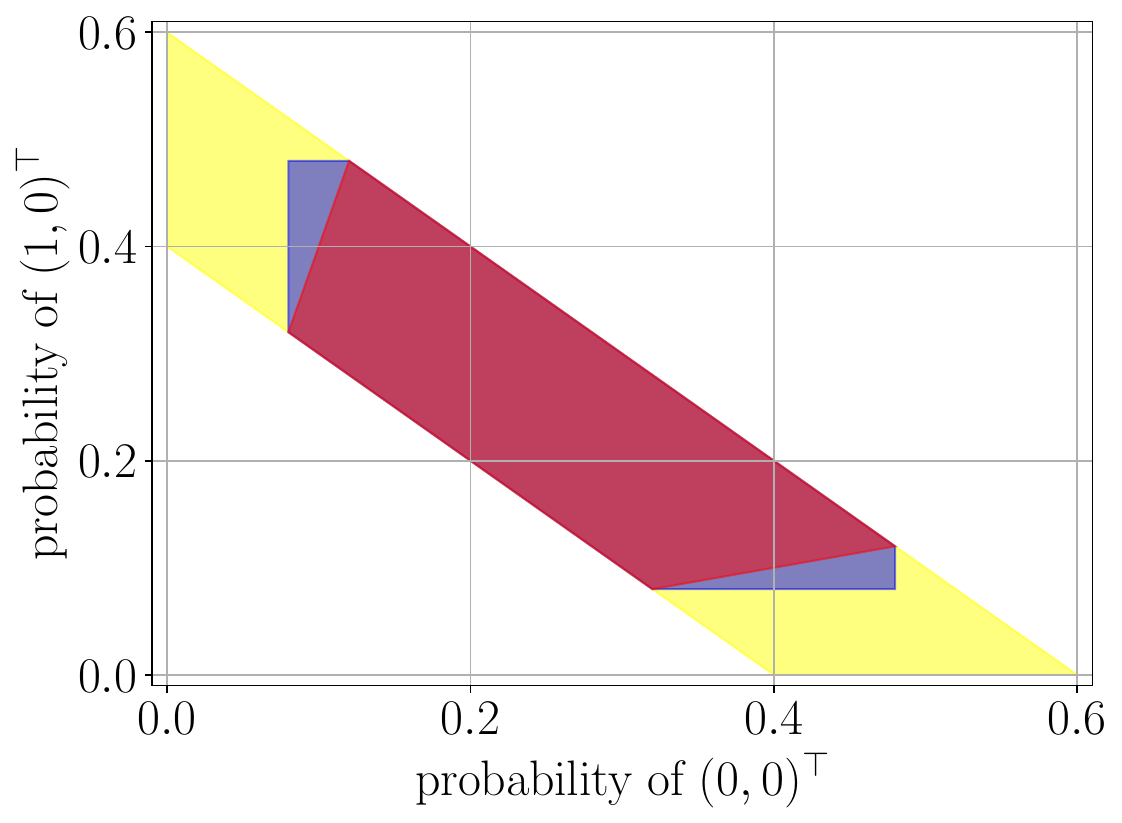} &
        \includegraphics[height = 3.8cm]{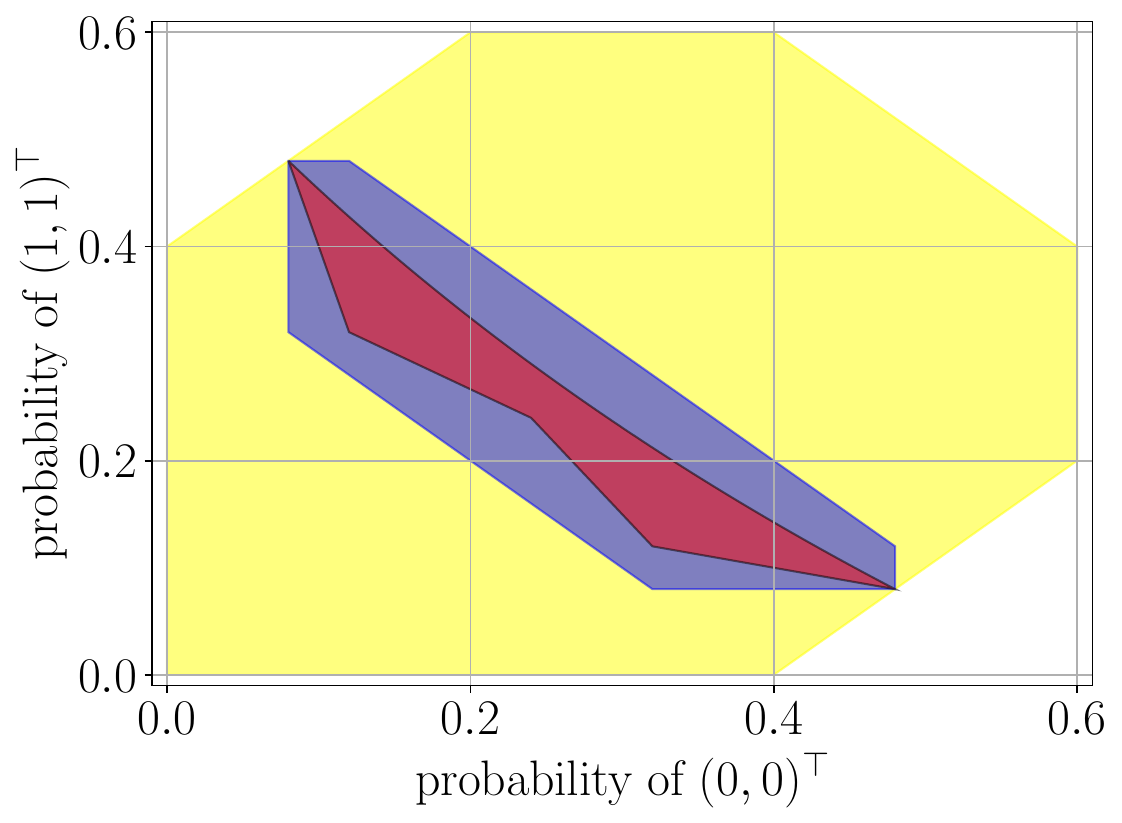} &
        \raisebox{-0.0cm}{\includegraphics[height = 3.8cm]{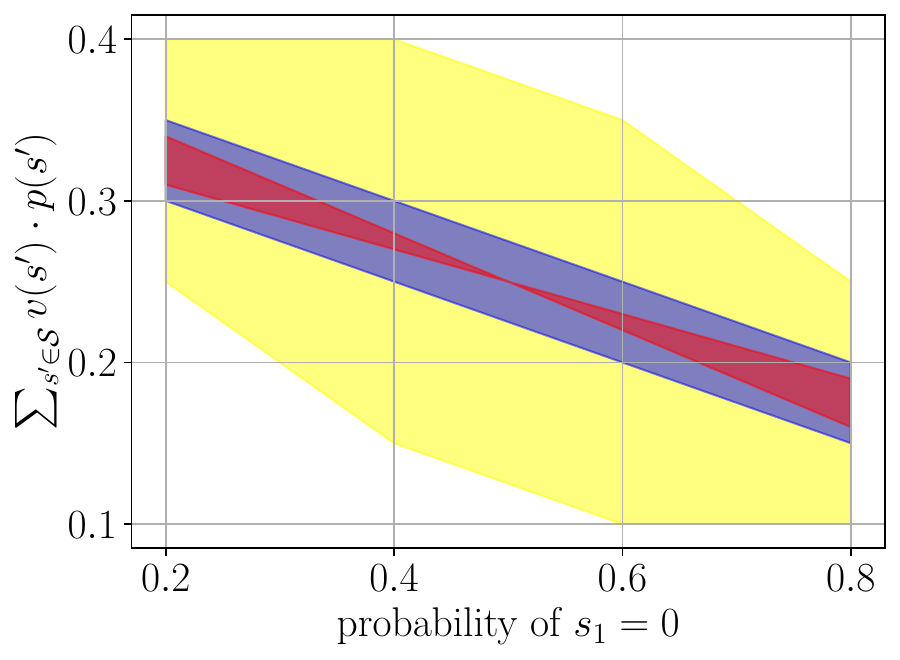}}
    \end{tabular}
    \caption{Different two-dimensional projections of the four-dimensional ambiguity sets $\mathcal{P} (s, a)$ from Example~\ref{ex:factored_nonfactored} (first three subplots) as well as the range of possible expectations $\sum_{s' \in \mathcal{S}} p(s') \cdot v(s')$ (fourth subplot) for the factored representation (red), the non-factored representation (yellow) and the lift defined later in this section (blue). \label{fig:factored_amb_sets}}
\end{figure}

\begin{ex}[Factored and Non-Factored Ambiguity Sets]\label{ex:factored_nonfactored}
    Consider an RFMDP with $N = 2$ sub-states, $\mathcal{S}_1, \mathcal{S}_2 = \mathbb{B}$, $\mathcal{A} = \{ 1 \}$, $\phi ((s, a)) = (s, a)$, and $\mathcal{P}^1 (s, a) = \{ p_1 (\cdot | s, a) \in \Delta (\mathcal{S}_1) \, : \, p_1 (s'_1 = 1 | s, a) \in$ $[0.2, 0.8] \}$ and $\mathcal{P}^2 (s, a) = \left\{ p_2 (\cdot | s, a) \in \Delta (\mathcal{S}_2) \, : \, p_2 (s'_2 = 1 | s, a) \in [0.4, 0.6] \right\}$. The first three plots in Figure~\ref{fig:factored_amb_sets} visualize different projections of the corresponding factored ambiguity set (in red) and non-factored ambiguity set (in yellow). Likewise, the fourth plot visualizes the range of possible expectations $\sum_{s' \in \mathcal{S}} p(s') \cdot v(s')$ across $p \in \mathcal{P} (s, a)$ for the value function $v (s') = 1/4$ if $s'_2 = 1$; $= 0.5 \cdot s'_1$ otherwise, in the factored and non-factored RFMDPs. The plots show that the factored ambiguity set is non-convex and substantially smaller than its non-factored counterpart.
\end{ex}

Recall from the Section~\ref{sec:robust_FMDP} that computing the worst-case expectation $\min_{p_{sa}} \sum_k w_k \sum_{s'} \nu_k (s') \cdot p_{sa} (s')$ over a non-factored ambiguity set amounts to solving an LP. This is no longer the case for factored ambiguity set since they involve products of marginal distributions. Indeed, computing worst-case expectations over factored ambiguity sets is computationally challenging.

\begin{prop}\label{prop:factored_ambiguity_hard}
    The worst-case expectation problem
    \begin{equation*}
        \begin{array}{l@{\quad}l}
            \displaystyle \mathop{\text{\emph{minimize}}}_p & \displaystyle \sum_{s' \in \mathcal{S}} p (s' | \phi (s, a)) \cdot v(s') \\
            \displaystyle \text{\emph{subject to}} & \displaystyle p \in \mathcal{P} (\phi (s, a))
        \end{array}
    \end{equation*}
    over a factored ambiguity set and a value function $v : \mathcal{S} \rightarrow \mathbb{R}$ is PPAD-hard even when $N = 2$.
\end{prop}

The complexity class PPAD (``polynomial parity arguments on directed graphs'') is closely related but different to the complexity class NP. Similar to NP-hard problems, PPAD-hard problems are believed to be intractable \citep{P94:ppad}. The complexity class PPAD is well-known for comprising the computation of Nash equilibria, which we use in our reduction proof as well.

\begin{figure}[tb]
    \vspace{-2.85cm}
    \includegraphics[width = \textwidth]{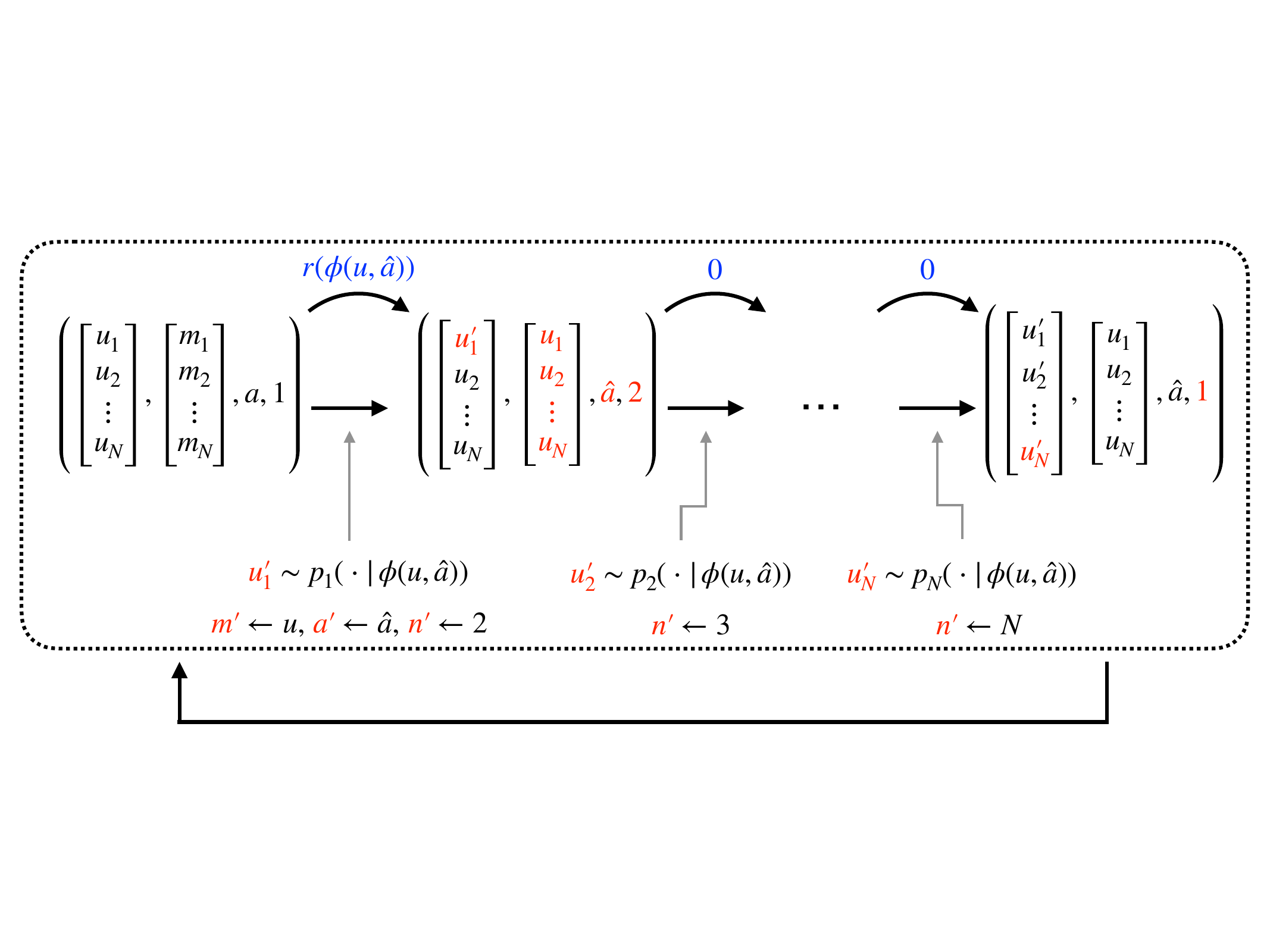}
    \vspace{-3cm}
    \caption{The lift of a factored RFMDP. The figure shows $N$ subsequent transitions in the lifted RFMDP that correspond to one transition in the non-lifted RFMDP. The one-step rewards are shown in blue, and any state updates are highlighted in red. \label{fig:lift}}
\end{figure}

Instead of solving factored RFMDPs directly, we will conservatively approximate them via \emph{lifted} non-factored RFMDPs. We present an informal description of this lift in the following and relegate its formal definition to Section~A of the e-companion. Fix a factored RFMDP $(\mathcal{S}, \mathcal{A}, \phi, q, \mathcal{P}, r, \gamma)$. We define its corresponding non-factored lifted RFMDP $(\hat{\mathcal{S}}, \hat{\mathcal{A}}, \hat{\phi}, \hat{q}, \hat{\mathcal{P}}, \hat{r}, \hat{\gamma})$---or, for short, its \emph{lift}---as follows. The lifted state space is $\hat{\mathcal{S}} = \mathcal{S}^2 \times \mathcal{A} \times [N]$, where each state $\hat{s} = (u, m, a, n)$ comprises an update block $u \in \mathcal{S}$, a memory block $m \in \mathcal{S}$, an action block $a \in \mathcal{A}$ and a timer block $n \in [N]$. The lifted action space is $\hat{\mathcal{A}} = \mathcal{A}$. Fix any state $s^0 \in \mathcal{S}$ and action $a^0 \in \mathcal{A}$ in the non-lifted factored RFMDP. The lifted initial state distribution is defined as $q (\hat{s}) = q (u)$ for any $\hat{s} = (u, s^0, a^0, 1)$ with $u \in \mathcal{S}$, whereas $q (\hat{s}) = 0$ otherwise. Assume that the system is in state $\hat{s} = (u, m, a, n)$ and that action $\hat{a}$ is taken. The system then enters the next state $\hat{s}' = (u', m', a', n')$ according to the following dynamics:
\begin{itemize}
    \item If $n = 1$, then $u'_1 \sim p_1 (\cdot | \phi (u, \hat{a}))$ for some marginal transition kernel $p_1 \in \mathcal{P}^1 (\phi (u, \hat{a}))$, $u'_\eta = u_\eta$ for $\eta \neq 1$, $m' = u$, $a' = \hat{a}$ and $n' = 2$.
    \item If $1 < n < N$, then $u'_n \sim p_n (\cdot | \phi (m, a))$ for some marginal transition kernel $p_n \in \mathcal{P}^n (\phi (m, a))$, $u'_\eta = u_\eta$ for $\eta \neq n$, $m' = m$, $a' = a$ and $n' = n + 1$.
    \item If $n = N$, then $u'_N \sim p_N (\cdot | \phi (m, a))$ for some marginal transition kernel $p_N \in \mathcal{P}^N (\phi (m, a))$, $u'_\eta = u_\eta$ for $\eta \neq N$, $m' = s^0$, $a' = a^0$ and $n' = 1$.
\end{itemize}
For $\hat{s} = (u, m, a, n)$, we define the lifted rewards $\hat{r} (\hat{\phi} (\hat{s}, \hat{a})) = r (\phi (u, \hat{a}))$ if $(m, a, n) = (s^0, a^0, 1)$ and $\hat{r} (\hat{\phi} (\hat{s}, \hat{a})) = 0$ otherwise, and we set the lifted discount factor to $\hat{\gamma} = \gamma^{1/N}$. The lifted RFMDP is non-factored since for each lifted state-action pair $(\hat{s}, \hat{a}) \in \hat{\mathcal{S}} \times \hat{\mathcal{A}}$, all but one of the marginal ambiguity sets $\mathcal{P}^n (\phi(\hat{s}, \hat{a}))$, $n \in [N]$, are singleton sets that contain Dirac distributions, and thus the factored and non-factored ambiguity set $\mathcal{P} (\phi (\hat{s}, \hat{a}))$ coincide. Moreover, Section~A of the e-companion argues that any low-scope properties of the factored RFMDP are preserved by its lift. Both properties imply that we can directly employ the methods developed in Section~\ref{sec:robust_FMDP} to compute near-optimal policies for the  lifted RFMDP. Figure~\ref{fig:lift} illustrates the lift of a factored RFMDP.

We now investigate how a factored RFMDP relates to its lift. To this end, fix any non-lifted RFMFP instance $\mathcal{R} = (\mathcal{S}, \mathcal{A}, \phi, q, \mathcal{P}, r, \gamma)$, let $\Pi = \{ \pi : \mathcal{S} \rightarrow \mathcal{A} \}$ and $\hat{\Pi} = \{ \hat{\pi} : \hat{\mathcal{S}} \rightarrow \hat{\mathcal{A}} \}$ denote the sets of non-lifted and lifted policies, respectively, and consider the following three objectives:
\begin{enumerate}
    \item[\emph{(i)}] The \textbf{non-factored objective} $\mathrm{NF} (\pi; \mathcal{R}) = \inf_{p \in \mathcal{P}^\mathrm{NF}} \mathbb{E} \big[ \sum_{t \in \mathbb{N}_0} \gamma^t \cdot r (\tilde{s}^t, \pi (\tilde{s}^t)) \, | \, \tilde{s}^0 \sim q \big]$, with optimal value $\displaystyle \mathrm{NF}^\star (\mathcal{R}) = \sup_{\pi \in \Pi} \mathrm{NF} (\pi; \mathcal{R})$, that combines the marginal ambiguity sets $\mathcal{P}^n (\phi (s, a))$, $(s, a) \in \mathcal{S} \times \mathcal{A}$, in a non-factored way (\emph{cf.}~Section~\ref{sec:robust_FMDP}).
    \item[\emph{(ii)}] The \textbf{factored objective} $\mathrm{F} (\pi; \mathcal{R}) = \inf_{p \in \mathcal{P}^\mathrm{F}} \mathbb{E} \big[ \sum_{t \in \mathbb{N}_0} \gamma^t \cdot r (\tilde{s}^t, \pi (\tilde{s}^t)) \, | \, \tilde{s}^0 \sim q \big]$, with optimal value $\displaystyle \mathrm{F}^\star (\mathcal{R}) = \sup_{\pi \in \Pi} \mathrm{F} (\pi; \mathcal{R})$, that combines the marginal ambiguity sets $\mathcal{P}^n (\phi (s, a))$, $(s, a) \in \mathcal{S} \times \mathcal{A}$, in a factored way (as presented in the beginning of this section).
    \item[\emph{(iii)}] The \textbf{lifted non-factored objective} $\mathrm{L} (\hat{\pi}; \mathcal{R}) = \inf_{\hat{p} \in \hat{\mathcal{P}}} \mathbb{E} \big[ \sum_{t \in \mathbb{N}_0} \hat{\gamma}^t \cdot \hat{r} (\tilde{\hat{s}}^t, \hat{\pi} (\tilde{\hat{s}}^t)) \, | \, \tilde{\hat{s}}^0 \sim \hat{q} \big]$, with optimal value $\displaystyle \mathrm{L}^\star (\mathcal{R}) = \sup_{\hat{\pi} \in \hat{\Pi}} \mathrm{L} (\hat{\pi}; \mathcal{R})$, where $\hat{\Pi} : \hat{\mathcal{S}} \rightarrow \hat{\mathcal{A}}$ (as just presented above).
\end{enumerate}
Note that the above objectives compute the exact worst-case expected total reward, that is, they do not employ value function approximations. This ensures comparability of the objectives, whose associated RFMDPs will employ different value function approximations due to their differing state space dimensions. For any policy $\pi : \mathcal{S} \rightarrow \mathcal{A}$, we define the following associated set of policies in the lifted non-factored RFMDP:
\begin{equation*}
    \hat{\Pi} (\pi) = \left\{
        \hat{\pi} : \hat{\mathcal{S} } \rightarrow \hat{\mathcal{A}} \, : \, \hat{\pi} (\hat{s}) = \pi (u) \;\; \forall \hat{s} = (u, s^0, a^0, 1) \in \hat{\mathcal{S}}
    \right\}.
\end{equation*}
We next show that $\hat{\Pi} (\pi)$, $\pi \in \Pi$, partitions $\hat{\Pi}$, and that all lifted policies in the same class $\hat{\Pi} (\pi)$ generate the same rewards in the lifted non-factored RFMDP.

\begin{obs}\label{obs:lifted_policies}
    The policy sets $\hat{\Pi} (\pi)$, $\pi \in \Pi$, satisfy the following properties:
    \begin{enumerate}
        \item[\emph{(i)}] $\{ \hat{\Pi} (\pi) : \pi \in \Pi \}$ partitions $\hat{\Pi}$.
        \item[\emph{(ii)}] $\mathrm{L} (\hat{\pi}; \mathcal{R}) = \mathrm{L} (\hat{\pi}'; \mathcal{R})$ whenever $\hat{\pi}, \hat{\pi}' \in \hat{\Pi} (\pi)$ for some $\pi \in \Pi$.
    \end{enumerate}
\end{obs}

We now study the tightness of the lift.

\begin{thm}[Relationship between Ambiguity Sets]\label{thm:rel_ambiguity_sets}
    ~\begin{enumerate}
        \item[\emph{(i)}] For any RFMDP instance $\mathcal{R}$, any non-lifted policy $\pi \in \Pi$ and any lifted policy $\hat{\pi} \in \hat{\Pi} (\pi)$, we have $\mathrm{NF} (\pi; \mathcal{R}) \leq \mathrm{L} (\hat{\pi}; \mathcal{R}) \leq \mathrm{F} (\pi; \mathcal{R})$, where none, either one or both inequalities can be tight.
        \item[\emph{(ii)}] There are RFMDP instances $\mathcal{R}$ with 0/1-rewards such that $\mathrm{NF}^\star (\mathcal{R}) = \mathrm{L}^\star (\mathcal{R}) = 0$ while $\mathrm{F}^\star (\mathcal{R}) = \gamma \cdot (1 - \gamma)^{-1} \cdot (1 - 2^{-\lfloor N / 2 \rfloor})$, where $N$ denotes the number of sub-states in $\mathcal{R}$.
        \item[\emph{(iii)}] There are RFMDP instances $\mathcal{R}$ with 0/1-rewards such that $\mathrm{NF}^\star (\mathcal{R}) = 0$ while $\mathrm{L}^\star (\mathcal{R}) = \mathrm{F}^\star (\mathcal{R}) = \gamma \cdot (1 - \gamma)^{-1} \cdot (1 - 2^{-(N-1)})$, where $N$ denotes the number of sub-states in $\mathcal{R}$.
    \end{enumerate}
\end{thm}

Note that for an RFMDP with 0/1-rewards that starts in a state where all actions have a reward of 0, the maximally achievable expected total reward is $\gamma \cdot (1 - \gamma)^{-1}$. Thus, the statements \emph{(ii)} and \emph{(iii)} of Theorem~\ref{thm:rel_ambiguity_sets} show that the non-factored, factored and lifted non-factored RFMDPs can essentially (up to an exponentially small factor) differ maximally in such problem classes. The hope is that for non-pathological problem instances, the lift behaves much more similarly to the factored RFMDP rather than its non-factored variant. This is illustrated in the next example and further explored in our numerical experiments in the next section.

\begin{ex}[Factored and Non-Factored Ambiguity Sets, Cont'd]
    For Example~\ref{ex:factored_nonfactored}, the blue shaded regions in Figure~\ref{fig:factored_amb_sets} visualize the conservative approximation offered by the corresponding non-factored lifted RFMDP. We can see that in this example, the lift provides a significantly tighter approximation of the factored RMFDP than the ordinary non-factored RFMDP.
\end{ex}

\section{Numerical Experiments}\label{sec:num_experiments}

We investigate the numerical performance of our proposed methodology on two problems: the classical system administrator problem from the FMDP literature \citep{SP01:direct_value, GKPV03:efficient_solution} and a more contemporary factored multi-armed bandit problem. Section~\ref{sec:num_experiments:fmdps} introduces the system administrator problem and compares the performance of our cutting plane scheme from Section~\ref{sec:FMDP} with the variable-elimination method of \cite{GKPV03:efficient_solution}. Section~\ref{sec:num_experiments:rfmdps} extends the problem to a data-driven setting where the transition kernel needs to be estimated from past data, and it elucidates how the RFMDPs from Section~\ref{sec:robust_FMDP} can construct robust policies that enjoy a superior out-of-sample performance. Section~\ref{sec:num_experiments:bandits}, finally, studies a factored multi-armed bandit problem and compares FMDPs with methods from the related literature on weakly coupled MDPs. All algorithms are implemented in C++ using Gurobi 9.5.1 and run on Intel Xeon 2.20GHz cluster nodes with 64 GB dedicated main memory in single-core mode. Our documented sourcecodes, data sets and additional results can be found online (\emph{cf.}~Footnote~\ref{fn:github}).

\subsection{Factored Markov Decision Processes}\label{sec:num_experiments:fmdps}

\begin{figure}[tb]
    \includegraphics[width = \textwidth]{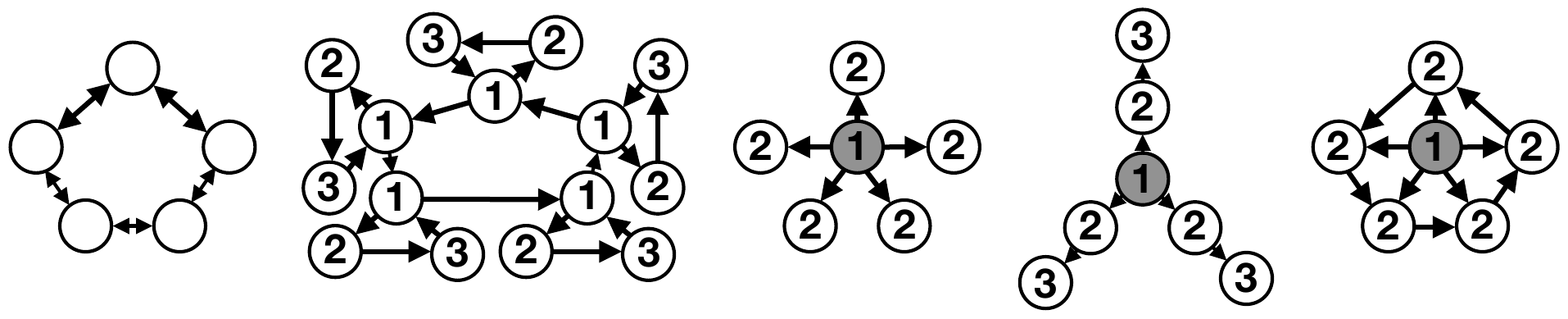}
    \caption{Network topologies for the system administrator problem (from left to right): bidirectional ring, ring of rings, star, 3 legs and ring and star. The first two topologies are decentralized, whereas the others rely on centralized servers (shown as shaded nodes). \label{fig:sysadmin:topologies}}
\end{figure}

Consider a system administrator that maintains a network of $N$ computers. The computers are organized into one of five distinct topologies, as illustrated in Figure~\ref{fig:sysadmin:topologies}.  In this figure, the computers are represented as nodes, and their connections are depicted as unidirectional or bidirectional arcs. In every time period, each computer is in one of three states: inoperative, semi-operational and fully operational. The state of each machine deteriorates stochastically over time, depending on its own current state as well as the state of its predecessors in the topology. To uphold the network's functionality, the system administrator can reboot and thus restore to a fully operational state up to $2$ computers in each time period. The administrator earns a per-period reward of 1 unit (2 units) for every computer in a semi-operational (fully operational) state.

The problem can be formulated as an FMDP with the sub-state spaces $\mathcal{S}_n = \{ \mathrm{e}_1, \mathrm{e}_2, \mathrm{e}_3 \} \subseteq \mathbb{B}^3$, $n \in [N]$, and the action space $\mathcal{A} = \{ a \in \mathbb{B}^N : \sum_{n \in [N]} a_n \leq 2 \}$. Under a na\"ive implementation, each set of transition probabilities $p_n$ has a scope between 4 and 10 (accounting for the computer's current state and those of the predecessor nodes, along with the action linked to the computer), whereas each reward $r_n$ has a scope of 3 (accounting for the state of the associated computer). In our experiments we employ scope-1, scope-2 or scope-3 value functions. The scope-1 value functions comprise $3$ basis functions for each computer $n \in [N]$: $\nu_{(n, i)} = \mathbf{1} [ s_n = \mathrm{e}_i ]$ for $i \in \{ 1, 2, 3 \}$. The scope-2 value functions comprise $9$ basis functions for each computer pair $(n, n \oplus 1)$, where $n \oplus 1 = n + 1$ for $n \in [N-1]$ as well as $N \oplus 1 = 1$: $\nu_{(n, i), (n \oplus 1, j)} = \mathbf{1} [ s_n = \mathrm{e}_i \text{ and } s_{n \oplus 1} = \mathrm{e}_j ]$ for $i, j \in \{ 1, 2, 3\}$. The scope-3 value functions, finally, comprise $27$ basis functions for each computer triplet $(n, n \oplus 1, n \oplus 2)$ that are defined analogously.

We solve our FMDP problem~\eqref{opt:fmdp:basis_fcts_lp_simplified} using Algorithm~\ref{alg:CP}. The master problem is an LP that is solved directly via Gurobi. The subproblem is an MILP, which is initially solved heuristically by identifying violated constraints $(s, a) \in \mathcal{S} \times \mathcal{A}$ via a coordinate-wise descent approach that starts at a randomly selected constraint $(s, a) \in \mathcal{S} \times \mathcal{A}$ and iteratively replaces each sub-state $s_n$ and sub-action $a_n$, $n \in [N]$, with a new one that maximally increases the constraint violation until a partial optimum is reached (\emph{cf.}~Remark~\ref{rem:practical_considerations} and \citealp{KPK07:biconvex} as well as \S 3.3 of \citealp{B09:conv_opt_theory}). If the coordinate-wise descent fails to identify a violated constraint, we solve the subproblem as an MILP (\emph{cf}.~Theorem~\ref{thm:solution_subproblem}) with early termination once Gurobi has \emph{(i)} identified a constraint with a sufficiently large violation; \emph{(ii)} concluded that all constraints are satisfied by the current basis function approximation; or \emph{(iii)} exceeded the time limit of $100 + 3N$ seconds.

\begin{figure}[tb]
    $\mspace{-100mu}$
    \begin{tabular}{c@{}c@{}c@{}c@{}c}
        \includegraphics[width=4cm]{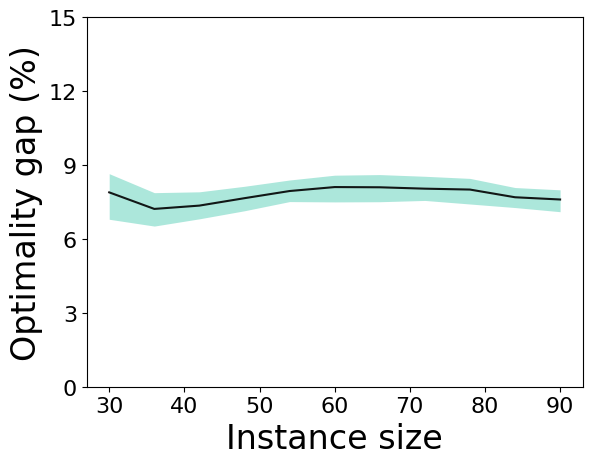} &
        \includegraphics[width=4cm]{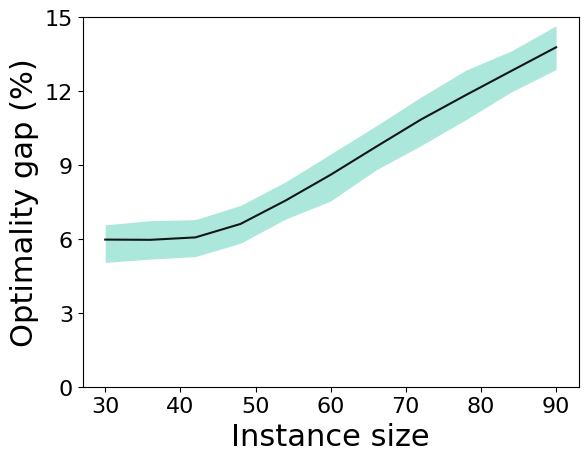} &
        \includegraphics[width=4cm]{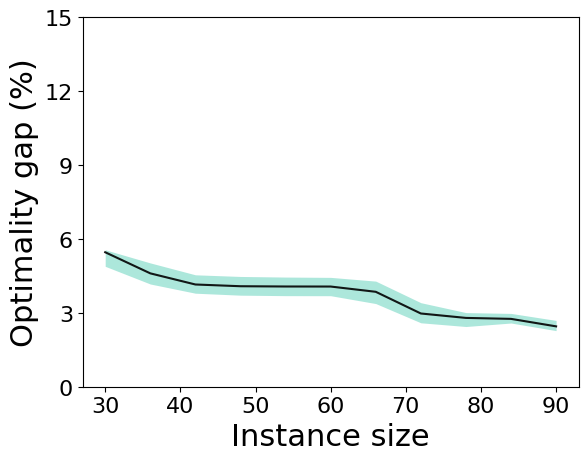} &
        \includegraphics[width=4cm]{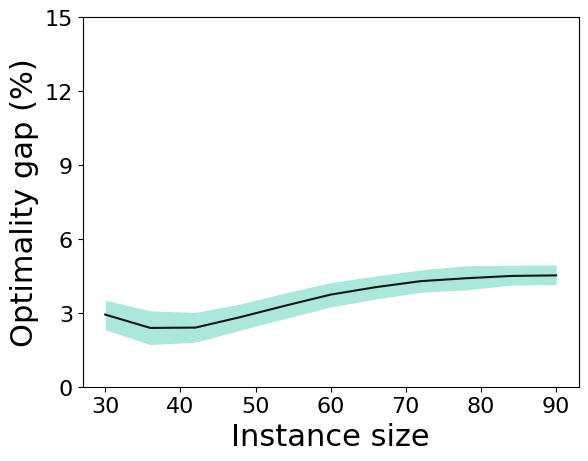} &
        \includegraphics[width=4cm]{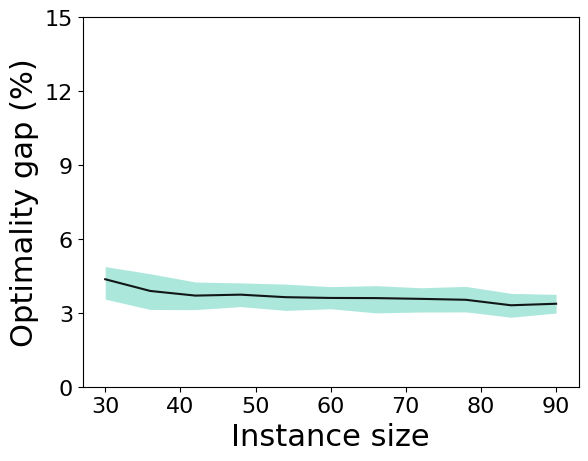}
    \end{tabular}
    \caption{Optimality gaps (in \%) for our cutting plane approach on the 5 network topologies of Figure~\ref{fig:sysadmin:topologies} as a function of the instance size (number $N$ of computers). \label{fig:opt_gaps}}
\end{figure}

\begin{figure}[tb]
$\mspace{-100mu}$
    \begin{tabular}{c@{}c@{}c@{}c@{}c}
        \includegraphics[width=4cm]{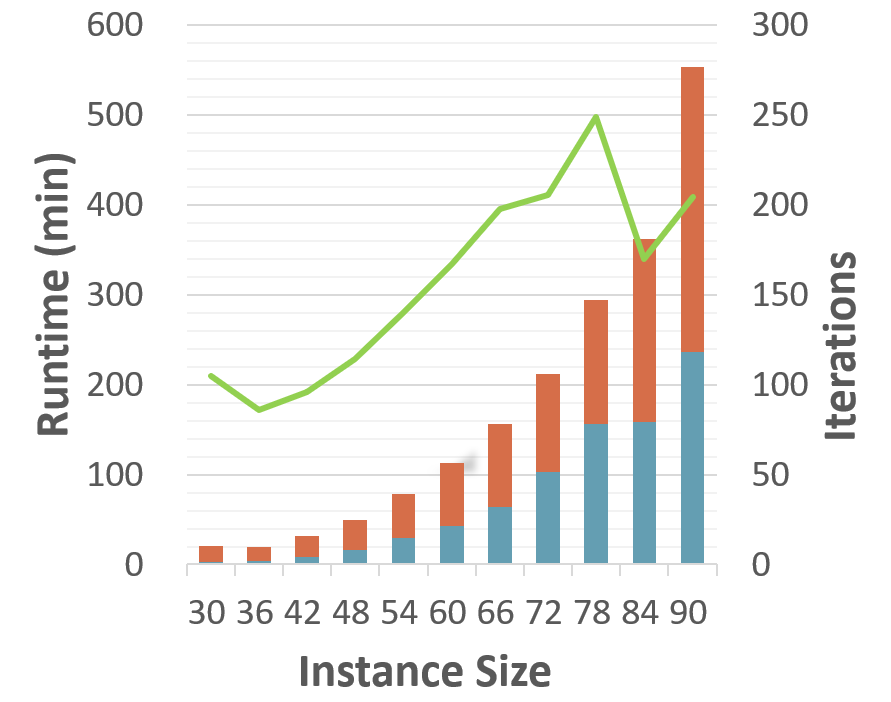} &
        \includegraphics[width=4cm]{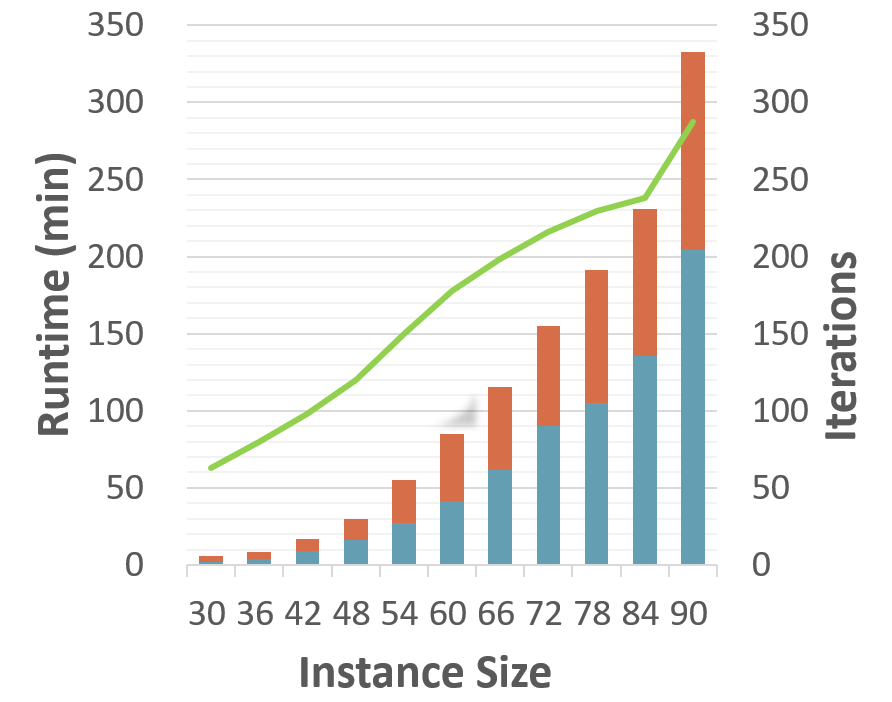} &
        \includegraphics[width=4cm]{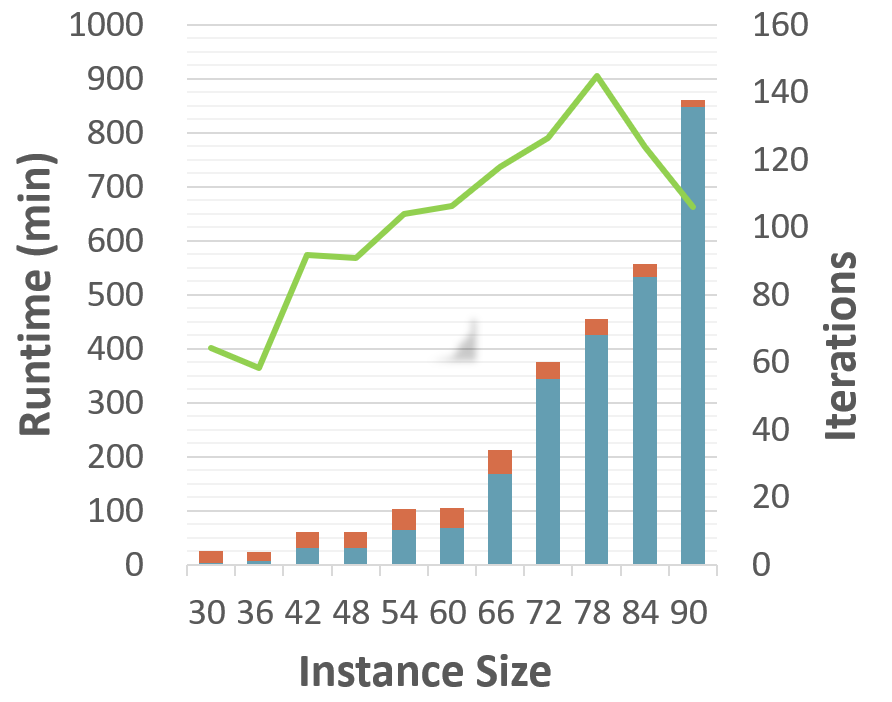} &
        \includegraphics[width=4cm]{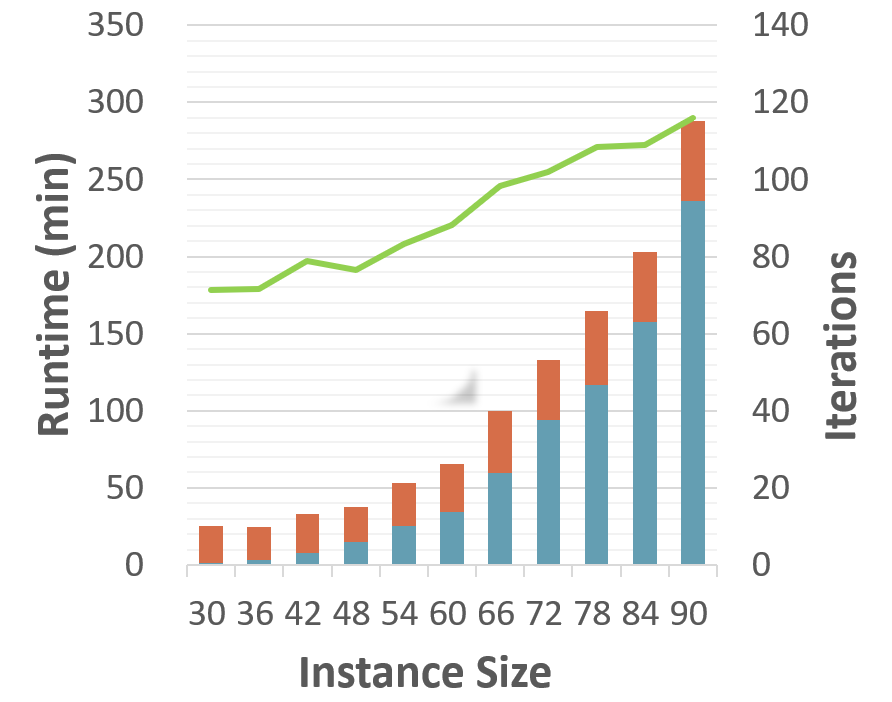} &
        \includegraphics[width=4cm]{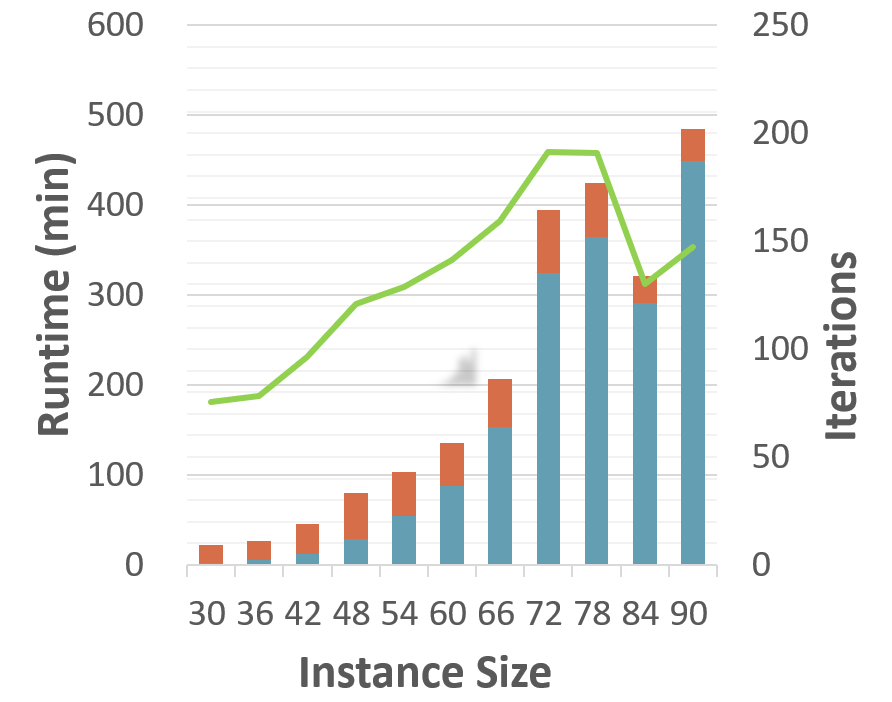}
    \end{tabular}
    \caption{Runtimes (bars) and iteration numbers (bold lines) for our cutting plane approach on the 5 network topologies of Figure~\ref{fig:sysadmin:topologies} as a function of the instance size (number $N$ of computers). The runtime bars comprise the runtimes for the subproblem (upper red bars) and the master problem (lower blue bars). \label{fig:runtimes}}
\end{figure}

Figure~\ref{fig:opt_gaps} reports the optimality gaps of our method for the different network topologies of Figure~\ref{fig:sysadmin:topologies} under scope-2 value functions. In the figure, the optimality gaps are computed as $100\% \cdot (\mathrm{UB} - \mathrm{LB}) / \mathrm{LB}$, where $\mathrm{UB}$ denotes the upper bound on the expected total reward of the unknown optimal policy that is provided by the optimal value of terminal master problem of Algorithm~\ref{alg:CP} (\emph{cf.}~the discussion after Proposition~\ref{prop:value_function_reformulation}), and  $\mathrm{LB}$ is an approximation of the expected total reward of the policy computed by Algorithm~\ref{alg:CP} that is determined via Monte Carlo simulation ($200$ statistically independent runs comprising $200$ transitions each). We repeat the experiment for every network topology and every instance size (measured by the number $N$ of computers) $250$ times, and we report the median optimality gaps (bold lines) as well as the 10\%- and 90\%-quantiles (shaded regions). The figure shows that for most network topologies, the optimality gaps remain well below 10\% and do not grow with the size of the instances. An exception is the ring of rings topology, where each computer can connect with as many as $4$ other computers. This complexity suggests that scope-2 value functions lack the flexibility to closely approximate the true value function.

Figure~\ref{fig:runtimes} presents the runtimes and iteration counts of our method for the different network topologies of Figure~\ref{fig:sysadmin:topologies} under scope-2 value functions. The figure reports the median runtimes across $250$ runs, split up into the times spent on solving the subproblem (upper red bars) and the master problem (lower blue bars). The median iteration numbers are reported as bold lines. The figure shows that both the number of iterations and the overall time spent increases with problem size. There is no clear trend concerning the computational bottleneck: for some network topologies, the majority of time is spent on solving the master problem, whereas in some topologies the subproblem consumes an increasing share of the runtime as instances become larger.

\begin{table}[tb]
    $\mspace{-28mu}$
    \begin{tabular}{cr|ccccccc}
        & & 10 & 15 & 20 & 25 & 30 & 35 & 40 \\ \hline
        \multirow{12}{*}{\shortstack{\\ \textbf{rewards} \\ \textbf{(gaps)}}} & random & 336.1 & 444.9 & 528.6 & 620.8 & 716.6 & 816.0& 914.2 \\[1mm]
        & priority & 338.1 & 449.8 & 536.8 & 633.8 & 734.0 & 835.1 & 938.1\\[1mm]
        & \multirow{2}{*}{VE (scope-1)} & 342.9 & 454.9 & 534.9  & 633.7  & 737.0 & 838.8 & 941.1\\
        && (9.2\%) & (23.2\%) & (33.6\%) & (35.9\%)  & (36.4\%)  &  (36.9\%) & (37.1\%)\\[1mm]
        & \multirow{2}{*}{VE (scope-2)} & \textbf{344.5} & 482.2  & \multirow{2}{*}{---} & \multirow{2}{*}{---} & \multirow{2}{*}{---} & \multirow{2}{*}{---} & \multirow{2}{*}{---} \\
        &&(8.7\%) & (7.7\%) &   &  &  &  &  \\[1mm]
        & \multirow{2}{*}{ours (scope-1)} & 340.4  &  454.0 & 533.2 & 619.4 & 726.1 & 821.0 & 929.3\\ 
        && (10.0\%) & (23.4\%) & (35.1\%) & (37.6\%) & (33.5\%) & (33.4\%)  & (31.9\%)\\ [1mm]
        & \multirow{2}{*}{ours (scope-2)} & 343.1  &  481.6 & 587.2 &682.2 & 778.2 & 877.3& 975.4\\ 
        && (9.1\%) & (8.2\%) & (7.3\%) & (8.5\%) & (8.2\%) & (8.6\%) & (8.5\%)\\ [1mm]
        & \multirow{2}{*}{ours (scope-3)} & \textbf{344.5}  &  \textbf{482.5} & \textbf{592.8} & \textbf{690.0} & \textbf{789.0} & \textbf{888.0} & \textbf{988.3}\\ 
        && (6.1\%)& (5.2\%) & (2.8\%) & (3.3\%) & (3.1\%) & (3.4\%) & (3.1\%)\\ [1mm]
        \hline
        \multirow{5}{*}{\textbf{runtimes} \textbf{(s)}} & VE (scope-1) & 4.1 & 59.2 & 235.2 & 1,149.3 & 2,404.0 & 16,134.2 & 56,914.1 \\[1mm]
        & VE (scope-2) & 899.0 & 19,627.1 & --- & --- & --- & ---& \\[1mm]
        & ours (scope-1) & 86.8 & 146.3 & 161.9 & 418.7 & 243.3 & 400.3 & 686.2\\[1mm] 
        & ours (scope-2) & 131.2 & 259.9 & 170.9 & 622.8 & 404.8 & 798.3 & 1,437.9\\[1mm]
        & ours (scope-3) & 140.8 & 469.3 & 901.6 & 1,571.5 & 2,134.6 & 4,063.5 & 7,323.8 \\
        \hline \hline
    \end{tabular}
    \caption{Expected total rewards, optimality gaps and runtimes for the variable-elimination method (`VE'), two na\"ive benchmark strategies (`random' and `priority') as well as our cutting plane approach (`ours') on the bidirectional ring topology. Absent results (``---'') indicate that the corresponding experiments have run out of memory. \label{tab:ve-comparison}}
\end{table}

\begin{table}[tb]
\centering
    \begin{tabular}{r|ccccccc}
       \textbf{rewards} & 12 & 24 & 36 & 48 & 60 \\ \hline
         random & 413.9 & 694.7 & 917.9 & 1,149.7 & 1,386.2  & \\[1mm]
         priority & 416.2 & 714.8 & 958.9  & 1,208.3 & 1,454.6 & \\[1mm]
         level & 428.5 & 787.2 & 1,071.1  & 1,297.9 & 1,524.8 & \\[1mm]
         \hline 
         \multirow{2}{*}{ours (scope-1)} & 429.4  & 801.6  & 1,122.1 & 1,401.4 &  1,670.8& \\ 
        & (5.8\%) & (10.1\%) &  (11.6\%) & (13.2\%)  & (15.2\%) & \\[1mm] 
        \multirow{2}{*}{ours (scope-2)} &  430.9 & 804.8  & 1,131.12 & 1,417.1 & 1,679.5 & \\ 
        & (4.4\%) &(6.6\%) & (7.1\%) & (7.9\%)  & (8.7\%) & \\[1mm] 
        \multirow{2}{*}{ours (scope-3)} & \textbf{431.5} &  \textbf{807.3} & \textbf{1,142.3} & \textbf{1,425.7} & \textbf{1,688.3} & \\ 
        & (3.2\%) & (4.3\%) & (4.7\%)& (5.2\%) & (5.8\%) & \\ 
        \hline  \hline
    \end{tabular}
    \caption{Expected total rewards for the three benchmark strategies as well as our cutting plane approach on the ring of rings topology. In contrast to the benchmark strategies, our approach additionally provides optimality gaps, which are listed as well. \label{tab:ve-comparison-2}}
\end{table}

To put our results into perspective, we compare the policies and runtimes of our cutting plane approach (`ours') with those of the variable-elimination (`VE') method by \citet{GKPV03:efficient_solution} on the bidirectional ring topology. We also include two na\"ive benchmark strategies: a `random' strategy that selects the computers to reboot randomly from those that are inoperative or semi-operational, and a `priority' strategy that prioritizes the rebooting of inoperative computers over that of semi-operational ones. The results are summarized in Table~\ref{tab:ve-comparison}. Missing entries correspond to experiments that ran out of memory. The table shows that our cutting plane approach results in similar optimality gaps as the variable-elimination method at substantially lower runtime and memory requirements. The table also shows that the two na\"ive benchmark strategies---apart from being unable to provide optimality gaps---may underperform by up to 10\% relative to our scope-3 strategy. This is further emphasized in Table~\ref{tab:ve-comparison-2}, which compares the benchmark strategies with our cutting plane approach for larger instances sizes on the ring of rings topology. Here, even the more sophisticated `priority' strategy is outperformed by our scope-1 value functions, and our scope-3 value functions generate up to 19\% more expected total rewards. The ring of rings topology allows us to to distinguish between nodes of different levels as shown in Figure~\ref{fig:sysadmin:topologies}, ranging from 1 (highest) to 3 (lowest). We therefore also include a `level' heuristic which breaks ties among inoperative and semi-operational computers based on the priorities of the involved computers. We see that this heuristic, too, is outperformed by all of our policies.

\begin{table}[tb]
    $\mspace{-47mu}$
    \begin{tabular}{crr@{}|@{}lcccr@{}|@{}lccr@{}|@{}lcccr@{}|@{}lccr}
        \multicolumn{4}{c}{} & \multicolumn{3}{c}{ring of rings} & \multicolumn{1}{r}{} && \multicolumn{2}{c}{star} & \multicolumn{1}{r}{} && \multicolumn{3}{c}{3 legs} & \multicolumn{1}{r}{} && \multicolumn{2}{c}{ring and star} \\ \cline{5-7} \cline{10-11} \cline{14-16} \cline{19-20}
        \multicolumn{4}{c}{} & 1 & 2 & \multicolumn{1}{c}{3} & \multicolumn{1}{r}{} && 1 & \multicolumn{1}{c}{2} & \multicolumn{1}{r}{} && 1 & 2 & \multicolumn{1}{c}{3} & \multicolumn{1}{r}{} && 1 & 2 \\ \hline
        \multirow{2}{*}{\textbf{inoperative}}& prob &&& 0.4 & 1.9 & 2.8 &&& 0.0 & 0.9 &&& 0.0 & 1.9 & 6.8 &&& 0.0 & 2.4 \\
        & reboot &&& 94.8 & 26.5 & 6.0 &&& --- & 88.2 &&& --- & 49.2 & 17.4 &&& -- & 36.9 \\ \hline
        \textbf{semi-} & prob &&& 10.6 & 21.2 & 39.2 &&& 5.0 & 17.1 &&& 5.1 & 20.7 & 44.1 &&& 5.0 & 37.6 \\
        \textbf{operational} & reboot &&& 63.1 & 17.1 & 2.2 &&& 100.0 & 18.8 &&& 100.0 & 21.8 & 3.6 &&& 100.0 & 8.7 \\ \hline \hline
    \end{tabular}
    \caption{Probability of a computer being inoperative or semi-operational (`prob'; in \%) and the corresponding conditional reboot probability (`reboot') for computers of different levels under different network topologies. \label{tab:reboot-probs}}
\end{table}

So far, our discussion has focused on a quantitative analysis of the suboptimality incurred and the runtimes as well as the memory consumed by different algorithms. Table~\ref{tab:reboot-probs} offers some insights into the qualitative structure of the determined policies. In particular, the table reveals the probabilities with which different computers are in an inoperative or semi-operational state in any given time period, as well as the conditional probabilities with which those computers are rebooted in such time periods, under the policies determined by our cutting plane approach. To this end, the table distinguishes between nodes of different levels as shown in Figure~\ref{fig:sysadmin:topologies}, ranging from 1 (highest) to 3 (lowest). The table omits the bidirectional ring topology as it does not contain nodes of different levels. In line with our intuition, the policies generated by our cutting plane approach reboot more promptly computers at lower levels, which in turn implies that those computers have a higher probability of being in a fully operational state.

\subsection{Robust Factored Markov Decision Processes}\label{sec:num_experiments:rfmdps}

We next consider a data-driven variant of the system administrator problem under the bidirectional ring topology with $N = 20$ computers. In this experiment, we assume that the transition kernel is unknown, and the decision maker only has access to a state-action history of length $T$ that has been generated by a historical policy. The historical policy coincides with the `priority' policy from the previous subsection, that is, it reboots computers in order of their state: inoperative computers are rebooted with highest priority, followed by semi-operational computers and, finally, fully operational computers. The historical policy breaks ties randomly.

\begin{figure}[tb]
    \centering
    \begin{tabular}{c@{}c}
        \includegraphics[width=8cm]{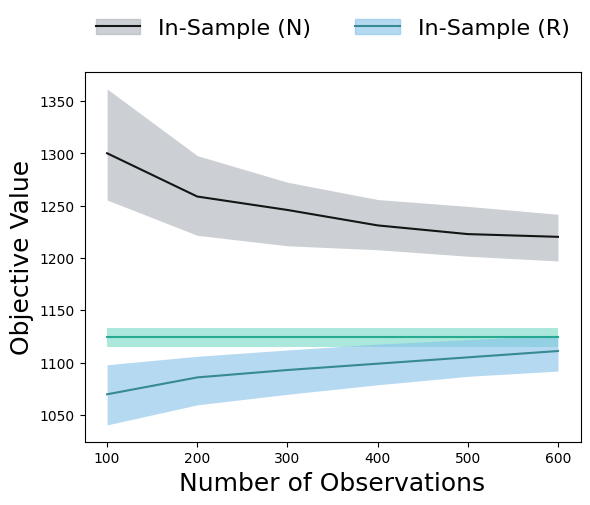} &
        \includegraphics[width=8cm]{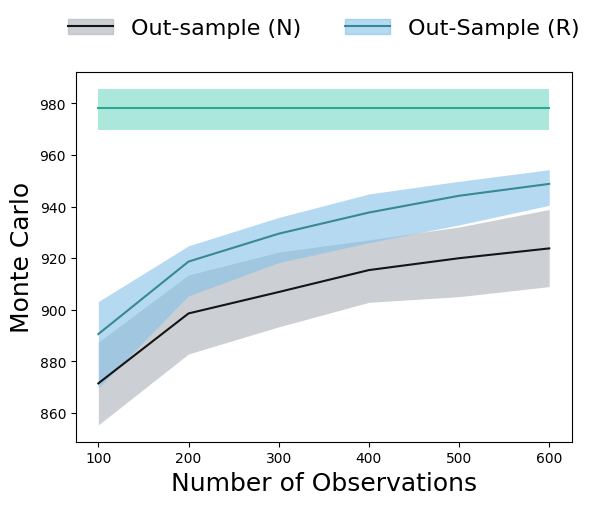}
    \end{tabular}
    \caption{In-sample (left) and out-of-sample (right) expected total rewards for the nominal (gray, `N') and robust (blue, `R') FMDP policies based on historical policies of different length (abscissa). The clairvoyant policy is shown in green. Bold lines represent median performances, and shaded regions correspond to the 25\%-75\% quantile ranges. \label{fig:robust-data-driven}}
\end{figure}

We again use Algorithm~\ref{alg:CP} to solve the nominal FMDP problem~\eqref{opt:fmdp:basis_fcts_lp_simplified}, where we now estimate the transition kernel from the state-action history using maximum likelihood estimation. We also use Algorithm~\ref{alg:CP-robust} to solve the RFMDP problem~\ref{opt:rfmdp:basis_fcts_lp_simplified} to stationarity. We use as marginal ambiguity sets the intersection of the $\infty$-norm ball of radius $\epsilon \in \{0, 0.03, 0.06, 0.09, 0.12, 0.15\}$ with the probability simplex $\Delta (\mathcal{S}_n)$. To this end, we split the the state-action history into a training history ($80\%$ of the transitions) and a validation history ($20\%$ of the transitions). We use maximum likelihood estimation to estimate a training and a validation transition kernel from the respective histories. We then use the training history to estimate the centers of our ambiguity sets, and we use the validation history to select the radii of our ambiguity sets. Figure~\ref{fig:robust-data-driven} compares the in-sample and out-of-sample performances of the nominal and the robust policy with the performance of the clairvoyant policy that has access to the true transition kernel. The results are based on 250 randomly generated problem instances, and the out-of-sample performances are approximated by 200 statistically independent Monte Carlo simulations over 150 transitions each. We observe that the nominal FMDP policies suffer from a high post-decision disappointment in that their in-sample expected total rewards substantially overestimate their out-of-sample counterparts. The post-decision disappointment of the robust FMDP policies, in contrast, tend to be negative: robust policies tend to `underpromise and overdeliver'. Note that the disappointment of robust policies can nevertheless become positive, which is owed to the fact that the radii of the ambiguity sets are chosen to maximize out-of-sample performance, rather than to provide statistically meaningful confidence regions. The figure also shows that the robust FMDP policies outperform their nominal counterparts out-of-sample by a statistically meaningful margin.

\subsection{Factored Multi-Armed Bandit Problems}\label{sec:num_experiments:bandits}

Our final experiment follows \citet{BZ22:dps_with_signals} and considers a multi-armed bandit problem with $M$ arms over an infinite time horizon. Each arm $m \in [M]$ transitions between $S_m$ different states. There are two actions associated with each arm $m \in [M]$, $a_m = 0$ (`do nothing') and $a_m = 1$ (`pull arm'), and we can pull up to $B$ arms in each time period. An arm remains in its current state under the `do nothing' action and does not generate any rewards, whereas it stochastically transitions into a new state and generates a random reward under the `pull arm' action. The transition probabilities and rewards of each arm depend on the arm's current state as well as a shared signal $z$ that evolves exogenously according to a Markov chain with $S_{M+1}$ states.

The problem can be naturally formulated as an FMDP with sub-state spaces $\mathcal{S}_m = \{ \mathrm{e}_1, \ldots, \mathrm{e}_{S_m} \} $, $m \in [M+1]$, where the first $M$ sub-states record the states of the arms and the last sub-state keeps track of the shared signal, respectively, and the action space $\mathcal{A} = \{ a \in \mathbb{B}^M : \sum_{m \in [M]} a_m \leq B \}$. In our experiments, we consider $M \in \{ 20, 50, 80 \}$ arms with a homogeneous number $S_m \in \{ 2, 3, 4 \}$ of states each. The shared signal, on the other hand, transitions between $S_{M+1} = 5$ states throughout our experiments. The initial distributions and transition probabilities of each arm $m$ and the shared signal $z$ are selected uniformly at random from the corresponding probability simplices, and the rewards of the `pull arm' actions are selected uniformly at random from the interval $[0, 1]$. We fix $B = \beta \cdot M$ where $\beta \in \{ 5\%, 10\%, 15\%, 20\% \}$ denotes the percentage of arms that can be pulled in each time period. Since a simple greedy policy turns out to perform well across this class of problems, we random select $\eta \cdot M$ arms, $\eta \in \{ 0\%, 10\%, 20\% \}$, as `special arms' whose rewards are zero under any state other than state $0$, and whose rewards in state $0$ under the `pull arm' action are inflated by a factor of $S_m$ to maintain a long term profitability that is comparable to that of the other arms. Throughout our experiments, we use the value functions $\nu_{mij} (s) = \mathbf{1} [ s_m = \mathrm{e}_i \text{ and } s_{M + 1} = \mathrm{e}_j ]$ for $m \in [M]$, $i \in [S_m]$ and $j \in [S_{M + 1}]$.

\begin{table}[tb]
    $\mspace{-25mu}$
    \renewcommand{\arraystretch}{1.1}
    \scriptsize
    \begin{tabular}{crp{0.7cm}p{0.6cm}p{0.6cm}p{0.6cm}p{0.4cm}cp{0.6cm}p{0.6cm}p{0.6cm}p{0.4cm}cp{0.6cm}p{0.6cm}p{0.6cm}p{0.7cm}}
        &&&\multicolumn{4}{c}{20 arms}&&\multicolumn{4}{c}{50 arms}&&\multicolumn{4}{c}{80 arms}\\
        \cline{4-7} \cline{9-12} \cline{14-17}        &&&5\%&10\%&15\%&20\%&&5\%&10\%&15\%&20\%&&5\%&10\%&15\%&20\%\\ \hline
        \multirow{9}{*}{\textbf{2 states}} &	\multirow{3}{*}{0\% S-arms}	& B\&Z &	7.56\% &	6.52\% & 	5.88\% &	5.62\% &	&	6.98\% &	5.76\% &	5.29\% &	4.90\% &	&	6.28\% &	5.48\% &	5.06\% &	4.66\% \\
        &&GRD&1.60\%& 	1.94\%& 	2.09\%& 	2.35\%&& 		\textbf{0.95\%} & 	\textbf{1.31\%} & 	\textbf{1.61\%} & 	\textbf{1.83\%} && 		\textbf{0.91\%} & 	\textbf{1.25\%} & 	\textbf{1.50\%}& 	\textbf{1.65\%} \\
        &&FMDP&	\textbf{1.13\%} & 	\textbf{1.15\%} & 	\textbf{1.33\%} & 	\textbf{1.50\%} && 		0.97\%& 	1.89\%& 	2.05\%& 	2.17\%&& 		1.89\%& 	2.40\%& 	2.45\%& 	2.47\%\\[1mm]
        & \multirow{3}{*}{10\% S-arms} &	B\&Z&	8.21\%& 	6.63\%& 	6.21\%& 	5.81\%&& 		7.07\%& 	5.91\%& 	5.52\%& 	5.03\%&& 		6.34\%& 	5.77\%& 	5.18\%& 	4.82\% \\ 
        &&GRD&	6.51\%& 	5.37\%& 	5.06\%& 	5.33\%&& 		5.99\%& 	4.75\%& 	4.55\%& 	4.73\%&& 		5.23\%& 	4.88\%& 	4.49\%& 	4.84\%\\
        &&FMDP&	\textbf{1.31\%} & 	\textbf{1.18\%} & 	\textbf{1.33\%} & 	\textbf{1.48\%} && 		\textbf{1.15\%} & 	\textbf{2.06\%} & 	\textbf{2.26\%} & 	\textbf{2.30\%} && 		\textbf{2.28\%} & 	\textbf{2.69\%} & 	\textbf{2.70\%} & 	\textbf{2.63\%} \\[1mm]
        &\multirow{3}{*}{20\% S-arms}&	B\&Z&	8.22\%& 	6.71\%& 	6.14\%& 	5.82\%&& 		7.39\%& 	6.05\%& 	5.78\%& 	5.37\%&& 		6.51\%& 	5.76\%& 	5.19\%& 	5.04\%\\ 
        && GRD&	8.23\%& 	8.17\%& 	8.40\%& 	8.83\%&& 		8.70\%& 	8.05\%& 	7.67\%& 	7.90\%&& 		8.83\%& 	8.12\%& 	7.44\%& 	7.81\%\\ 
        &&FMDP&	\textbf{1.65\%} & 	\textbf{1.20\%} & 	\textbf{1.22\%} & 	\textbf{1.53\%} && 		\textbf{1.73\%} & 	\textbf{2.30\%} & 	\textbf{2.34\%} & 	\textbf{2.38\%} && 		\textbf{2.90\%} & 	\textbf{3.04\%} & 	\textbf{3.00\%} & 	\textbf{3.00\%} \\[1mm] \hline \\[-2mm]
        \multirow{9}{*}{\textbf{3 states}} &	\multirow{3}{*}{0\% S-arms} &	B\&Z&	10.62\%& 	9.22\%& 	8.46\%& 	7.89\%&& 		9.40\%& 	8.10\%& 	7.60\%& 	6.77\%&& 		8.49\%& 	7.56\%& 	6.90\%& 	6.55\%\\ 
        &&		GRD&	2.25\%& 	2.74\%& 	3.12\%& 	3.44\%&& 		1.42\%& 	\textbf{2.04\%} & 	\textbf{2.40\%} & 	2.73\%&& 		\textbf{1.38\%} &	\textbf{1.89\%} & 	\textbf{2.22\%} & 	\textbf{2.59\%} \\ 
        &&FMDP&	\textbf{1.73\%} & 	\textbf{1.82\%} & 	\textbf{1.97\%} & 	\textbf{2.02\%} && 		\textbf{1.37\%} & 	2.23\%& 	2.41\%& 	\textbf{2.48\%} && 		2.22\%& 	2.72\%& 	2.83\%& 	3.06\%\\[1mm]
        &	\multirow{3}{*}{10\% S-arms}&	B\&Z &	11.51\%& 	9.62\%& 	8.73\%& 	8.02\%&& 		10.65\%& 	8.55\%& 	7.81\%& 	7.18\%&& 		9.17\%& 	8.02\%& 	7.37\%& 	6.83\%\\ 
        &&		GRD&	7.13\%& 	7.73\%& 	7.03\%& 	7.39\%& &		9.66\%& 	7.68\%& 	7.04\%& 	6.77\%&& 		9.33\%& 	6.77\%& 	6.76\%& 	6.77\%\\ 
        &&		FMDP&	\textbf{1.84\%} & 	\textbf{1.68\%} & 	\textbf{1.79\%} & 	\textbf{1.89\%} &&	\textbf{1.68\%} & 	\textbf{2.17\%} & 	\textbf{2.38\%} & 	\textbf{2.52\%} && 		\textbf{2.36\%} & 	\textbf{2.72\%} & 	\textbf{3.02\%} & 	\textbf{3.11\%} \\[1mm]
        &	\multirow{3}{*}{20\% S-arms} &	B\&Z&	12.10\%& 	10.13\%& 	9.16\%& 	8.67\%&& 		10.85\%& 	9.02\%& 	8.25\%& 	7.49\%&& 		9.42\%& 	8.32\%& 	7.72\%& 	7.10\%\\ 
        &&GRD&	13.06\%& 	12.01\%& 	11.69\%& 	11.55\%&& 		14.12\%& 	11.36\%& 	10.71\%& 	10.61\%&& 		13.91\%& 	11.67\%& 	11.19\%& 	11.20\%\\ 
        &&		FMDP&	\textbf{2.74\%}& 	\textbf{1.76\%}& 	\textbf{1.71\%}& 	\textbf{1.90\%}&& 		\textbf{2.44\%}& 	\textbf{2.20\%}& 	\textbf{2.38\%}& 	\textbf{2.53\%}&& 		\textbf{2.92\%}& 	\textbf{2.87\%}& 	\textbf{3.00\%}& 	\textbf{3.08\%} \\[1mm] \hline \\[-2mm]
        \multirow{9}{*}{\textbf{4 states}} &	\multirow{3}{*}{0\% S-arms} &	B\&Z &	12.29\%& 	10.66\%& 	9.62\%& 	8.72\%&& 		10.85\%& 	9.12\%& 	8.58\%& 	7.88\%&& 		9.59\%& 	8.62\%& 	7.97\%& 	7.51\%\\ 
        &&		GRD	&2.71\%& 	3.29\%& 	3.66\%& 	3.94\%&& 		\textbf{1.68\%}& 	2.46\%& 	2.87\%& 	3.32\%&& 		\textbf{1.60\%}& 	\textbf{2.25\%}& 	\textbf{2.78\%}& 	\textbf{3.12\%}\\ 
        &&		FMDP&	\textbf{2.24\%}& 	\textbf{2.21\%}& 	\textbf{2.28\%}& 	\textbf{2.29\%}&& 		1.76\%& 	\textbf{2.37\%}& 	\textbf{2.57\%}& 	\textbf{2.70\%}&& 		2.45\%& 	2.92\%& 	3.05\%& 	3.24\%\\[1mm] 
        & \multirow{3}{*}{10\% S-arms}	&B\&Z	&13.21\%& 	11.57\%& 	10.27\%& 	9.73\%&& 		11.86\%& 	9.88\%& 	9.18\%& 	8.28\%&& 		10.61\%& 	9.23\%& 	8.47\%& 	7.91\%\\ 
        &&	GRD&	10.94\%& 	9.70\%& 	9.46\%& 	8.96\%&& 		12.71\%& 	8.83\%& 	8.54\%& 	8.02\%&& 		11.66\%& 	8.13\%& 	7.83\%& 	8.13\%\\ 
        &&		FMDP&	\textbf{2.66\%}& 	\textbf{2.13\%}& 	\textbf{2.08\%}& 	\textbf{2.21\%}&& 		\textbf{2.41\%}& 	\textbf{2.33\%}& 	\textbf{2.50\%}& 	\textbf{2.62\%}&& 		\textbf{2.58\%}& 	\textbf{2.80\%}& 	\textbf{3.02\%}& 	\textbf{3.15\%}\\[1mm] 
        & \multirow{3}{*}{20\% S-arms}&	B\&Z &	15.06\%& 	11.90\%& 	10.71\%& 	9.97\%&& 		12.77\%& 	9.86\%& 	9.48\%& 	8.82\%&& 		11.16\%& 	9.48\%& 	8.98\%& 	8.40\%\\
        &&		GRD&	15.58\%& 	14.14\%& 	12.19\%& 	13.22\%&& 		17.53\%& 	13.79\%& 	12.79\%& 	12.89\%&& 		17.32\%& 	13.85\%& 	12.31\%& 	12.65\%\\ 
        &&		FMDP&	\textbf{3.65\%}&  	\textbf{2.28\%}&  	\textbf{2.14\%}&  	\textbf{2.12\%}& & 		\textbf{3.23\%}&  	\textbf{2.29\%}&  	\textbf{2.48\%}&  	\textbf{2.59\%}& & 		\textbf{3.24\%}&  	\textbf{2.85\%}&  	\textbf{3.09\%}&  	\textbf{3.21\%}\\ 
        \hline \hline
    \end{tabular}
    \caption{Optimality gaps of our FMDP approach (`FMDP'), the dynamic fluid policy (`B\&Z') and the greedy policy (`GRD') across different sub-state space dimensions $S_m$ and percentages $\eta$ of special arms (rows) as well as numbers of arms $M$ and percentages $\beta$ of arms to be pulled (columns). All gaps are averages over 100 random problem instances. The best performance for each instance class is highlighted in bold print. \label{tab:my_label}}
\end{table}

Similar to our earlier experiments in Section~\ref{sec:num_experiments:fmdps}, we solve our FMDP problem~\eqref{opt:fmdp:basis_fcts_lp_simplified} using Algorithm~\ref{alg:CP}. The factored multi-armed bandit problem also allows for a representation as a weakly coupled MDP, and we thus compare our FMDP approach with the dynamic fluid policy of \cite{BZ22:dps_with_signals}. Since this policy is restricted to  finite horizon problems, we fix its time horizon to $150$ periods and discount the rewards geometrically. We also consider a greedy policy that pulls the $B$ arms with largest one-step rewards in any particular time period. Table~\ref{tab:my_label} compares the out-of-sample rewards generated by these three policies in terms of relative gap to the (typically unachievable) upper reward bound provided by the fluid relaxation of \citet{BZ22:dps_with_signals}. The table shows that our FMDP approach dominates both benchmark approaches over the considered class of instances. That said, we expect the dynamic fluid policy to outperform our FMDP strategy for large problem sizes since it enjoys asymptotic optimality guarantees.

\newpage
{
    \singlespacing
    \small
    \bibliographystyle{plainnat}
    \bibliography{bibliography}
}

\newpage

\section*{E-Companion A: Formal Construction of the Lift}

\renewcommand{\thepage}{EC-\arabic{page}}
\setcounter{page}{1}

For a factored RFMDP $(\mathcal{S}, \mathcal{A}, \phi, q, \mathcal{P}, r, \gamma)$ with any fixed state $s^0 \in \mathcal{S}$ and action $a^0 \in \mathcal{A}$, we define the corresponding non-factored lifted RFMDP $(\hat{\mathcal{S}}, \hat{\mathcal{A}}, \hat{\phi}, \hat{q}, \hat{\mathcal{P}}, \hat{r}, \hat{\gamma})$ as follows. The lifted state space is $\hat{\mathcal{S}} = \mathcal{S}^2 \times \mathcal{A} \times [\lceil \log_2 N \rceil]$, where each state $\hat{s} = (u, m, a, n)$ comprises an update block $u \in \mathcal{S}$, a memory block $m \in \mathcal{S}$, an action block $a \in \mathcal{A}$ and a timer block $n \in [N]$ (represented in binary notation to adhere to our Definitions~\ref{def:fmdp} and~\ref{def:fb-fmdp} of FMDPs). The lifted action space is $\hat{\mathcal{A}} = \mathcal{A}$. The lifted initial state distribution is defined as $q (\hat{s}) = q (u)$ for any $\hat{s} = (u, s^0, a^0, 1)$ with $u \in \mathcal{S}$, whereas $q (\hat{s}) = 0$ otherwise. For $\hat{s} = (u, m, a, n)$, we define the lifted rewards
\begin{equation*}
    \hat{r} (\hat{\phi} (\hat{s}, \hat{a})) = \mathbf{1} [(m, a, n) = (s^0, a^0, 1)] \cdot r_j (\phi (u, \hat{a})).    
\end{equation*}
The lifted discount factor is $\hat{\gamma} = \gamma^{1/N}$.

We next define the marginal ambiguity sets of the lift. For $(\hat{s}, \hat{a}) \in \hat{\mathcal{S}} \times \hat{\mathcal{A}}$ and $\hat{n} \in [N]$, that is, the ambiguity sets relating to the update block $u$ in state $\hat{s} = (u, m, a, n)$, we fix
\begin{equation*}
    \hat{\mathcal{P}}^{\hat{n}} (\phi (\hat{s}, \hat{a})) =
    \begin{cases}
        \mathcal{P}^n ( \phi (u, \hat{a}) ) & \text{if } \hat{n} = n = 1, \\
        \mathcal{P}^n ( \phi (m, a) ) & \text{if } \hat{n} = n \neq 1 \\
        \{ p_{\hat{n}} \} & \text{otherwise},
    \end{cases}
\end{equation*}
where $p_{\hat{n}} (u'_{\hat{n}} | \phi(\hat{s}, \hat{a})) = \mathbf{1} [u'_{\hat{n}} = u_{\hat{n}}]$ for $s' = (u', m', a', n')$. In other words, the currently active update block $u_{\hat{n}}$ is updated according to some distribution from $\mathcal{P}^{\hat{n}}$, accounting for the state and action recorded at time $n = 1$ and stored in $(u, \hat{a})$ (if $\hat{n} = 1$) or $(m, a)$ (if $\hat{n} > 1$). For $\hat{n} \in [2N] \setminus [N]$, that is, the ambiguity sets relating to the memory block $m$ in state $\hat{s} = (u, m, a, n)$, we fix $\hat{\mathcal{P}}^{\hat{n}} (\phi (\hat{s}, \hat{a})) = \{ p_{\hat{n}} \}$ with
\begin{equation*}
    p_{\hat{n}} (m'_{\hat{n} - N} | \phi(\hat{s}, \hat{a})) = \mathbf{1} \left[
    \begin{array}{l}
        \displaystyle \mspace{19mu} (m'_{\hat{n} - N}, n) = (u_{\hat{n} - N}, 1) \\
        \displaystyle \vee \; (m'_{\hat{n} - N}, n) \in \{ m_{\hat{n} - N} \} \times ([N - 1] \setminus \{ 1 \}) \\
        \displaystyle \vee \; (m'_{\hat{n} - N}, n) = (s^0, N)
    \end{array}
    \right]
\end{equation*}
for $s' = (u', m', a', n')$. Thus, the content of the update block $u$ is stored in $m'$ at time $n = 1$, $m'$ is reset to $s^0$ at time $n = N$, and $m'$ is left unchanged at other times. For $\hat{n} = 2N + 1$, that is, the ambiguity set related to the action block $a$ in state $\hat{s} = (u, m, a, n)$, we fix $\hat{\mathcal{P}}^{\hat{n}} (\phi (\hat{s}, \hat{a})) = \{ p_{\hat{n}} \}$ with
\begin{equation*}
    p_{\hat{n}} (a' | \phi(\hat{s}, \hat{a})) = \mathbf{1}
    \left[
    \begin{array}{l}
        \displaystyle \mspace{19mu} (a', n) = (\hat{a}, 1) \\
        \displaystyle \vee \; (a', n) \in \{ a \} \times ([N - 1] \setminus \{ 1 \}) \\
        \displaystyle \vee \; (a', n) = (a^0, N)
        \end{array}
    \right]
\end{equation*}
for $s' = (u', m', a', n')$, that is, the new action $\hat{a}$ is stored in $a'$ at time $n = 1$, $a'$ is reset to $a^0$ at time $n = N$, and $a'$ is left unchanged at other times. For $\hat{n} = 2N + 2$, that is, the ambiguity set related to the timer block $n$ in state $\hat{s} = (u, m, a, n)$, finally, we fix $\hat{\mathcal{P}}^{\hat{n}} (\phi (\hat{s}, \hat{a})) = \{ p_{\hat{n}} \}$ with
\begin{equation*}
    p_{\hat{n}} (\text{bin} (n') | \phi(\hat{s}, \hat{a})) = \mathbf{1} \big[ n' = n \oplus 1 \big],
\end{equation*}
where $\text{bin} (n')$ refers to the binary representation of $n'$ and $n \oplus 1 = n + 1$ for $n < N$ and $N \oplus 1 = 1$. Thus, the timer $n'$ is incremented by one unit in each transition and reset to $1$ every $N$ transitions.

The above transition and reward structure can be implemented with a lifted feature mapping $\hat{\phi}$ that contains \emph{(i)} the lifted state $\hat{s}$; \emph{(ii)} the feature mappings $\phi (u, \hat{a})$ and $\phi (m, a)$ corresponding to the non-lifted state-action pairs at times $n = 1$ and $n > 1$, respectively; as well as the \emph{(iii)} $\mathcal{O} (| \mathcal{A} | + \sum_{n \in [N]} | \mathcal{S}_n |)$ indicator functions used by the rewards and marginal transition probabilities. With this mapping $\hat{\phi}$, any of the above singleton-valued ambiguity set mappings can be expressed as scope-1 mappings. If the $\hat{n}$-th ambiguity set mapping $\mathcal{P}^{\hat{n}}$ of the non-lifted factored RFMDP, $\hat{n} \in [N]$, has scope $\sigma$, then the $\hat{n}$-th ambiguity set mapping $\hat{\mathcal{P}}^{\hat{n}}$ of the lifted non-factored RFMDP has a scope of at most $\sigma + \lceil \log_2 N \rceil + \max_{n \in [N]} S_n$ due to its additional dependence on the timer $n$ as well as one component of the update block $u$ in the current state $\hat{s} = (u, m, a, n)$. Likewise, if the $j$-th reward component of the non-lifted factored RFMDP has scope $\sigma$, then the $j$-th reward component of the lifted non-factored RFMDP has a scope of $\sigma + 1$ due to its additional dependence on the indicator function $\mathbf{1} [(m, a, n) = (s^0, a^0, 1)]$.

\newpage

\section*{E-Companion B: Proofs}

% \noindent \textbf{Proof of Observation~\ref{obs:scope_of_trans_probs}.}
%     Let $Q:\mathcal{S}\times \mathcal{A} \to \Delta (\mathcal{S}_n)$ be the map defined by $ Q(s,a) = p_n (\cdot\, | \, s, a)$. By assumption, there are $C$ possible images of $Q$, denoted by $q^1, \dots, q^C$.
%     Let $L = \ceil[\big]{\log C} $. For $\ell = 1,\dots, L$, define the binary variable
%     \[ b_\ell (s,a) = \text{ the $\ell$-th digit of the binary representation of the number } c(s,a), \]
%     where $c(s,a)$ is defined by $Q(s,a) = q^{c(s,a)}$. Taking $b = (b_1,\dots,b_L)^\top$ as the feature map, we can write
%     \[ p_n (s'_n\, | \, s, a) = p_n (s'_n\, | \, b(s, a)). \]
%     By Lemma~\ref{lem:shiny_new_lemma}, each $b_\ell (s,a)$ admits a mixed-integer linear representation. Thus, the scope of $p_n$ is $L = \ceil[\big]{\log C}$.
% \qed

% ~\\[-8mm]

\noindent \textbf{Proof of Theorem~\ref{thm:complexity_representing_csi}.}    
    In view of statement~\emph{(i)}, consider first a representation of the function $f : \mathbb{B}^I \rightarrow \mathbb{R}$ as a rule system. Concretely, fix a representation of the form
    \begin{equation*}
        \bigwedge_{(i, v) \in \mathfrak{A}_r} (x_i = v) \;\; \Rightarrow \;\; f(x) = c_r \qquad \forall r \in \mathfrak{R},
    \end{equation*}
    where $\mathfrak{R}$ is the set of rules to determine $f$,  $\mathfrak{A}_r$ is the set of antecedents (or premises) for rule $r$ and $c_r$ is the consequent (or conclusion) of $r$. This rule-based representation is equivalent to the MILP-based feature representation under which $f(x) = c_r$ if and only if $\phi_r (x) = 1$, $r \in \mathfrak{R}$, where the MILP-based features $\phi_r (x) = \mathbf{1} [\bigwedge_{(i, v) \in \mathfrak{A}_r} (x_i = v)]$ have the polynomial-size representation
    \begin{equation*}
        \left.
        \begin{array}{ll}
            \displaystyle \phi_r (x) \leq 1 - x_i \;\; \forall (i, 0) \in \mathfrak{A}_r; &
            \displaystyle \phi_r (x) \leq x_i \;\; \forall (i, 1) \in \mathfrak{A}_r \\
            \multicolumn{2}{l}{\displaystyle \phi_r (x) \geq \sum_{(i, 0) \in \mathfrak{A}_r} (1 - x_i) + \sum_{(i, 1) \in \mathfrak{A}_r} x_i - \left[ | \mathfrak{A}_r | - 1 \right]}
        \end{array}
        \right\}
        \;\; \forall r \in \mathfrak{R}.
    \end{equation*}
    This concludes statement \emph{(i)} for rule systems. Since every decision tree can be represented as a rule system with one rule for each path from the root node to one of the leaf nodes, our proof for rule systems immediately applies to decision tree representations as well.

    To see that statement \emph{(i)} also holds for ADDs, we first observe that the empty function $f : \emptyset \rightarrow \mathbb{R}$ satisfying $f() = c$ has a trivial MILP-based feature representation. Assume now that the two functions $f', f'' : \mathbb{B}^{I - i}$ are ADDs of binary vectors $(x_{i+1}, \ldots, x_I)$ with polynomial-sized MILP-based feature representations $f' (x_{i+1}, \ldots, x_I) = c'_r$ if $\phi'_r (x_{i+1}, \ldots, x_I) = 1$, $r \in \mathfrak{R}'$, and $f'' (x_{i+1}, \ldots, x_I) = c''_r$ if $\phi''_r (x_{i+1}, \ldots, x_I) = 1$, $r \in \mathfrak{R}''$. This implies that there are polyhedra $\mathcal{F}', \mathcal{F}''$ as well as auxiliary vectors $\zeta' \in \mathbb{R}^{F'_l} \times \mathbb{B}^{F'_b}$ and $\zeta'' \in \mathbb{R}^{F''_l} \times \mathbb{B}^{F''_b}$, all of polynomial size, such that $\phi' (x_{i+1}, \ldots, x_I) = \varphi \in \mathbb{B}^{| \mathfrak{R}' |}$ precisely when $(\varphi, \zeta; x_{i+1}, \ldots, x_I) \in \mathcal{F}'$ and $\phi'' (x_{i+1}, \ldots, x_I) = \varphi  \in \mathbb{B}^{| \mathfrak{R}'' |}$ precisely when $(\varphi, \zeta; x_{i+1}, \ldots, x_I) \in \mathcal{F}''$, respectively. We show that in this case, the function $f : \mathbb{B}^{I - i + 1}$ satisfying $f (x_i, \ldots, x_I) = x_i \cdot f' (x_{i+1}, \ldots, x_I) + (1 - x_i) \cdot f'' (x_{i+1}, \ldots, x_I)$ has an MILP-based feature representation $f (x_i, \ldots, x_I) = c_r$ if $\phi_r (x_i, \ldots, x_I) = 1$, $r \in \mathfrak{R}$, of polynomial size. Indeed, fix $\mathfrak{R} = \{ (1, r') \, : \, r' \in \mathfrak{R}' \} \cup \{ (0, r'') \, : \, r'' \in \mathfrak{R}'' \}$,
    $c_r = c'_{r'}$ for $r = (1, r') \in \mathfrak{R}$, $r' \in \mathfrak{R}'$, and $c_r = c''_{r''}$ for $r = (0, r'') \in \mathfrak{R}$, $r'' \in \mathfrak{R}''$. The claim follows if we can show that the feature map $\phi$ satisfying
    \begin{equation*}
        \begin{array}{rll}
            & \displaystyle \phi_r (x_i, \ldots, x_I) = \mathbf{1} [ x_i = 1 \text{ and } \phi'_{r'} (x_{i+1}, \ldots, x_I) = 1 ] & \displaystyle \text{for } r = (1, r') \in \mathfrak{R}, r' \in \mathfrak{R}' \\
            \text{and} & \displaystyle \phi_r (x_i, \ldots, x_I) = \mathbf{1} [ x_i = 0 \text{ and } \phi''_{r''} (x_{i+1}, \ldots, x_I) = 1 ] & \displaystyle \text{for } r = (0, r'') \in \mathfrak{R}, r'' \in \mathfrak{R}''
        \end{array}
    \end{equation*}
    has an MILP representation of polynomial size. This is the case since we have
    \begin{equation*}
        \phi (x_i, \ldots, x_I) = \varphi \in \mathbb{B}^{| \mathfrak{R} |}
        \quad \Longleftrightarrow \quad
        \left[
        \begin{array}{ll}
            \displaystyle \varphi \geq \varphi' - (1 - x_i) \cdot \mathrm{e},
            & \displaystyle \varphi \leq \varphi' + (1 - x_i) \cdot \mathrm{e} \\
            \displaystyle \varphi \geq \varphi'' - x_i \cdot \mathrm{e}, 
            & \displaystyle \varphi \leq \varphi'' + x_i \cdot \mathrm{e} \\
            \multicolumn{2}{l}{\displaystyle (\varphi', \zeta'; x_{i+1}, \ldots, x_I) \in \mathcal{F}_{r'}} \\
            \multicolumn{2}{l}{\displaystyle (\varphi'', \zeta''; x_{i+1}, \ldots, x_I) \in \mathcal{F}_{r''}}
        \end{array}
        \right].
    \end{equation*}
    
    As for statement \emph{(ii)}, consider first the parity function $f (x_1, \ldots, x_I) = 1$ if $\sum_{i=1}^I x_i$ is even; $= 0$ otherwise. This function cannot be computed without considering the value $x_i$ of each of its inputs $i \in [I]$, and thus any representation of this function as a tree or rule system necessarily scales exponentially in $I$. In contrast, we have $f (x_1, \ldots, x_I) = 1$ if $\phi (x_1, \ldots, x_I) = 1$; $= 0$ otherwise for the feature map $\phi (x_1, \ldots, x_I) = \mathbf{1} [\sum_{i=1}^I x_i \text{ is even}]$, which has the MILP representation
    \begin{equation*}
        \phi (x_1, \ldots, x_I) = \varphi
        \quad \Longleftrightarrow \quad
        \left[
        \begin{array}{l}
            \displaystyle \varphi = 1 - z_I, \quad z_1 = x_1 \\
            \displaystyle \mspace{-10mu}
            \left.
            \begin{array}{l}
                \displaystyle z_{i+1} = z_i + x_{i+1} - 2 y_{i+1} \\
                \displaystyle y_{i+1} \leq z_i, \quad y_{i+1} \leq z_{i+1} \\
                \displaystyle y_{i+1} \geq z_i + z_{i+1} - 1
            \end{array}
            \right\} \;\; \forall i \in [I - 1]
        \end{array}
        \right],
    \end{equation*}
    $\varphi \in \mathbb{B}$, that is of polynomial size. Note that in this formulation, we have $y_{i+1} = \mathbf{1} [ z_i = 1 \wedge z_{i+1} = 1 ]$ for all $i \in [I - 1]$ and $z_i = \mathbf{1} [\sum_{i'=1}^i x_{i'} \text{ is odd}]$ for all $i \in [I]$.
    
    The argument from the previous paragraph does not extend to ADDs; in fact, ADDs can also represent the parity function efficiently. To see that statement \emph{(ii)} nevertheless extends to ADDs, consider the function $f (a_1, \ldots, a_I; b_1, \ldots, b_I)$ that computes one of the $2I$ bits of the result of multiplying the binary numbers $(a_1, \ldots, a_I)$ and $(b_1, \ldots, b_I)$. It follows from \cite{B91:complexity} that ADDs cannot represent this function in polynomial size. In contrast, we show next that this function can be represented efficiently via MILP-based features. To this end, consider the following feature map, which computes the entire result of the multiplication:
    \begin{equation*}
        \phi (a_1, \ldots, a_I; b_1, \ldots, b_I) = \varphi
        \;\; \Longleftrightarrow \;\;
        \left[
        \begin{array}{l}
            \displaystyle \sum_{i \in [2I]} 2^{i-1} \cdot \varphi_i
            \; = \;
            \sum_{i, j \in [I]} 2^{i + j - 2} \cdot \zeta_{ij} \\
            \displaystyle \zeta_{ij} \leq a_i, \;\; \zeta_{ij} \leq b_j, \;\; \zeta_{ij} \geq a_i + b_j - 1 \quad \forall i, j \in [I]
        \end{array}
        \right]
    \end{equation*}
    To see that the above representation is correct, observe first that for all $i, j \in [I]$ we have $\zeta_{ij} = a_i \cdot b_j$. The right-hand side of the first constraint therefore evaluates to
    \begin{equation*}
        \sum_{i, j \in [I]} 2^{i + j - 2} \cdot a_i \cdot b_j
        \;\; = \;\;
        \left( \sum_{i \in [I]} 2^{i-1} a_i \right) \left( \sum_{j \in [I]} 2^{j-1} b_j \right),    
    \end{equation*}
    where $\sum_{i \in [I]} 2^{i-1} \cdot a_i$ and $\sum_{j \in [I]} 2^{j-1} \cdot b_j$ are the numbers represented by the binary vectors $(a_1, \ldots, a_I)$ and $(b_1, \ldots, b_I)$, respectively. The left-hand side of the first constraint thus ensures that the number represented by the binary vector $(\varphi_1, \ldots, \varphi_{2I})$ equals the product of $(a_1, \ldots, a_I)$ and $(b_1, \ldots, b_I)$, as desired. Our constraints defining the feature map $\phi$ involve coefficients that are exponential in the length $2I$ of its input; this can be avoided by emulating a binary multiplier through MILP constraints. We omit the details of this tedious but otherwise straightforward representation.
\qed

~\\[-8mm]

\noindent \textbf{Proof of Theorem~\ref{thm:generic_and_nasty_features}.}
    In view of the first assertion, fix any feature map $\phi: \mathcal{S} \times \mathcal{A} \rightarrow \mathbb{B}^F$ and define $\Phi_l = \{ (s, a, \varphi_l) \in \mathcal{S} \times \mathcal{A} \times \mathbb{B} \, : \, \phi_l (s, a) = \varphi_l \}$ as the graph of the $l$-th component of $\phi$. We then note that for all $(s, a) \in \mathcal{S} \times \mathcal{A}$ and $l \in [F]$, we have that
    \begin{equation*}
        \phi_l (s, a) = \varphi_l
        \quad \Longleftrightarrow \quad
        (s, a, \varphi_l) \in \Phi_l
        \quad \Longleftrightarrow \quad
        (s, a, \varphi_l) \neq (s', a', \varphi'_l) \;\; \forall (s', a', \varphi'_l) \in [ \mathcal{S} \times \mathcal{A} \times \mathbb{B} ] \setminus \Phi_l.
    \end{equation*}
    The latter condition holds precisely when
    \begin{equation*}
        \mspace{-25mu}
        \begin{array}{rl@{\quad}l}
            & \displaystyle
            \mspace{2mu} \left[ s \neq s' \; \vee \; a \neq a' \; \vee \; \varphi_l \neq \varphi'_l \right]
            & \displaystyle \forall (s', a', \varphi'_l) \in [ \mathcal{S} \times \mathcal{A} \times \mathbb{B} ] \setminus \Phi_l \\[0mm]
            \Longleftrightarrow \quad & \displaystyle
            \left[ \bigvee_{n \in [N]} \bigvee_{i \in [S_n]} s_{ni} \neq s'_{ni} \; \vee \; \bigvee_{m \in [A]} a_m \neq a'_m \; \vee \; \varphi_l \neq \varphi'_l \right]
            & \displaystyle \forall (s', a', \varphi'_l) \in [ \mathcal{S} \times \mathcal{A} \times \mathbb{B} ] \setminus \Phi_l \\[5mm]
            \Longleftrightarrow \quad & \displaystyle
            \left[ \sum_{n \in [N]} \sum_{i \in [S_n]} \frac{s'_{ni} - s_{ni}}{2 s'_{ni} - 1} \geq 1 \; \vee \; \sum_{m \in [A]} \frac{a'_m - a_m}{2 a'_m - 1} \geq 1 \; \vee \; \frac{\varphi'_l - \varphi_l}{2 \varphi'_l - 1} \geq 1 \right]
            & \displaystyle \forall (s', a', \varphi'_l) \in [ \mathcal{S} \times \mathcal{A} \times \mathbb{B} ] \setminus \Phi_l \\[5mm]
            \Longleftrightarrow \quad & \displaystyle
            \sum_{n \in [N]} \sum_{i \in [S_n]} \frac{s'_{ni} - s_{ni}}{2 s'_{ni} - 1} + \sum_{m \in [A]} \frac{a'_m - a_m}{2 a'_m - 1} \; + \; \frac{\varphi'_l - \varphi_l}{2 \varphi'_l - 1} \geq 1
            & \displaystyle \forall (s', a', \varphi'_l) \in [ \mathcal{S} \times \mathcal{A} \times \mathbb{B} ] \setminus \Phi_l.
        \end{array}
    \end{equation*}
    In summary, for each $l \in [F]$ the constraint $\phi_l (s,a) = \varphi_l$ can be represented by a system of at most $2 \cdot | \mathcal{S} \times \mathcal{A}|$ linear constraints, which proves the first assertion.
    
    As for the second assertion, consider a class of FMDP instances $(\mathcal{S}, \mathcal{A}, q, p, r, \gamma)$ whose state spaces $\mathcal{S}$ encode 3SAT instances in binary representation, while the definition of their action spaces $\mathcal{A}$, their initial state and transition probabilities $q$ and $p$, their rewards $r$ as well as their discount factors $\gamma$ are not relevant for our construction and are thus omitted. We choose the class of features $\phi$ that evaluate the solvability of these 3SAT instances, that is, $\phi$ evaluates to $1$ if and only if the answer to the 3SAT instance associated with the FMDP instance is affirmative. Since 3SAT is known to be strongly NP-hard \citep{GJ79:ComputersIntractability}, the evaluation of $\phi$ must be strongly NP-hard. Thus, unless P = NP, any system of linear inequalities evaluating $\phi$ must have a number of variables or constraints that is exponential in the bit length of $\mathcal{S}$ (since linear programs can be solved in polynomial time), and any system of mixed-binary linear constraints computing $\phi$ must have exponentially many continuous variables or constraints, or it must contain a number of binary variables that grows at least logarithmically in $\mathcal{S}$ (since we can enumerate all possible values of the binary variables and solve linear programs in the remaining variables and constraints).
\qed

~\\[-8mm]

\noindent \textbf{Proof of Proposition~\ref{prop:example_feature_mapping}.}
    In view of the first assertion, note that
    \begin{equation*}
        \phi (s, a)
        \; = \;
        \mathbf{1} [s_n \in \mathcal{S}'_n]
        \; = \;
        \sum_{s'_n \in \mathcal{S}'_n} \mathbf{1} [s_n = s'_n],
    \end{equation*}
    and that each indicator function $\mathbf{1} [s_n = s'_n]$ can be expressed by a variable $\zeta_{s'_n} \in \mathbb{R}_+$ satisfying
    \begin{equation*}
        1 + \sum_{j \in [S_n]} \frac{s_{n,j} - s'_{n,j}}{2s'_{n,j} - 1}
        \;\; \leq \;\;
        \zeta_{s'_n}
        \;\; \leq \;\;
        \min_{j \in [S_n]} \left\{ \frac{s'_{n,j} + s_{n,j} - 1}{2s'_{n,j} - 1} \right\}
        \qquad \forall s'_n \in \mathcal{S}'_n.
    \end{equation*}
    This representation has $F_l = | \mathcal{S}'_n |$ auxiliary variables and $F_c = 1 + | \mathcal{S}'_n | + S_n \cdot | \mathcal{S}'_n |$ constraints.

    As for the second assertion, note that
    \begin{equation*}
        \phi (s, a)
        \; = \;
        \mathbf{1} [s_n \in \mathcal{S}'_n \;\; \forall n \in \mathcal{N}]
        \quad \Longleftrightarrow \quad
        \left( \sum_{n \in \mathcal{N}} \mathbf{1} [s_n \in \mathcal{S}'_n] \right) - (| \mathcal{N} | - 1)
        \; \leq \;
        \phi (s, a)
        \; \leq \;
        \min_{n \in \mathcal{N}} \, \mathbf{1} [s_n \in \mathcal{S}'_n],
    \end{equation*}
    and that each indicator function $\mathbf{1} [s_n \in \mathcal{S}'_n]$ can be represented using $| \mathcal{S}'_n |$ auxiliary variables and $(1 + S_n) \cdot | \mathcal{S}'_n |$ constraints (\emph{cf.}~Assertion 1). The overall representation thus has $F_l = \sum_{n \in \mathcal{N}} | \mathcal{S}'_n |$ auxiliary variables and $F_c = \sum_{n \in \mathcal{N}} (1 + S_n) \cdot | \mathcal{S}'_n | + 1 + | \mathcal{N} |$ constraints as claimed.

    In view of the third assertion, note that
    \begin{equation*}
        \phi (s, a)
        \; = \;
        \mathbf{1} [\exists n \in \mathcal{N} \, : \, s_n \in \mathcal{S}'_n]
        \; = \;
        1 - \mathbf{1} [s_n \in \mathcal{S}_n \setminus \mathcal{S}'_n \;\; \forall n \in \mathcal{N}],
    \end{equation*}
    and we can thus reuse the representation from Assertion~2.

    As for the fourth assertion, note that
    \begin{equation*}
        \phi (s, a)
        \; = \;
        \mathbf{1}[s_n \in \mathcal{S}'_n \text{ for at least } \nu \text{ different } n \in \mathcal{N}]
        \quad \Longleftrightarrow \quad
        \phi (s, a)
        \; = \;
        \mathbf{1} \left[ \sum_{n \in \mathcal{N}} \mathbf{1} [s_n \in \mathcal{S}'_n] \geq \nu \right],
    \end{equation*}
    and that the latter equation holds if and only if $\phi (s, a) = \eta$ for the auxiliary variable $\eta \in \mathbb{B}$ satisfying
    \begin{equation*}
        \frac{\displaystyle \left( \sum_{n \in \mathcal{N}} \mathbf{1} [s_n \in \mathcal{S}'_n] \right) - (\nu - 1)}{| \mathcal{N} | + 1}
        \;\; \leq \;\;
        \eta
        \;\; \leq \;\;
        1 + \frac{\displaystyle \left( \sum_{n \in \mathcal{N}} \mathbf{1} [s_n \in \mathcal{S}'_n] \right) - \nu}{| \mathcal{N} |}.
    \end{equation*}
    Since each indicator function $\mathbf{1} [s_n \in \mathcal{S}'_n]$ can be represented using $| \mathcal{S}'_n |$ auxiliary variables and $(1 + S_n) \cdot | \mathcal{S}'_n |$ constraints (\emph{cf.}~Assertion 1), the overall representation has $F_l = \sum_{n \in \mathcal{N}} | \mathcal{S}'_n |$ and $F_b = 1$ auxiliary variables as well as $F_c = \sum_{n \in \mathcal{N}} (1 + S_n) \cdot | \mathcal{S}'_n | + 2$ constraints.

    In view of the last assertion, finally, note that
    \begin{equation*}
        \phi (s, a)
        \; = \;
        \mathbf{1}[s_n \in \mathcal{S}'_n \text{ for at most } \nu \text{ different } n \in \mathcal{N}]
        \quad \Longleftrightarrow \quad
        \phi (s, a)
        \; = \;
        \mathbf{1} \left[ \sum_{n \in \mathcal{N}} \mathbf{1} [s_n \in \mathcal{S}'_n] \leq \nu \right],
    \end{equation*}
    and that the latter equation holds if and only if $\phi (s, a) = \eta$ for the auxiliary variable $\eta \in \mathbb{B}$ satisfying
    \begin{equation*}
        \frac{\displaystyle \nu - \left( \sum_{n \in \mathcal{N}} \mathbf{1} [s_n \in \mathcal{S}'_n] - 1 \right)}{| \mathcal{N} | + 1}
        \;\; \leq \;\;
        \eta
        \;\; \leq \;\;
        1 + \frac{\displaystyle \nu - \sum_{n \in \mathcal{N}} \mathbf{1} [s_n \in \mathcal{S}'_n]}{| \mathcal{N} |}.
    \end{equation*}
    Since each indicator function $\mathbf{1} [s_n \in \mathcal{S}'_n]$ can be represented using $| \mathcal{S}'_n |$ auxiliary variables and $(1 + S_n) \cdot | \mathcal{S}'_n |$ constraints (\emph{cf.}~Assertion 1), the overall representation has $F_l = \sum_{n \in \mathcal{N}} | \mathcal{S}'_n |$ and $F_b = 1$ auxiliary variables as well as $F_c = \sum_{n \in \mathcal{N}} (1 + S_n) \cdot | \mathcal{S}'_n | + 2$ constraints.
\qed

~\\[-8mm]

Similar to Theorem~1 of \cite{L97:prop_planning}, our proof of Theorem~\ref{thm:complexity_factored_mdps} relies on a reduction to a combinatorial game proposed by \cite{stockmeyer1979provably}. In contrast to our setting, however, \cite{L97:prop_planning} studies propositional planning problems that are structurally different from the FMDPs considered in this paper, and as a result our reduction, while relying on the same combinatorial game, requires a different proof strategy.

~\\[-12mm]

\noindent \textbf{Proof of Theorem~\ref{thm:complexity_factored_mdps}.}
    Since the feature map $\phi$ is assumed to be the identity, we will use Definition~\ref{def:fmdp} (relating to FMDPs) rather than Definition~\ref{def:fb-fmdp} (relating to feature-based FMDPs) throughout this proof. To see that computing the optimal expected total reward of an FMDP is EXPTIME-hard, we consider the following game proposed by \cite{stockmeyer1979provably}: \\
    
    \fbox{\parbox{15cm}{ {\centering \textsc{Two-Person Combinatorial Game $G_4$.}
      \\}
    \textbf{Instance.} Given a 13-DNF formula\footnotemark $F : \mathbb{B}^k \times \mathbb{B}^k \rightarrow \mathbb{B}$ and a starting position $(x^0, y^0) \in \mathbb{B}^k \times \mathbb{B}^k$. Player 1 takes the first turn. \\
    \textbf{Game.} The players take turns in switching at most one of their variables $\{ x_i \}_{i=1}^k$ (player 1) and $\{ y_i \}_{i=1}^k$ (player 2). The player whose switch causes $F (x, y)$ to evaluate to $1$ wins. \\
    \textbf{Question.} Does the starting position admit a winning strategy for player 1, that is, a strategy under which player 1 always wins, no matter what player 2 does?
  }} \footnotetext{13-DNF: $C_1 \vee \ldots C_m$, $m \in \mathbb{N}$, where each $C_i$ is a conjunction of at most 13 variables (or their negations).} \\
  
  \noindent Game $G_4$ is known to be EXPTIME-complete \citep[Theorem~3.1]{stockmeyer1979provably}.
  
  Without loss of generality, we can assume that each player switches exactly one variable in each turn. Indeed, the possibility to not switch any variable in a particular move can be catered for by adding an auxiliary variable that does not impact the value of $F$. Moreover, we will use the fact that if there is a winning strategy for player 1, then there is a winning strategy under which the player always wins in at most $2^{\mathcal{O} (k)}$ steps \citep[p.~160]{stockmeyer1979provably}.
  
  For a given instance of $G_4$, we first construct a simplified FMDP with scope $\mathcal{O} (k)$. We then outline how the FMDP can be modified so that is has scope $\mathcal{O} (\log k)$.
  
  \noindent \textbf{State Space.} The state space is $\mathcal{S} = \mathbb{B}^2 \times \mathbb{B}^k \times \mathbb{B}^k \times \mathbb{B}^{\lceil \log_2 T \rceil}$. A state $(w, x, y, z) \in \mathcal{S}$ records the status $w$ of the game ($(0, 0)$ = open, $(1, 0)$ = player~1 won, $(0, 1)$ = player~2 won), the current position $(x, y)$ as well as the number $z$ of completed turns (in binary representation).
  
  \noindent \textbf{Action Space.} We set $\mathcal{A} = \{ \mathrm{e}_i \, : \, i \in [k] \} \subseteq \mathbb{B}^k$, where $a \in \mathcal{A}$ records for each variable $i \in [k]$ whether $x_i$ is being flipped (if $a_i = 1$) or not (if $a_i = 0$) by player~1.
  
  \noindent \textbf{Initial Distribution.} The initial state is $(w^0, x^0, y^0, z^0)$ with $w^0 = (0, 0)$ and $z^0 = (0, \ldots, 0)$.
  
  \noindent \textbf{Transition Probabilities.} If the game status is $(1, 0)$ or $(0, 1)$, then the FMDP remains at the current state with probability $1$, independent of the selected action. If the game status is $(0, 0)$, on the other hand, then the position $x$ of player~1 transitions deterministically to $x'_i = x_i + (1 - 2 x_i) \cdot a_i$, $i \in [k]$, where $a \in \mathcal{A}$ is the current action. The position $y$ of player~2 transitions to each of the states $y' = y + (1 - 2 y_i) \cdot \mathrm{e}_i$, $i \in [k]$, with probability $1 / k$ each. The turn counter $z$ is increased by one (in binary representation). The game status $w$, finally, is set to $w' = (1, 0)$ if $F (x', y) = 1$, it is set to $w' = (0, 1)$ if $F (x, y) = 1$ or $z = (1, \ldots, 1)$, and it remains $w' = w$ otherwise. Note that the game status is backward looking with respect to the move of player~2, that is, a winning move of player~2 in the $z$-th turn is only recognized by the new game status $w'$ at the end of turn $z + 1$.
  
  \noindent \textbf{Reward Function.} We set $r ((w, x, y, z), a) = -1$ if $w = (0, 1)$ and $r ((w, x, y, z), a) = 0$ otherwise.
  
  \noindent \textbf{Discount Factor.} We choose any discount factor $\gamma \in (0, 1)$.

  We claim that the optimal expected total reward of the above FMDP vanishes if and only if player~1 has a winning strategy in game $G_4$. Indeed, one readily verifies that the FMDP tracks the game's evolution in a canonical way. Moreover, if player~1 has a winning strategy, then player~2 cannot win the game, which implies that no non-zero rewards are earned under the winning strategy, thus resulting in zero expected total rewards. If there is a policy earning zero expected total rewards, finally, then it corresponds to a strategy under which player~1 wins in no more than $T$ turns, irrespective of the moves of player~2 -- this is, by definition, a winning strategy.
  
  The FMDP described above can be implemented with a scope of $\mathcal{O} (k)$, that is, a scope that does not depend on the duration $T$ of the game. Indeed, through a judicious choice of features, the transition probabilities for each position $x'_i$ of player~1, $i \in [k]$, can be implemented with a scope of $4$ as they depend on the game status $w$, the previous position $x_i$ as well as the action $a_i$ selected for position $x_i$. The transition probabilities for each position $y'_i$ of player~2, $i \in [k]$, on the other hand, require a scope of $k + 2$ as they depend on the game status $w$ as well as all previous positions $\{ y_j \}_{j=1}^k$; this is needed to ensure that exactly one position is switched. The transition probabilities for each bit $z'_i$ of the turn counter, $i \in [\lceil \log_2 T \rceil]$, as well as the status $w'$ of the game, finally, can be implemented with a scope of $2$ since they evolve deterministically.

  The above FMDP can be modified so that it has a  scope of $\mathcal{O} (\log k)$. To this end, we augment the state space to $\mathcal{S} = \mathbb{B}^2 \times \mathbb{B}^k \times \mathbb{B}^k \times \mathbb{B}^{\lceil \log_2 k \rceil} \times \mathbb{B}^{\lceil \log_2 T \rceil}$. A state $(w, x, y, u, z)$ now additionally records in the game status $w$ whose turn it is ($(0, 0)$ = player 1's turn, $(0, 1)$ = player 2's turn, $(1, 0)$ = player 1 has won, $(1, 1)$ = player 2 has won) as well as, if player~2 moves next, which of her variables $y_j$, $j = \sum_{l = 1}^{\lceil \log_2 k \rceil} 2^{l - 1} \cdot u_l$, is switched. In that case, the scope of the transition probabilities for each position $y'_j$ of player~2, $j \in [k]$, can be reduced to $\log (k)$ since they become deterministic functions of the components of $u$, while the transition probabilities of $u$ have scope $0$ since they do not depend on any sub-state or action. For the sake of brevity, we omit the details of this tedious but otherwise straightforward extension.
\qed

~\\[-8mm]

\noindent \textbf{Proof of Proposition~\ref{prop:value_function_reformulation}.}
Replacing the value function $v$ in~\eqref{opt:fmdp:exact_lp} with our basis function approximation results in the following LP:
\begin{equation*}
    \begin{array}{l@{\quad}l}
        \displaystyle \mathop{\text{minimize}}_{w} & \displaystyle \sum_{s \in \mathcal{S}} \left[ \prod_{n \in [N]} q_n (s_n) \right] \left[ \sum_{k \in [K]} w_k \cdot \nu_k (s) \right]  \\
        \displaystyle \text{subject to} & \displaystyle \sum_{k \in [K]} w_k \cdot \nu_k (s) \geq \mspace{8mu} \sum_{j \in [J]} r_j (\phi (s, a)) + \gamma \sum_{s' \in \mathcal{S}} \left[ \prod_{n \in [N]} p_n (s'_n \, | \, \phi (s, a) )\right] \left[ \sum_{k \in [K]} w_k \cdot \nu_k ( s') \right] \\
        & \displaystyle \mspace{535mu} \forall (s, a) \in \mathcal{S} \times \mathcal{A} \\
        & \displaystyle w \in \mathbb{R}^K.
     \end{array}
\end{equation*}

In view of the objective function of this problem, we note that
\begin{align}
    & \sum_{s \in \mathcal{S}} \left[ \prod_{n \in [N]} q_n (s_n) \right] \left[ \sum_{k \in [K]} w_k \cdot \nu_k (s) \right]
    \;\; = \;\;
    \sum_{k \in [K]} w_k \sum_{s \in \mathcal{S}} \nu_k (s) \cdot \left[ \prod_{n \in [N]} q_n (s_n) \right] \nonumber \\
    = \;\; &
    \sum_{k \in [K]} w_k \sum_{\substack{s_n \in \mathcal{S}_n : \\ n \in \overline{\mathfrak{s}} [\nu_k]}} \; \sum_{\substack{s_{n'} \in \mathcal{S}_{n'} : \\ n' \in [N] \setminus \overline{\mathfrak{s}} [\nu_k]}} \nu_k (s) \cdot \left[ \prod_{n \in \overline{\mathfrak{s}} [\nu_k]} q_n (s_n) \right] \cdot \left[ \prod_{n' \in [N] \setminus \overline{\mathfrak{s}} [\nu_k]} q_{n'} (s_{n'}) \right] \nonumber \\
    = \;\; &
    \sum_{k \in [K]} w_k \sum_{s \in \overline{\mathfrak{S}} [\nu_k]} \sum_{\substack{s'_{n'} \in \mathcal{S}_{n'} : \\ n' \in [N] \setminus \overline{\mathfrak{s}} [\nu_k]}} \nu_k (s) \cdot \left[ \prod_{n \in \overline{\mathfrak{s}} [\nu_k]} q_n (s_n) \right] \cdot \left[ \prod_{n' \in [N] \setminus \overline{\mathfrak{s}} [\nu_k]} q_{n'} (s'_{n'}) \right], \label{eq:basis_fct_proof_eq1}
\end{align}
where the first identity applies standard algebraic manipulations, the second identity decomposes each state $s \in \mathcal{S}$ into its components $s_n$, $n \in \overline{\mathfrak{s}} [\nu_k]$, and $s_{n'}$, $n' \in [N] \setminus \overline{\mathfrak{s}} [\nu_k]$, and the third identity utilizes our definition of $\overline{\mathfrak{S}} [\nu_k]$ and the fact that $\nu_k$ does not depend on the components $s_{n'}$, $n' \in [N] \setminus \overline{\mathfrak{s}} [\nu_k]$. We next observe that~\eqref{eq:basis_fct_proof_eq1} can be re-expressed as
\begin{align*}
    & \sum_{k \in [K]} w_k \left( \sum_{s \in \overline{\mathfrak{S}} [\nu_k]} \nu_k (s) \cdot \left[ \prod_{n \in \overline{\mathfrak{s}} [\nu_k]} q_n (s_n) \right] \right) \cdot \left( \sum_{\substack{s'_{n'} \in \mathcal{S}_{n'} : \\ n' \in [N] \setminus \overline{\mathfrak{s}} [\nu_k]}} \prod_{n' \in [N] \setminus \overline{\mathfrak{s}} [\nu_k]} q_{n'} (s'_{n'}) \right) \\
    = \;\; &
    \sum_{k \in [K]} w_k \left( \sum_{s \in \overline{\mathfrak{S}} [\nu_k]} \nu_k (s) \cdot \left[ \prod_{n \in \overline{\mathfrak{s}} [\nu_k]} q_n (s_n) \right] \right),
\end{align*}
where the first row uses the distributive property of multiplication over addition, and the identity holds since the second multiplier in the first row evaluates to one. Similar arguments show that
\begin{equation*}
    \mspace{-12mu}
    \sum_{s' \in \mathcal{S}} \left[ \prod_{n \in [N]} p_n (s'_n \, | \, \phi (s, a) )\right] \left[ \sum_{k \in [K]} w_k \cdot \nu_k ( s') \right]
    \;\; = \;\;
    \sum_{k \in [K]} w_k \left( \sum_{s' \in \overline{\mathfrak{S}} [\nu_k]} \nu_k (s') \left[ \prod_{n \in \overline{\mathfrak{s}} [\nu_k]} p_n (s'_n \, | \, \phi (s, a)) \right] \right),
\end{equation*}
and we thus conclude that under our value function approximation,~\eqref{opt:fmdp:exact_lp} indeed simplifies to~\eqref{opt:fmdp:basis_fcts_lp_simplified}.

To prove the second part, we show that problem~\eqref{opt:fmdp:basis_fcts_lp_simplified} is feasible and not unbounded. The statement then follows from the fact that every feasible finite-dimensional LP attains its optimal value if the latter is finite. To see that~\eqref{opt:fmdp:basis_fcts_lp_simplified} is feasible,
fix $C = (1 - \gamma)^{-1} \cdot \max \{ r (\phi(s, a)) \, : \, (s, a) \in \mathcal{S} \times \mathcal{A} \}$ and let $w^\mathbf{1} \in \mathbb{R}^K$ be such that $\sum_{k \in [K]} w^\mathbf{1}_k \cdot \nu_k (s) = 1$ for all $s \in \mathcal{S}$. We then observe that
\begin{align*}
    \sum_{k \in [K]} (C \cdot w^\mathbf{1}_k) \cdot \nu_k (s)
    \;\; &= \;\;
    C
    \;\; = \;\;
    \frac{1}{1 - \gamma} \cdot \max \left\{ r (\phi(s', a)) \, : \, (s', a) \in \mathcal{S} \times \mathcal{A} \right\} \\
    &= \;\;
    \max \left\{ r (\phi(s', a)) \, : \, (s', a) \in \mathcal{S} \times \mathcal{A} \right\} + \gamma \cdot C \\
    &\geq \;\;
    \max_{a \in \mathcal{A}} \left[ r (\phi (s, a)) + \gamma \sum_{k \in [K]} (C \cdot w^\mathbf{1}_k) \sum_{s' \in \overline{\mathfrak{S}} [\nu_k]} \nu_k (s') \prod_{n \in \overline{\mathfrak{s}} [\nu_k]} p_n (s'_n \, | \, \phi (s, a)) \right]
\end{align*}
for all $s \in \mathcal{S}$, where the first two identities are due to the definitions of $w^\mathbf{1}$ and $C$, respectively, the third identity applies basic algebraic manipulations, and the inequality follows from the definition of $w^\mathbf{1}$, the fact that $\max \{ r (\phi (s, a)) \, : \, a \in \mathcal{A} \} \leq \max \{ r (\phi (s', a)) \, : \, (s', a) \in \mathcal{S} \times \mathcal{A} \}$, as well as the inequality $\sum_{s' \in \overline{\mathfrak{S}} [\nu_k]} \prod_{n \in \overline{\mathfrak{s}} [\nu_k]} p_n (s'_n \, | \, \phi (s, a)) \leq 1$, which holds for all $(s, a) \in \mathcal{S} \times \mathcal{A}$.

To see that problem~\eqref{opt:fmdp:basis_fcts_lp_simplified} is not unbounded, finally, note that problem~\eqref{opt:fmdp:exact_lp} bounds~\eqref{opt:fmdp:basis_fcts_lp_simplified} from below. 
As we have discussed in Section~\ref{sec:FMDP:cutting_plane}, however, the optimal value of~\eqref{opt:fmdp:exact_lp} coincides with the expected total reward of an optimal policy to the FMDP, which is bounded by construction.
\qed

~\\[-8mm]

The proof of Theorem~\ref{thm:Algo_1} relies on Hoffman's error bound, which we state next.

\begin{lem}[\citealt{hoffman2003approximate}]\label{lem:Hoffman}
    Let $\mathcal{K} = \{ x \in \mathbb{R}^n \, : \, A x \leq b \}$ be a non-empty polyhedron. Then there exists a constant $\kappa > 0$ such that for all $x \in \mathbb{R}^n$, the Euclidean distance between $x$ and $\mathcal{K}$ is bounded from above by $\kappa \cdot \left \lVert [A x - b]_+ \right \rVert_\infty$.
\end{lem}

\noindent \textbf{Proof of Theorem~\ref{thm:Algo_1}.}
To see that Algorithm~\ref{alg:CP} terminates in finite time, we note that by construction of the algorithm, every constraint of problem~\eqref{opt:fmdp:basis_fcts_lp_simplified} is added at most once. The claim now follows from the fact that problem~\eqref{opt:fmdp:basis_fcts_lp_simplified} comprises finitely many constraints.

To see that $F (w^\epsilon) \geq F^\star - \kappa_1 \cdot \epsilon$ for some $\kappa_1 > 0$ independent of $\epsilon$, denote by $\overline{w}^\epsilon \in \mathbb{R}^K$ the (unique) projection of $w^\epsilon$ onto the feasible region of problem~\eqref{opt:fmdp:basis_fcts_lp_simplified} and note that
\begin{align*}
    F^\star - F (w^\epsilon)
    \;\; &\leq \;\;
    F (\overline{w}^\epsilon) - F (w^\epsilon)
    \;\; = \;\;
    \sum_{s \in \mathcal{S}} q(s) \cdot \left[ \sum_{k \in [K]} (\overline{w}^\epsilon_k - w^\epsilon_k) \cdot \nu_k (s) \right] \\
    &\leq \;\;
    \max_{s \in \mathcal{S}} \left[ \sum_{k \in [K]} (\overline{w}^\epsilon_k - w^\epsilon_k) \cdot \nu_k (s) \right] \;\; \leq \;\;
    \max_{s \in \mathcal{S}} \Big[ \left \lVert \overline{w}^\epsilon - w^\epsilon \right \rVert_2 \cdot \left \lVert \nu (s) \right \rVert_2 \Big] \\
    &= \;\;
    \left \lVert \overline{w}^\epsilon - w^\epsilon \right \rVert_2 \cdot \max_{s \in \mathcal{S}} \Big[ \left \lVert \nu (s) \right \rVert_2 \Big],
\end{align*}
where the first inequality follows from the feasibility of $\overline{w}^\epsilon$ and the optimality of $F^\star$, the first identity employs formulation~\eqref{opt:fmdp:exact_lp} and Proposition~\ref{prop:value_function_reformulation} to evaluate $F (\overline{w}^\epsilon)$ and $F (w^\epsilon)$, the second inequality exploits that $q \in \Delta (\mathcal{S})$, the third inequality introduces the notation $\nu (s) = (\nu_1 (s), \dots, \nu_K (s))^\top$ and uses the Cauchy-Schwarz inequality, and the last identity reorders terms. We can now use Lemma~\ref{lem:Hoffman} to bound the term $\left \lVert \overline{w}^\epsilon - w^\epsilon \right \rVert_2$ from above by $\kappa_1 \cdot \epsilon$ if we identify $\mathcal{K}$ with the feasible region of problem~\eqref{opt:fmdp:basis_fcts_lp_simplified} and $x$ with $w^\epsilon$, respectively. The lemma then shows that $\left \lVert \overline{w}^\epsilon - w^\epsilon \right \rVert_2$, which amounts to the Euclidean distance between $x$ and $\mathcal{K}$, is bounded from above by a multiple of the maximum constraint violation of $w^\epsilon$ in problem~\eqref{opt:fmdp:basis_fcts_lp_simplified}, which itself is guaranteed to be bounded from above by $\epsilon$ by construction of Algorithm~\ref{alg:CP}. The statement then follows from the fact that $\max_{s \in \mathcal{S}} \big[ \left \lVert \nu (s) \right \rVert_2 \big]$ does not depend on $\epsilon$.

To see that $\hat{w}$ is feasible in problem~\eqref{opt:fmdp:basis_fcts_lp_simplified}, fix any $(s, a) \in \mathcal{S} \times \mathcal{A}$ and observe that
\begin{align*}
    & \sum_{k \in [K]} \hat{w}_k \cdot \nu_k (s)
    \;\; = \;\;
    \left[ \sum_{k \in [K]} w^\epsilon_k \cdot \nu_k (s) \right] + \frac{\epsilon}{1 - \gamma} \\
    \;\; \geq \;\;
    & r (\phi (s, a)) + \gamma \sum_{s' \in \mathcal{S}} p (s' \, | \, \phi (s, a)) \left[ \sum_{k \in [K]} w^\epsilon_k \cdot \nu_k (s') \right] - \epsilon + \frac{\epsilon}{1 - \gamma} \\
    \;\; = \;\;
    & r (\phi (s, a)) + \gamma \sum_{s' \in \mathcal{S}} p (s' \, | \, \phi (s, a)) \left[ \left( \sum_{k \in [K]} w^\epsilon_k \cdot \nu_k (s') \right) + \frac{\epsilon}{1 - \gamma} \right],
\end{align*}
where the first identity employs the definition of $\hat{w}$, the first inequality exploits the $\epsilon$-feasibility of $w^\epsilon$, and the second identity reorders terms and uses that $\sum_{s' \in \mathcal{S}} p (s' \, | \, \phi (s, a)) = 1$.
%Note next that
%\begin{equation*}
%    \left( \sum_{k \in [K]} w^\epsilon_k \cdot \nu_k (s) \right) + \frac{\epsilon}{1 - \gamma}
%    \;\; = \;\;
%    \sum_{k \in [K]} \left(w^\epsilon_k + \frac{\epsilon}{1 - \gamma} \cdot w^\mathbf{1}_k \right) \cdot \nu_k (s)
%    \;\; = \;\;
%    \sum_{k \in [K]} \hat{w}_k \cdot \nu_k (s),
%\end{equation*}
%where the first and second identity use the definition of $w^\mathbf{1}$ and $\hat{w}$, respectively. 
Applying the identity from the first row again to the term inside the square bracket in the last row, we obtain
\begin{equation*}
    \sum_{k \in [K]} \hat{w}_k \cdot \nu_k (s)
    \;\; \geq \;\;
    r (\phi (s, a)) + \gamma \sum_{s' \in \mathcal{S}} p (s' \, | \, \phi (s, a)) \left[ \sum_{k \in [K]} \hat{w}_k \cdot \nu_k (s) \right]
    \qquad \forall (s, a) \in \mathcal{S} \times \mathcal{A},
\end{equation*}
that is, $\hat{w}$ is indeed feasible in problem~\eqref{opt:fmdp:basis_fcts_lp_simplified}. The feasibility of $\hat{w}$ in~\eqref{opt:fmdp:basis_fcts_lp_simplified} and the fact that $w^\epsilon$ optimally solves a relaxation of~\eqref{opt:fmdp:basis_fcts_lp_simplified} then implies that $F (w^\epsilon) \leq F^\star \leq F (\hat{w})$.

To see that $F (\hat{w}) \leq F^\star + \kappa_2 \cdot \epsilon$ for some $\kappa_2 > 0$ independent of $\epsilon$, finally, we note that
\begin{equation*}
    F (\hat{w}) - F^\star
    \;\; \leq \;\;
    F (\hat{w}) - F (w^\epsilon)
    \;\; = \;\;
    \sum_{s \in \mathcal{S}} q(s) \cdot \left[ \sum_{k \in [K]} (\hat{w}_k - w^\epsilon_k) \cdot \nu_k (s) \right]
    \;\; = \;\;
    \frac{\epsilon}{1 - \gamma},
\end{equation*}
where the first inequality holds since $w^\epsilon$ optimally solves the master problem of the last iteration of Algorithm~\ref{alg:CP}, which itself is a relaxation of problem~\eqref{opt:fmdp:basis_fcts_lp_simplified}, the first identity employs formulation~\eqref{opt:fmdp:exact_lp} and Proposition~\ref{prop:value_function_reformulation} to evaluate $F (\hat{w})$ and $F (w^\epsilon)$, and the last identity follows from the definitions of $\hat{w}$ and $w^\mathbf{1}$ as well as the fact that $\sum_{s \in \mathcal{S}} q(s) = 1$.
\qed

~\\[-8mm]

\noindent \textbf{Proof of Theorem~\ref{thm:solution_subproblem}.}
    As discussed in Algorithm~\ref{alg:CP}, a maximally violated constraint can be extracted from an optimal solution to the optimization problem
    \begin{equation}\label{eq:proof_subproblem}
        \begin{array}{l@{\quad}l@{\qquad}l}
            \displaystyle \mathop{\text{maximize}}_{s, \, a, \, \varphi, \, \zeta} & \displaystyle \sum_{j \in [J]} r_j (\varphi) + \gamma \sum_{k \in [K]} w^\star_k \cdot \overline{\nu}_k (\varphi) - \sum_{k \in [K]} w^\star_k \cdot \nu_k (s) \\
            \text{subject to} & \displaystyle (s, a) \in \mathcal{S} \times \mathcal{A}, \;\; (\varphi, \zeta; s, a) \in \mathcal{F}, \;\; \zeta \in \mathbb{R}^{F_l} \times \mathbb{B}^{F_b}.
        \end{array}
    \end{equation}
    The reward components $r_j (\varphi)$, $j \in [J]$, in this problem can be reformulated as 
    \begin{equation*}
        r_j (\varphi)
        \;\; = \;\;
        \sum_{f \in \mathfrak{S} [r_j]} r_j(f) \cdot \mathbf{1} \left[ \varphi_i = f_i \;\; \forall i \in \mathfrak{s} [r_j]  \right]
        \;\; = \;\;
        \sum_{f \in \mathfrak{S} [r_j] } r_j (f) \cdot \eta_{jf}
    \end{equation*}
    for $\eta_{jf} = \mathbf{1} [\varphi_i = f_i \;\; \forall i \in \mathfrak{s} [r_j] ]$, $j \in [J]$ and $f \in \mathfrak{S} [r_j]$. One readily verifies that $\eta_{jf} = \mathbf{1} [\varphi_i = f_i \;\; \forall i \in \mathfrak{s} [r_j] ]$ if and only if $\eta_{jf} \in [0,1]$ and
    \begin{equation*}
        \eta_{jf} \leq (2 f_i-1) \varphi_i + 1 - f_i \;\; \forall i \in \mathfrak{s} [r_j] \quad \text{as well as} \quad \eta_{jf} \geq 1 + \sum_{i \in \mathfrak{s} [r_j]}\frac{ \varphi_i - f_i}{2 f_i-1}.
    \end{equation*}
    Similar arguments show that
    \begin{equation*}
        \overline{\nu}_k(\varphi)
        \;\; = \;\;
        \sum_{f \in \mathfrak{S} [\overline{\nu}_k]} \overline{\nu}_k (f) \cdot \mathbf{1} \left[ \varphi_i = f_i \;\; \forall i \in \mathfrak{s} [\overline{\nu}_k] \right]
        \;\; = \;\;
        \sum_{f \in \mathfrak{S} [\overline{\nu}_k]} \overline{\nu}_k (f) \cdot \xi_{kf}
    \end{equation*}
    for $\xi_{kf} = \mathbf{1} [\varphi_i = f_i \;\; \forall i \in \mathfrak{s} [\overline{\nu}_k] ]$, $k \in [K]$ and $f \in \mathfrak{S} [\overline{\nu}_k]$, which holds if and only if $\xi_{kf} \in [0,1]$ and
    \begin{equation*}
        \xi_{kf} \leq (2 f_i-1) \varphi_i + 1 - f_i \;\; \forall i \in \mathfrak{s} [\overline{\nu}_k] \quad \text{as well as} \quad \xi_{kf} \geq 1 + \sum_{i \in \mathfrak{s} [\overline{\nu}_k]} \frac{\varphi_i - f_i }{2 f_i-1},
    \end{equation*}
    as well as
    \begin{equation*}
        \nu_k(s)
        \;\; = \;\;
        \sum_{s' \in \mathfrak{S} [\nu_k]} \nu_k(s') \cdot \mathbf{1} \left[ [s]_i = [s']_i \;\; \forall i \in \mathfrak{s} [\nu_k] \right]
        \;\; = \;\;
        \sum_{s' \in \mathfrak{S} [\nu_k]} \nu_k(s') \cdot \beta_{ks'}
    \end{equation*}
    for $\beta_{ks'} = \mathbf{1} [[s]_i = [s']_i \;\; \forall i \in \mathfrak{s} [\nu_k]]$, $k \in [K]$ and $s' \in \mathfrak{S} [\nu_k]$, which holds if and only if $\beta_{ks'} \in [0, 1]$ and
    \begin{equation*}
        \beta_{ks'} \leq (2 [s']_i-1) [s]_i +1-[s']_i \;\; \forall i \in \mathfrak{s} [\nu_k] \quad \text{as well as} \quad \beta_{ks'} \geq 1+\sum_{i \in \mathfrak{s} [\nu_k]}\frac{ [s]_i - [s']_i}{2 [s']_i-1}.
    \end{equation*}
    The statement now follows if we replace $r_j (\varphi)$, $\overline{\nu}_k (\varphi)$ and $\nu_k (s)$ in~\eqref{eq:proof_subproblem} with their reformulations.
\qed

~\\[-8mm]

\noindent \textbf{Proof of Proposition~\ref{prop:complexity_fmdp_subproblem}.}
    Theorem~\ref{thm:solution_subproblem} determines a maximally violated constraint through the solution of an MILP whose size is polynomial in the description of the feature-based FMDP. We thus conclude that determining a maximally violated constraint is in NP. To see that the problem is also strongly NP-hard, consider the following integer feasibility problem:
    \begin{center}
        \fbox{\parbox{10cm}{ {\centering \textsc{0/1 Integer Programming.}
      \\}
        \textbf{Instance.} Given are $G \in \mathbb{Z}^{L \times N}$ and $h \in \mathbb{Z}^L$. \\
        \textbf{Question.} Is there a vector $x \in \left\{ 0, 1 \right\}^N$ such that $G x \leq h$?
        }} \\
    \end{center}
    For a given instance of the integer feasibility problem, we construct a feature-based FMDP with the state space $\mathcal{S} = \bigtimes_{n \in [N]} \mathcal{S}_n$ with $\mathcal{S}_n = \mathbb{B}$, $n \in [N]$, the action space $\mathcal{A} = \mathbb{B}$, the feature map $\phi : \mathcal{S} \times \mathcal{A} \rightarrow \mathbb{B}$ satisfying $\phi (s, a) = \mathbf{1} [ G s \leq h ]$, the reward function satisfying $r (s, a) = \phi (s, a)$, any scope-1 initial distribution $q$, any scope-1 transition kernel $p$ as well as any discount factor; we not specify those components further since they will not impact our argument. We set the basis function weights $w$ to zero so that the choice of basis functions for the value function is irrelevant as well. In this case, any maximally violated constraint $(s^\star, a^\star)$ simply maximizes $\phi (s, a)$ over $\mathcal{S} \times \mathcal{A}$, and we conclude that the maximum constraint violation is $1$ precisely when the answer to the integer feasibility problem is affirmative. Since the integer feasibility problem is strongly NP-complete \citep{GJ79:ComputersIntractability}, the statement of the proposition then follows.

    It remains to be shown that the feature map $\phi : \mathcal{S} \times \mathcal{A} \rightarrow \mathbb{B}$ admits a representation that is polynomial in the length of the integer feasibility instance. To this end, we claim that
    \begin{equation}\label{eq:np_hard_phi}
        \phi (s, a) = \varphi
        \quad \Longleftrightarrow \quad
        \left[
        \begin{array}{l@{\quad}l}
        \displaystyle g_l{}^\top s \geq h_l + \frac{1}{2} - \mathrm{M} \cdot \zeta_l & \displaystyle \forall l \in [L] \\
        \displaystyle g_l{}^\top s \leq h_l + \mathrm{M} \cdot (1 - \zeta_l) & \displaystyle \forall l \in [L] \\
        \displaystyle \varphi \leq \zeta_l & \displaystyle \forall l \in [L] \\
        \displaystyle \varphi \geq \left( \sum_{l \in [L]} \zeta_l \right) - (L - 1)
        \end{array}
        \right]
        \text{ for some } \zeta \in \mathbb{B}^L,
    \end{equation}
    where $g_l{}^\top \in \mathbb{Z}^N$ is the $l$-th row of matrix $G$ and $\mathrm{M} = \max_{l \in [L]} \left[ 1 + | h_l | + \sum_{n \in [N]} | g_{ln} | \right]$. To see this, note that thanks to the integrality of $G$ and $h$, the equation system~\eqref{eq:np_hard_phi} ensures that $\zeta_l = \mathbf{1} [ g_l{}^\top s \leq h_l ]$ for all $l \in [L]$ as well as $\varphi = \min \{ \zeta_l \, : \, l \in [L] \}$. Thus, whenever $\phi (s, a) = 1$, that is, whenever $G s \leq h$, we have $\zeta = \mathbf{1}$ and thus $\varphi = 1$. In contrast, whenever $\phi (s, a) = 0$, that is, whenever $G s \not\leq h$, we have $\zeta_l = 0$ for at least one $l \in [L]$, and thus $\varphi = 0$ as desired.
\qed

~\\[-8mm]

\noindent \textbf{Proof of Corollary~\ref{cor:determine_best_action}.}
    The statement follows immediately from the optimization problem in Theorem~\ref{thm:solution_subproblem} if we fix the state $s$, disregard the last expression from the objective function and remove the auxiliary variables $\{ \beta_{ks'} \}$.
\qed

~\\[-8mm]

\noindent \textbf{Proof of Observation~\ref{obs:rip}.}
    To see that the first set of expressions in the proposition is a subset of the second set of expressions, fix any joint distribution $p_{sa} \in \mathcal{P} (\phi (s, a))$ and construct the marginal distributions $\{ p_{sak} \}_k$ via $p_{sak} (s') = p_{sa} (s')$ for all $k \in [K]$ and $s' \in \overline{\mathfrak{S}} [\nu_k]$. One readily verifies that
    \begin{equation}\label{eq:rip:equivalence}
        \sum_{k \in [K]} \sum_{s' \in \overline{\mathfrak{S}} [\nu_k]} w_k \cdot \nu_k (s') \cdot p_{sa} (s')
        \;\; = \;\;
        \sum_{k \in [K]} \sum_{s' \in \overline{\mathfrak{S}} [\nu_k]} w_k \cdot \nu_k (s') \cdot p_{sak} (s')
    \end{equation}
    and that the marginal distributions $\{ p_{sak} \}_k$ satisfy the conditions of the second set of expressions.

    To see that both sets of expressions coincide whenever the basis functions $\nu_k$, $k \in [K]$, satisfy the RIP, we refer to \citet[\S~2.1]{doan2015robustness}, whose results imply that under the RIP, any set of marginal distributions $\{ p_{sak} \}_k$ satisfying the conditions of the second set of expressions gives rise to a joint distribution $p_{sa} \in \mathcal{P} (\phi (s, a))$ satisfying~\eqref{eq:rip:equivalence}.
\qed

~\\[-8mm]

The proof of Theorem~\ref{thm:lambda-epsilon-kkt-point} relies on the following generalization of Hoffman's error bound to polynomial systems, which we state first.

\begin{lem}[\cite{li2015new}]\label{lem:poly_err_bound}
    Consider the set $\mathcal{K} = \{ x \in \mathbb{R}^n : g_i (x) \leq 0 \; \forall i \in [m_1], \; h_j(x) = 0 \; \forall j \in [m_2] \}$ defined by the polynomials $g_i : \mathbb{R}^n \rightarrow \mathbb{R}$, $i\in [m_1]$, and $h_j : \mathbb{R}^n \rightarrow \mathbb{R}$, $j\in[m_2]$, of degree at most $d$. For any compact set $\mathcal{K}' \subseteq \mathbb{R}^n$ there exists a constant $\kappa >0$ such that
    \begin{align*}
        \mathrm{dist}(x, \mathcal{K})
        \;\; \leq \;\;
        \kappa \left( \sum_{i = 1}^{m_1} [ g_i(x) ]_+ +\sum_{j=1}^{m_2} |h_j(x)| \right)^\tau \qquad \forall x \in \mathcal{K}',
    \end{align*} 
    where $\mathrm{dist}$ is the standard point-to-set distance in $\mathbb{R}^n$ and
    \begin{align*}
        \tau = \max\left\{  \frac{1}{(d+1)(3d)^{ n+m_1+m_2 - 1 }}, \; \frac{1}{d(6d-3)^{n+m_1 - 1}} \right\}.
    \end{align*}
\end{lem}

\noindent \textbf{Proof of Theorem~\ref{thm:lambda-epsilon-kkt-point}.} $\;$
Since problem~\eqref{opt:rfmdp:basis_fcts_lp_simplified} comprises finitely many bilinear constraints parametrized by $(s,a) \in \mathcal{S} \times \mathcal{A}$, Algorithm~\ref{alg:CP-robust} terminates after a finite number of iterations.

To see that $(\hat{w}, \hat{p})$ is feasible in~\eqref{opt:rfmdp:basis_fcts_lp_simplified}, note that $\hat{p}_{sa} \in \mathcal{P} (\phi (s, a))$ for all $(s, a) \in \mathcal{S} \times \mathcal{A}$ by construction. \mbox{Moreover, $(\hat{w}, \hat{p})$ satisfies the first constraint set in~\eqref{opt:rfmdp:basis_fcts_lp_simplified} since for any $(s, a) \in \mathcal{S} \times \mathcal{A}$, we have}
\begin{align*}
    \sum_{k \in [K]} \hat{w}_k \cdot \nu_k (s)
    \;\; &= \;\;
    \left[ \sum_{k \in [K]} w^\epsilon_k \cdot \nu_k (s) \right] + \frac{\epsilon}{1 - \gamma} \\
    &\geq \;\;
    \sum_{j \in [J]} r_j (\phi (s, a)) + \gamma \sum_{k \in [K]} w^\epsilon_k \sum_{s' \in \overline{\mathfrak{S}} [\nu_k]} \nu_k (s') \cdot \hat{p}_{sa} (s') + \frac{\epsilon \gamma}{1 - \gamma} \\
    &= \;\;
    \sum_{j \in [J]} r_j (\phi (s, a)) + \gamma \sum_{k \in [K]} \hat{w}_k \sum_{s' \in \overline{\mathfrak{S}} [\nu_k]} \nu_k (s') \cdot \hat{p}_{sa} (s'),
\end{align*}
where the first identity follows from the definition of $\hat{w}$. The first inequality holds for all $(s, a) \in \mathcal{C}$ even if $\epsilon \gamma / (1 - \gamma)$ on the right-hand side is replaced with the larger constant $\epsilon / (1 - \gamma)$ since $\hat{p}_{sa} = p^\epsilon_{sa}$ in that case and the solution $(w^\epsilon, p^\epsilon)$ satisfies all constraints $(s, a) \in \mathcal{C}$ by construction. The first inequality also holds for all $(s, a) \in [\mathcal{S} \times \mathcal{A}] \setminus \mathcal{C}$ since $\hat{p}_{sa}$ minimizes the right-hand side across all $p_{sa} \in \mathcal{P} (\phi (s, a))$ in that case and the solution $(w^\epsilon, p^\epsilon)$ is $\epsilon$-feasible when Algorithm~\ref{alg:CP-robust}
terminates. The second identity, finally, follows again from the definition of $\hat{w}$.

We complete the proof by showing that the Euclidean distance of $(\hat{w}, \hat{p})$ to the stationary points of problem~\eqref{opt:rfmdp:basis_fcts_lp_simplified} is bounded from above as claimed. To this end, we proceed in three steps: \emph{(i)} we derive the Karush-Kuhn-Tucker (KKT) system for the terminal master problem of Algorithm~\ref{alg:CP-robust}; \emph{(ii)} we show that the solutions to the KKT system lie in a bounded set; \emph{(iii)} we employ Lemma~\ref{lem:poly_err_bound} to bound the  distance of $(\hat{w}, \hat{p})$ to the stationary points of~\eqref{opt:rfmdp:basis_fcts_lp_simplified}.

In view of the first step, assume that $\mathcal{P} (\phi (s, a)) = \{ p_{sa} \in \mathbb{R}^{| \mathcal{S} |} : \Psi_{sa} p_{sa} \leq \psi_{sa} \}$ for some $\Psi_{sa} \in \mathbb{R}^{ M_{sa}\times |\mathcal{S}| }$ with full row rank and $\psi_{sa} \in \mathbb{R}^{ M_{sa} }$. In that case, the KKT system of the terminal master problem of Algorithm~\ref{alg:CP-robust} presents itself as follows.
\begin{equation*}
    \boxed{
    \begin{array}{c}
        \text{\textbf{Stationarity:}} \\
        \displaystyle \sum_{s \in \overline{\mathfrak{S}} [\nu_k]} \nu_k (s) \prod_{n \in \overline{\mathfrak{s}} [\nu_k]} q_n (s_n) + \sum_{(s, a) \in \mathcal{C}} \Lambda_{sa} \left( \gamma \sum_{s' \in \overline{\mathfrak{S}} [\nu_k]} \nu_k (s') \cdot p_{sa} (s') - \nu_k (s) \right) = 0 \qquad \forall k \in [K] \\
        \displaystyle \gamma \cdot \Lambda_{sa} \sum_{k \in [K]} w_k \cdot \nu_k (s') + [\Psi_{sa}]_{s'}{}^\top \Gamma_{sa} = 0 \qquad \forall (s, a) \in \mathcal{C}, \; \forall s' \in \mathcal{S} \\
        \text{\textbf{Primal feasibility:}} \\
        \left.
        \begin{array}{c}
            \displaystyle \sum_{k \in [K]} w_k \cdot \nu_k (s) \geq r (\phi (s, a)) + \gamma \sum_{k \in [K]} w_k \sum_{s' \in \overline{\mathfrak{S}} [\nu_k]} \nu_k (s') \cdot p_{sa} (s') \\
            \displaystyle \Psi_{sa} p_{sa} \leq \psi_{sa}
        \end{array}
        \quad \right\} \quad \forall (s, a) \in \mathcal{C} \\
        \text{\textbf{Dual feasibility:}} \\
        \Lambda_{sa} \in \mathbb{R}_+ \;\; \text{and} \;\; \Gamma_{sa} \in \mathbb{R}^{ M_{sa} }_+ \qquad \forall (s, a) \in \mathcal{C} \\
        \text{\textbf{Complementary slackness:}} \\
        \left.
        \begin{array}{c}
            \displaystyle \Lambda_{sa} \left( r (\phi (s, a)) + \gamma \sum_{k \in [K]} w_k \sum_{s' \in \overline{\mathfrak{S}} [\nu_k]} \nu_k (s') \cdot p_{sa} (s') - \sum_{k \in [K]} w_k \cdot \nu_k (s) \right) = 0 \\
            \displaystyle \Gamma_{sa}{}^\top \left( \Psi_{sa} p_{sa} - \psi_{sa} \right) = 0
        \end{array}
        \quad \right\} \quad \forall (s, a) \in \mathcal{C}
    \end{array}}
\end{equation*}

\noindent Here, $\Lambda_{sa}$ and $\Gamma_{sa}$ are the dual variables corresponding to the primal Bellman constraints and the constraints defining the ambiguity sets, respectively.

For the second step, we bound each of the components $(w, p, \Lambda, \Gamma)$ of a KKT point separately. To bound $w$, observe that the first stationarity, primal feasibility, dual feasibility and complementary slackness condition imply that $w$ optimally solves the linear program
\begin{equation}\label{eq:lp_subproblem_from_kkt}
    \begin{array}{l@{\quad}l}
        \displaystyle \mathop{\text{minimize}}_{w} & \displaystyle \sum_{k \in [K]} w_k \sum_{s \in \overline{\mathfrak{S}} [\nu_k]} \nu_k (s) \prod_{n \in \overline{\mathfrak{s}} [\nu_k]} q_n (s_n) \\
        \displaystyle \text{subject to} & \displaystyle \sum_{k \in [K]} w_k \cdot \nu_k (s) \geq \sum_{j \in [J]} r_j (\phi (s, a)) + \gamma \sum_{k \in [K]} w_k \sum_{s' \in \overline{\mathfrak{S}} [\nu_k]} \nu_k (s') \cdot p_{sa} (s') \qquad \forall (s, a) \in \mathcal{C}.
     \end{array}
\end{equation}
Our earlier assumption that the initial master problem is not unbounded immediately implies a lower bound on the optimal value of~\eqref{eq:lp_subproblem_from_kkt}. Likewise, one readily verifies that $\overline{w}$ with $\overline{w} = (1 - \gamma)^{-1} \cdot \max \{ r(s, a) : (s, a) \in \mathcal{S} \times \mathcal{A} \} \cdot w^\mathbf{1}$ is feasible in~\eqref{eq:lp_subproblem_from_kkt} and thus provides the upper bound $(1 - \gamma)^{-1} \cdot \max \{ r(s, a) : (s, a) \in \mathcal{S} \times \mathcal{A} \}$ on the optimal value of~\eqref{eq:lp_subproblem_from_kkt}. The boundedness of this optimal value does not \emph{per se} imply boundedness of the associated optimal cost to-go functions $\sum_{k \in [K]} w_k \cdot \nu_k (\cdot)$, however, since some states $s \in \mathcal{S}$ may never be reached in the FMDP represented by problem~\eqref{eq:lp_subproblem_from_kkt}. Assume to the contrary that there is a sequence of optimal solutions to problem~\eqref{eq:lp_subproblem_from_kkt} whose cost to-go would diverge for some states. Then the same cost to-go would diverge in the associated sequence of KKT points to the terminal master problem of Algorithm~\ref{alg:CP-robust}. The assumed irreducibility of the RFMDP under some $p \in \mathcal{P}$ would then imply, however, that the associated sequence of objective values in the master problem would diverge, which contradicts our earlier assumption that the initial master problem is not unbounded.
%these cost to-go functions would diverge for \emph{all} states due to the optimal choice of $p_{sa}$ in the subproblems, which contradicts our earlier assumption that the initial master problem is not unbounded.
The boundedness of $w$ in any KKT point now follows from the fact that the mapping $w \mapsto \sum_{k \in [K]} w_k \cdot \nu_k (\cdot)$ is injective due to the assumed linear independence of the basis functions. The second primal feasibility condition immediately implies that the probabilities $p_{sa}$, $(s, a) \in \mathcal{C}$, are also contained in a bounded set. To see that $\Lambda$ is bounded as well, we multiply each of the $K$ equations in the stationarity condition with $w^\mathbf{1}_k$, $k \in [K]$, and sum them up to obtain
\begin{equation*}
    1 + \sum_{(s, a) \in \mathcal{C}} \Lambda_{sa} (\gamma - 1) = 0
    \quad \Longleftrightarrow \quad
    \sum_{(s, a) \in \mathcal{C}} \Lambda_{sa} = \frac{1}{1 - \gamma}.
\end{equation*}
The non-negativity of $\Lambda$ (\emph{cf.}~the dual feasibility conditions) then establishes the boundedness of $\Lambda$. The boundedness of $\Gamma$, finally, follows from the boundedness of $\Lambda$ and $w$, the full row rank of each $\Psi_{sa}$, the non-negativity of $\Gamma$ and the second stationarity condition.

In view of the third step, finally, observe that the KKT system for problem~\eqref{opt:rfmdp:basis_fcts_lp_simplified} coincides with our above KKT system for the terminal master problem of Algorithm~\ref{alg:CP-robust}, except that the index set $\mathcal{C}$ is replaced with $\mathcal{S} \times \mathcal{A}$. Recall that by definition, $(w^\epsilon, p^\epsilon)$ is a stationary point to the terminal master problem of Algorithm~\ref{alg:CP-robust}. In other words, together with some $\Lambda^\epsilon = \{\Lambda^\epsilon_{sa}\}_{(s,a)\in \mathcal{C}}$ and $\Gamma^\epsilon = \{\Gamma^\epsilon_{sa}\}_{(s,a)\in \mathcal{C}}$, $(w^\epsilon, p^\epsilon, \Lambda^\epsilon, \Gamma^\epsilon)$ is a KKT point of the terminal master problem. We can extend this KKT point from $(s, a) \in \mathcal{C}$ to $(s, a) \in \mathcal{S} \times \mathcal{A}$ by selecting $p^\epsilon_{sa} \in \arg\min \big\{ \sum_{k \in [K]} w^\epsilon_k \sum_{s' \in \overline{\mathfrak{S}} [\nu_k]} \nu_k (s') \cdot p_{sa} (s') \, : \, p_{sa}\in \mathcal{P} (\phi (s, a)) \big\}$ and  $(\Lambda^\epsilon_{sa}, \Gamma^\epsilon_{sa}) = (0, 0)$ for $(s, a) \in [\mathcal{S} \times \mathcal{A}] \setminus \mathcal{C}$. The resulting solution $(w^\epsilon, p^\epsilon, \Lambda^\epsilon, \Gamma^\epsilon)$ is an 
\emph{approximate} KKT point for~\eqref{opt:rfmdp:basis_fcts_lp_simplified} in that it violates the first primal feasibility condition by up to $\epsilon$ but satisfies all other KKT conditions. By Lemma~\ref{lem:poly_err_bound}, there exists a constant $\kappa (\mathcal{C}) > 0$ such that
\begin{equation}
    \label{ineq:poly_err_bound}
    \mathrm{dist} ((w^\epsilon, p^\epsilon, \Lambda^\epsilon, \Gamma^\epsilon), \mathcal{K}) \le \kappa(\mathcal{C}) \cdot (| \mathcal{S} \times \mathcal{A} | \cdot \epsilon)^\tau,
\end{equation}
where $\mathcal{K}$ and $\mathcal{K}'$ represent the compact sets that bound the set of KKT points for problem~\eqref{opt:rfmdp:basis_fcts_lp_simplified} and the terminal master problem of Algorithm~\ref{alg:CP-robust}, respectively.
Let $\kappa = \max \{ \kappa (\mathcal{C}') \, : \, \mathcal{C}' \subseteq \mathcal{S} \times \mathcal{A} \}$. Note that although the terminal index set $\mathcal{C}$ depends on $\epsilon$, the constant $\kappa$ does not. We then have
\begin{align*}
    \mathrm{dist} ((\hat{w}, \hat{p}), \mathcal{X})
    \;\; &\leq \;\;
    \mathrm{dist} ((w^\epsilon, p^\epsilon), \mathcal{X}) + \frac{\epsilon}{1 - \gamma} \left \lVert w^\mathbf{1} \right \rVert_2
    \;\; \leq \;\;
    \mathrm{dist} ((w^\epsilon, p^\epsilon, \Lambda^\epsilon, \Gamma^\epsilon), \mathcal{K}) + \frac{\epsilon}{1 - \gamma} \left \lVert w^\mathbf{1} \right \rVert_2 \\
    &\leq \;\;
    \kappa \cdot (| \mathcal{S} \times \mathcal{A} | \cdot \epsilon)^\tau + \frac{\epsilon}{1 - \gamma} \left \lVert w^\mathbf{1} \right \rVert_2,
\end{align*}
where $\mathcal{X}$ denotes the set of stationary points for problem~\eqref{opt:rfmdp:basis_fcts_lp_simplified}.
Here, the first inequality follows from the definition of $(\hat{w}, \hat{p})$ and the triangle inequality. The second inequality follows from the fact the the distance of a point to a set is at least its the distance of its projection to the projection of the set. The final inequality uses inequality~\eqref{ineq:poly_err_bound}. The claim of the theorem now follows if we redefine $\kappa$ as $\kappa \cdot (| \mathcal{S} \times \mathcal{A} |)^\tau$.
\qed

~\\[-8mm]

\noindent \textbf{Proof of Theorem~\ref{thm:kkt-convergence}}
We prove the statement in three steps: \emph{(i)} we argue that the sequence $\{ (w^t, p^t) \}_{t \in \mathbb{N}}$ produced by Algorithm~\ref{eq:update} is bounded and thus has at least one limit point $(w^\infty, p^\infty)$; \emph{(ii)} we show that for any such limit point $(w^\infty, p^\infty)$, $w^\infty$ and $p^\infty$ are minimizers in the first and second linear program in Algorithm~\ref{eq:update} if we replace $p^{t-1}$ with $p^\infty$ and $w^t$ with $w^\infty$, respectively; and \emph{(iii)} we show that the KKT conditions of the two linear programs from step \emph{(ii)} imply that $(w^\infty, p^\infty)$ can be augmented to a KKT point $(w^\infty, p^\infty, \Lambda, \Gamma)$ to the restricted master problem~\eqref{opt:rfmdp:basis_fcts_lp_simplified}, that is,  $(w^\infty, p^\infty)$ is a stationary point to the restricted master problem~\eqref{opt:rfmdp:basis_fcts_lp_simplified} as claimed.

In view of the first step, we note that similar arguments as in step~\emph{(ii)} of the proof of Theorem~\ref{thm:lambda-epsilon-kkt-point} show that the sequence $\{ w^t \}_{t \in \mathbb{N}}$ is bounded. Since the sequence $\{ p^t \}_{t \in \mathbb{N}}$ is bounded by construction, we conclude that the sequence $\{ (w^t, p^t) \}_{t \in \mathbb{N}}$ produced by Algorithm~\ref{eq:update} is bounded as well.

%By going over to subsequences if necessary, we assume without loss of generality that $\{ (w^t, p^t) \}_{t \in \mathbb{N}}$ itself converges to $(w^\infty, p^\infty)$.

As for the second step, fix any limit point $(w^\infty, p^\infty)$ of the  sequence $\{ (w^t, p^t) \}_{t \in \mathbb{N}}$ produced by Algorithm~\ref{eq:update}, and let $\{ (w^{t_j}, p^{t_j}) \}_{j \in \mathbb{N}}$ be a subsequence that converges to $(w^\infty, p^\infty)$. To see that $w^\infty$ minimizes the first LP in Algorithm~\ref{eq:update} if we replace $p^{t-1}$ with $p^\infty$, note first that for all $j \in \mathbb{N}$,
\begin{align*}
    \sum_{k \in [K]} w^{t_j}_k \cdot \nu_k (s)
    &\geq
    \sum_{j \in [J]} r_j (\phi (s, a)) + \gamma \sum_{k \in [K]} w^{t_j}_k \sum_{s' \in \overline{\mathfrak{S}} [\nu_k]} \nu_k (s') \cdot p^{t_j-1}_{sak} (s') \\
    &\geq
    \sum_{j \in [J]} r_j (\phi (s, a)) + \gamma \sum_{k \in [K]} w^{t_j}_k \sum_{s' \in \overline{\mathfrak{S}} [\nu_k]} \nu_k (s') \cdot p^{t_j}_{sak} (s')    
    \qquad \forall (s, a) \in \mathcal{C},
\end{align*}
where the first and second inequality follow directly from the construction of the first and second optimization problem in Algorithm~\ref{eq:update}, respectively. Taking the limit $j \longrightarrow \infty$ on both sides of this equation shows that
\begin{equation*}
    \sum_{k \in [K]} w^\infty_k \cdot \nu_k (s) \geq \sum_{j \in [J]} r_j (\phi (s, a)) + \gamma \sum_{k \in [K]} w^\infty_k \sum_{s' \in \overline{\mathfrak{S}} [\nu_k]} \nu_k (s') \cdot p^\infty_{sak} (s')
    \qquad \forall (s, a) \in \mathcal{C},
\end{equation*}
that is, $w^\infty$ is \emph{feasible} in the first LP in Algorithm~\ref{eq:update} if we replace $p^{t-1}$ with $p^\infty$. To see that $w^\infty$ is also \emph{optimal} in the first LP in Algorithm~\ref{eq:update} if we replace $p^{t-1}$ with $p^\infty$, we note that
\begin{align*}
    f (w^{t+1})
    \;\; \leq \;\;
    f (w^{t+1}) + \lambda_t \left \lVert w^{t+1} - w^t \right \rVert_1
    \;\; \leq \;\;
    f (w^t)
    \qquad \forall t \in \mathbb{N},
\end{align*}
where
\begin{equation*}
    f (w)
    \; = \;
    \sum_{k\in [K]} w_k \sum_{s \in \overline{\mathfrak{S}} [\nu_k]} \nu_k (s) \prod_{n \in \overline{\mathfrak{s}} [\nu_k]} q_n (s_n)
\end{equation*}
denotes the unregularized objective function of the first LP in Algorithm~\ref{eq:update}. Here, the first inequality holds by construction, and the second inequality is due to the feasibility (but typically suboptimality) of $w^t$ in the first LP in Algorithm~\ref{eq:update} at iteration $t + 1$. We thus conclude that
\begin{align}
    f (w^{t_{j+1}}) + \lambda_{t_j + 1} \left \lVert w^{t_j+1} - w^{t_j} \right \rVert_1
    \;\; &\leq \;\;
    f (w^{t_j + 1}) + \lambda_{t_j + 1} \left \lVert w^{t_j+1} - w^{t_j} \right \rVert_1 \nonumber \\
    &= \;\;
    \min_{w \in \Omega (p^{t_j})} \left\{ f (w) + \lambda_{t_j + 1} \left \lVert w - w^{t_j} \right \rVert_1 \right\}, \label{eq:limit_condition}
\end{align}
where
\begin{equation*}
    \mspace{-10mu}
    \Omega (p)
    \; = \;
    \left \{
    w \in \mathbb{R}^K \, : \,
    \sum_{k \in [K]} w_k \cdot \nu_k (s) \geq \sum_{j \in [J]} r_j (\phi (s, a)) + \gamma \sum_{k \in [K]} w^{t_j}_k \sum_{s' \in \overline{\mathfrak{S}} [\nu_k]} \nu_k (s') \cdot p_{sak} (s') \;\; \forall (s, a) \in \mathcal{C}
    \right\}
\end{equation*}
denotes the feasible region of the first LP in Algorithm~\ref{eq:update}. Here, the inequality holds since $f (w^{t+1}) \leq f (w^t)$ for all $t \in \mathbb{N}$ implies that $f (w^{t_{j+1}}) \leq f (w^{t_j + 1})$ and the equality is due to the definition of $f (w^{t_j + 1})$, respectively.
Since the sequence $\{w^t\}_{t\in\mathbb{N}}$ is bounded, we can construct a compact set $\mathcal{W}$ that contains it.
Proposition~4.4 of \citet{bonnans2013perturbation} as well as the subsequent discussion in the same reference then imply that the optimal value map
\begin{equation*}
    (w, p, \lambda)
    \; \mapsto \;
    \min_{w' \in \Omega(p)} \left\{ f(w') + \lambda \left \lVert w' - w \right \rVert_1 \right\}
\end{equation*}
is continuous over its domain $\mathcal{W} \times \mathcal{P} \times [0, \, \lambda_1]$. Therefore, we have
\begin{equation*}
    \lim_{j \to \infty} \; \min_{w \in \Omega(p^{t_j})}  \left\{ f(w) + \lambda_{t_j + 1} \left \lVert w - w^{t_j} \right \rVert_1 \right\}
    \;\; = \;\;
    \min_{w \in \Omega(p^\infty)} \left\{ f(w) + \lambda_{\infty} \left \lVert w- w^\infty \right \rVert_1 \right\}
    \;\; = \;\;
    \min_{w \in \Omega(p^\infty)} f(w),
\end{equation*}
and taking the limit $j \to \infty$ on both sides of~\eqref{eq:limit_condition} shows that $f(w^\infty) \leq \min \{ f(w) \, : \, w \in \Omega (p^\infty) \}$, that is, $w^\infty$ is indeed optimal in the first LP in Algorithm~\ref{eq:update} when $p^{t-1}$ is replaced with $p^\infty$.

We next show that $p^\infty$ minimizes the second LP in Algorithm~\ref{eq:update} if we replace $w^t$ with $w^\infty$. To this end, we observe that for all $t \in \mathbb{N}$, we have
\begin{equation*}
    \sum_{k\in[K]} \sum_{s' \in \overline{\mathfrak{S}} [\nu_k]} w^t_k \cdot \nu_k (s') \cdot p^t_{sak} (s')
    \;\; \leq \;\;
    \sum_{k\in[K]} \sum_{s' \in \overline{\mathfrak{S}} [\nu_k]} w^t_k \cdot \nu_k (s') \cdot p_{sak} (s')
    \qquad \forall p \in \mathcal{P} (\phi (s, a)),
\end{equation*}
where $\mathcal{P} (\phi (s, a))$ denotes the feasible region of the second LP in Algorithm~\ref{eq:update}. Taking the limit $t \longrightarrow \infty$ on both sides of this inequality implies that
\begin{equation*}
    \sum_{k\in[K]} \sum_{s' \in \overline{\mathfrak{S}} [\nu_k]} w^\infty_k \cdot \nu_k (s') \cdot p^\infty_{sak} (s')
    \;\; \leq \;\;
    \sum_{k\in[K]} \sum_{s' \in \overline{\mathfrak{S}} [\nu_k]} w^\infty_k \cdot \nu_k (s') \cdot p_{sak} (s')
    \qquad \forall p \in \mathcal{P} (\phi (s, a)),
\end{equation*}
which shows that $p^\infty$ minimizes the second LP in Algorithm~\ref{eq:update} if we replace $w^t$ with $w^\infty$. 

In view of the third step, finally, note that we can readily convert the probabilities $p^\infty_{sak} \in \Delta (\overline{\mathfrak{S}} [\nu_k])$, $k \in [K]$, to probabilities $p^\infty_{sa} \in \Delta (\mathcal{S})$ that act as decision variables in problem~\eqref{opt:rfmdp:basis_fcts_lp_simplified}. The KKT conditions of the first LP in Algorithm~\ref{eq:update} then imply the existence of $\Lambda_{sa}$, $(s, a) \in \mathcal{C}$, such that
\begin{equation*}
    \boxed{
    \begin{array}{c}
        \text{\textbf{Stationarity:}} \\
        \displaystyle \sum_{s \in \overline{\mathfrak{S}} [\nu_k]} \nu_k (s) \prod_{n \in \overline{\mathfrak{s}} [\nu_k]} q_n (s_n) + \sum_{(s, a) \in \mathcal{C}} \Lambda_{sa} \left( \gamma \sum_{s' \in \overline{\mathfrak{S}} [\nu_k]} \nu_k (s') \cdot p^\infty_{sa} (s') - \nu_k (s) \right) = 0 \qquad \forall k \in [K] \\
        \text{\textbf{Primal feasibility:}} \\
        \displaystyle \sum_{k \in [K]} w^\infty_k \cdot \nu_k (s) \geq r (\phi (s, a)) + \gamma \sum_{k \in [K]} w^\infty_k \sum_{s' \in \overline{\mathfrak{S}} [\nu_k]} \nu_k (s') \cdot p^\infty_{sa} (s') \qquad \forall (s, a) \in \mathcal{C} \\
        \text{\textbf{Dual feasibility:}} \\
        \displaystyle \Lambda_{sa} \in \mathbb{R}_+ \qquad \forall (s, a) \in \mathcal{C} \\
        \text{\textbf{Complementary slackness:}} \\
        \displaystyle \Lambda_{sa} \left( r (\phi (s, a)) + \gamma \sum_{k \in [K]} w^\infty_k \sum_{s' \in \overline{\mathfrak{S}} [\nu_k]} \nu_k (s') \cdot p^\infty_{sa} (s') - \sum_{k \in [K]} w^\infty_k \cdot \nu_k (s) \right) = 0 \qquad \forall (s, a) \in \mathcal{C}.
    \end{array}}
\end{equation*}
Similarly, the KKT conditions of the second LP in Algorithm~\ref{eq:update} imply the existence of $\Theta_{sa}$, $(s, a) \in \mathcal{C}$, such that
\begin{equation*}
    \boxed{
    \begin{array}{c}
        \text{\textbf{Stationarity:}} \\
        \displaystyle \gamma \cdot \Lambda_{sa} \sum_{k \in [K]} w^\infty_k \cdot \nu_k (s') + [\Psi_{sa}]_{s'}{}^\top \Theta_{sa} = 0 \qquad \forall (s, a) \in \mathcal{C}, \; \forall s' \in \mathcal{S} \\
        \text{\textbf{Primal feasibility:}} \\
        \displaystyle \Psi_{sa} p^\infty_{sa} \leq \psi_{sa} \qquad \forall (s, a) \in \mathcal{C} \\
        \text{\textbf{Dual feasibility:}} \\
        \displaystyle \Theta_{sa} \in \mathbb{R}^{ M_{sa} }_+ \qquad \forall (s, a) \in \mathcal{C} \\
        \text{\textbf{Complementary slackness:}} \\
        \displaystyle \Theta_{sa}{}^\top \left( \Psi_{sa} p^\infty_{sa} - \psi_{sa} \right) = 0 \qquad \forall (s, a) \in \mathcal{C}.
    \end{array}}
\end{equation*}
We can now see that $(w^\infty, p^\infty, \Lambda, \Gamma)$ with $\Gamma_{sa} = \gamma \Lambda_{sa} \Theta_{sa}$, $(s, a) \in \mathcal{C}$, satisfies the KKT conditions of the restricted master problem~\eqref{opt:rfmdp:basis_fcts_lp_simplified} as reported in the proof of Theorem~\ref{thm:lambda-epsilon-kkt-point}.
\qed

~\\[-8mm]

\noindent \textbf{Proof of Corollary~\ref{coro:robust:solution_subproblem}.} $\;$
    To determine a maximally violated constraint in problem~\eqref{opt:rfmdp:basis_fcts_lp_simplified}, we need to solve the max-min problem from Algorithm~\ref{alg:CP-robust}. We follow the standard robust optimization approach and first dualize the inner minimization problem,
    \begin{equation*}
        \begin{array}{l@{\quad}l@{\qquad}l}
            \displaystyle \mathop{\text{minimize}}_{p} & \displaystyle \sum_{k\in[K]} w^\star_k \sum_{s' \in \overline{\mathfrak{S}} [\nu_k]} \nu_k (s') \cdot p_{sak} (s') \\
            \displaystyle \text{subject to} & \displaystyle \mathrm{Marg}_n (p_{sak}) \in \mathcal{P}^n (\phi (s, a)) & \displaystyle \forall k \in [K], \; \forall n \in \overline{\mathfrak{s}} [\nu_k] \\
            & \displaystyle \mathrm{Marg}_{\mathcal{N} (k, k')} (p_{sak}) = \mathrm{Marg}_{\mathcal{N} (k, k')} (p_{sak'}) & \displaystyle \forall 1 \leq k < k' \leq K \\
            & \displaystyle p_{sak} \in \Delta (\overline{\mathfrak{S}} [\nu_k]), \, k \in [K].
        \end{array}
    \end{equation*}
    This problem can be re-expressed as
    \begin{equation}\label{eq:robust_subproblem:primal}
        \begin{array}{l@{\quad}l@{\qquad}l}
            \displaystyle \mathop{\text{minimize}}_{p} & \displaystyle \sum_{k \in [K]} w^\star_k \cdot \mu_k{}^\top p_{sak} \\
            \displaystyle \text{subject to} & \displaystyle \Psi_n (\phi(s,a)) \mathrm{M}_{k, n} p_{sak} \leq \psi_n (\phi(s,a)) & \displaystyle \forall k \in [K], \; \forall n \in \overline{\mathfrak{s}} [\nu_k] \\
            & \displaystyle \mathrm{M}_{k, \mathcal{N} (k, k')} p_{sak} = \mathrm{M}_{k', \mathcal{N} (k, k')} p_{sak'} & \displaystyle \forall 1 \leq k < k' \leq K \\
            & \displaystyle p_{sak} \geq 0, \, k \in [K],
        \end{array}
    \end{equation}
    where we introduce $| \overline{\mathfrak{S}} [\nu_k] |$-dimensional vectors $\mu_k$ satisfying $\mu_k{}^\top p_{sak} = \sum_{s' \in \overline{\mathfrak{S}} [\nu_k]} \nu_k (s') \cdot p_{sak} (s')$, where we describe the ambiguity set via $\mathcal{P}^n (\phi (s, a)) = \{ p_n \in \mathbb{R}^{S_n}_+ : \Psi_n (\phi(s,a)) p_n \leq \psi_n (\phi(s,a)) \}$ for some $\Psi_n (\phi(s,a)) \in \mathbb{R}^{M_{sa} \times S_n}$ and $\psi_n (\phi(s,a)) \in \mathbb{R}^{M_{sa}}$, and where we exploit that $\mathrm{Marg}_n (p_{sa \kappa}) = \mathrm{M}_{\kappa, \{ n \}} p_{sa \kappa}$ and $\mathrm{Marg}_{\mathcal{N} (k, k')} (p_{sa \kappa}) = \mathrm{M}_{\kappa, \mathcal{N} (k, k')} p_{sa \kappa}$ for
    \begin{equation*}
        \mathrm{M}_{k,\Omega} = \mathop{\otimes}\limits_{n \in \overline{\mathfrak{s}} [\nu_k]} \mathrm{E}_{k,\Omega, n} \in \mathbb{R}^{M_{\Omega} \times | \overline{\mathfrak{S}} [\nu_k] |}
        \quad \text{with} \quad
        E_{k,\Omega, n} = \begin{cases}
            \mathbf{I}_{2^{S_n}} & \text{if } n \in \Omega, \\
            \mathbf{1}_{2^{S_n}}{}^\top & \text{otherwise,}
        \end{cases}
    \end{equation*}
    where $\mathbf{I}_d$ and $\mathbf{1}_d$ are the $d \times d$ identity matrix and the $d$-vector of all ones, respectively, $\otimes$ denotes the Kronecker product, and $M_{\Omega} = 2^{\sum_{n \in \Omega} S_n}$. Here and in the following, we omit the brackets around $n$ in the subscript of $\mathrm{M}_{k, \{ n \}}$ for ease of notation.

    Since $\mathcal{P} (\phi (s, a))$ is non-empty and bounded by construction, problem~\eqref{eq:robust_subproblem:primal} is feasible and bounded, and it therefore affords the strong LP dual
    \begin{equation*}
        \begin{array}{l@{\quad}l@{\qquad}l}
            \displaystyle \mathop{\text{maximize}}_{\chi, \, \delta} & \displaystyle - \sum_{k \in [K]} \sum_{n \in \overline{\mathfrak{s}} [\nu_k]} \psi_n (\phi(s,a))^\top \delta_{kn} \\
            \displaystyle \text{subject to} & \displaystyle w^\star_k \mu_k + \sum_{n \in \overline{\mathfrak{s}} [\nu_k]} \mathrm{M}_{k, n}{}^\top \Psi_n (\phi(s,a))^\top \delta_{kn} \geq \\
            & \displaystyle \mspace{60mu} \sum_{1 \leq k' < k} \mathrm{M}_{k, \mathcal{N} (k, k')}{}^\top \chi_{k'k} - \sum_{k < k' \leq K} \mathrm{M}_{k, \mathcal{N} (k, k')}{}^\top \chi_{kk'} & \displaystyle \forall k \in [K]\\
            & \multicolumn{2}{l}{\displaystyle \chi_{kk'} \in \mathbb{R}^{M_{\mathcal{N} (k, k')}}, \, 1 \leq k < k' \leq K, \;\; \delta_{kn} \in \mathbb{R}^{M_{sa}}_+, \, k \in [K] \text{ and } n \in \overline{\mathfrak{s}} [\nu_k].}
        \end{array}
    \end{equation*}
    Note that the functions $\Psi_n$ and $\psi_n$ inherit the low-scope property of $\mathcal{P}^n$, and let $\mathfrak{s} [\mathcal{P}^n] = \mathfrak{s} [\Psi_n] \cup \mathfrak{s} [\psi_n]$ and $\mathfrak{S} [\mathcal{P}^n] = \{ \sum_{i \in [I]} \mathbf{1} [i \in \mathfrak{s} [\mathcal{P}^n] ] \cdot x_i \cdot \mathrm{e}_i \, : \, x \in \mathcal{X} \}$ describe the scope of $\mathcal{P}^n$ in terms of indices and sub-vectors, respectively. We then obtain that
    \begin{equation*}
        \Psi_n (\phi(s,a)) = \sum_{f \in \mathfrak{S} [\mathcal{P}^n]} \Psi_n (f) \cdot \mathbf{1} [f = \phi (s, a)]
        \quad \text{ and } \quad
        \psi_n (\phi(s,a)) = \sum_{f \in \mathfrak{S} [\mathcal{P}^n]} \psi_n (f) \cdot \mathbf{1} [f = \phi (s, a)].
    \end{equation*}
    Following the same reasoning as in the proof of Theorem~\ref{thm:solution_subproblem}, we have $\xi_{nf} = \mathbf{1} [\varphi_i = f_i \; \forall i \in \mathfrak{s} [\mathcal{P}^n]]$ if and only if $\xi_{nf} \in [0, 1]$ and
    \begin{equation*}
        \xi_{nf} \leq (2 f_i - 1) \varphi_i + 1 - f_i \;\; \forall i \in \mathfrak{s} [\mathcal{P}^n]
        \quad \text{as well as} \quad
        \xi_{nf} \geq 1 + \sum_{i \in \mathfrak{s} [\mathcal{P}^n]} \frac{\varphi_i - f_i}{2 f_i - 1}.
    \end{equation*}
    Hence, the dual problem has the equivalent reformulation
    \begin{equation*}
        \mspace{-45mu}
        \begin{array}{l@{\quad}l@{\quad}l}
            \displaystyle \mathop{\text{maximize}}_{\chi, \, \delta, \, \xi} & \displaystyle - \sum_{k \in [K]} \sum_{n \in \overline{\mathfrak{s}} [\nu_k]} \sum_{f \in \mathfrak{S} [\mathcal{P}^n]} \psi_n (f)^\top \delta_{kn} \cdot \xi_{nf} \\
            \displaystyle \text{subject to} & \displaystyle w^\star_k \mu_k + \sum_{n \in \overline{\mathfrak{s}} [\nu_k]} \sum_{f \in \mathfrak{S} [\mathcal{P}^n]} \mathrm{M}_{k, n}{}^\top \Psi_n (f)^\top \delta_{kn} \cdot \xi_{nf} \geq \\
            & \displaystyle\displaystyle \mspace{60mu} \sum_{1 \leq k' < k} \mathrm{M}_{k, \mathcal{N} (k, k')}{}^\top \chi_{k'k} - \sum_{k < k' \leq K} \mathrm{M}_{k, \mathcal{N} (k, k')}{}^\top \chi_{kk'} & \displaystyle \forall k \in [K] \\
            & \displaystyle \xi_{nf} \in [0, 1], \;\; 1 + \sum_{i \in \mathfrak{s} [\mathcal{P}^n]} \frac{\varphi_i - f_i}{2 f_i - 1} \leq \xi_{nf} \leq (2 f_l - 1) \varphi_l + 1 - f_l & \displaystyle \forall n \in [N], \; \forall f \in \mathfrak{S} [\mathcal{P}^n], \; \forall l \in \mathfrak{s} [\mathcal{P}^n] \\
            & \multicolumn{2}{l}{\mspace{-8mu} \displaystyle \chi_{kk'} \in \mathbb{R}^{M_{\mathcal{N} (k, k')}}, \, 1 \leq k < k' \leq K, \;\; \delta_{kn} \in \mathbb{R}^{M_{sa}}_+, \, k \in [K] \text{ and } n \in \overline{\mathfrak{s}} [\nu_k],}
        \end{array}
    \end{equation*}
    where we use the abbreviation $\varphi = \phi (s, a)$. We next linearize this problem by replacing both the decision vectors $\delta_{kn}$ and the bilinear terms $\delta_{kn} \cdot \xi_{nf}$ with triple-indexed auxiliary vectors $\delta_{knf}$ that satisfy $\delta_{knf} \geq 0$ and $\delta_{knf \ell} \leq \mathrm{M} \cdot \xi_{nf}$, $\ell \in [M_{sa}]$, for some sufficiently large constant $\mathrm{M} > 0$. Under this reformulation, the dual problem becomes
    \begin{equation*}
        \mspace{-90mu}
        \begin{array}{l@{\quad}l@{\quad}l}
            \displaystyle \mathop{\text{maximize}}_{\chi, \, \delta, \, \xi} & \displaystyle - \sum_{k \in [K]} \sum_{n \in \overline{\mathfrak{s}} [\nu_k]} \sum_{f \in \mathfrak{S} [\mathcal{P}^n]} \psi_n (f)^\top \delta_{knf} \\
            \displaystyle \text{subject to} & \displaystyle w^\star_k \mu_k + \sum_{n \in \overline{\mathfrak{s}} [\nu_k]} \sum_{f \in \mathfrak{S} [\mathcal{P}^n]} \mathrm{M}_{k, n}{}^\top \Psi_n (f)^\top \delta_{knf}  \geq \\
            & \displaystyle \mspace{60mu} \sum_{1 \leq k' < k} \mathrm{M}_{k, \mathcal{N} (k, k')}{}^\top \chi_{k'k} - \sum_{k < k' \leq K} \mathrm{M}_{k, \mathcal{N} (k, k')}{}^\top \chi_{kk'} & \displaystyle \forall k \in [K] \\
            & \displaystyle \xi_{nf} \in [0, 1], \;\; 1 + \sum_{i \in \mathfrak{s} [\mathcal{P}^n]} \frac{\varphi_i - f_i}{2 f_i - 1} \leq \xi_{nf} \leq (2 f_l - 1) \varphi_l + 1 - f_l & \displaystyle \forall n \in [N], \; \forall f \in \mathfrak{S} [\mathcal{P}^n], \; \forall l \in \mathfrak{s} [\mathcal{P}^n] \\
            & \multicolumn{2}{l}{\mspace{-8mu} \displaystyle \chi_{kk'} \in \mathbb{R}^{M_{\mathcal{N} (k, k')}}, \, 1 \leq k < k' \leq K, \;\;  \delta_{knf} \in \mathbb{R}^{M_{sa}}_+, \, \delta_{knf \ell} \leq \mathrm{M} \cdot \xi_{nf}, \, k \in [K], \, n \in \overline{\mathfrak{s}} [\nu_k], f \in \mathfrak{S} [\mathcal{P}^n] \text{ and } \ell \in [M_{sa}].}
        \end{array}
    \end{equation*}
    Finally, re-inserting this problem into the outer maximization of the max-min problem from Algorithm~\ref{alg:CP-robust} and applying the same transformations as in the proof of Theorem~\ref{thm:solution_subproblem} results in the formulation from the statement of the corollary.

    Following the standard theory of LP duality, the transition probabilities $p_{s^\star a^\star}$ associated with an optimal solution $(s^\star, a^\star, \varphi^\star, \zeta^\star, \eta^\star, \xi^\star, \beta^\star, \chi^\star, \delta^\star)$ of the MILP from the statement of the corollary correspond to the shadow prices of the first $\sum_{k \in [K]} | \overline{\mathfrak{s}} [\nu_k] |$ constraints in that problem.
\qed

~\\[-8mm]

The worst-case expectation in the statement of Proposition~\ref{prop:factored_ambiguity_hard} can be formulated as a bilinear program, and bilinear programming is known to be NP-hard. However, our worst-case expectation problem has a specific structure where each set of variables $p_n$ resides in a probability simplex. This structure appears to render straightforward reductions onto NP-hard classes of bilinear programs, such as the quadratic assignment problem \citep{SG76:qap} or bilinear set separation \citep{BM93:bilinear}, difficult. Instead, our proof of Proposition~\ref{prop:factored_ambiguity_hard} combines a classical representation of Nash equilibria due to \cite{MS64:nash} with a recent PPAD-hardness result of \cite{CD06:ppad}.

~\\[-8mm]

\noindent \textbf{Proof of Proposition~\ref{prop:factored_ambiguity_hard}.}
    We construct a reduction to the two-player Nash equilibrium problem:
    \begin{center}
        \fbox{\parbox{13cm}{ {\centering \textsc{Two-Player Nash Equilibrium.} \\}
        \textbf{Instance.} Given are payoff matrices $A, B \in \mathbb{R}^{m_1 \times m_2}$ for two players. \\
        \textbf{Question.} Find a profile of mixed strategies $x^\star \in \Delta([m_1])$ and $y^\star \in \Delta([m_2])$ that constitute mutual best responses, that is, $(x^\star)^\top A y^\star \geq x^\top A y^\star$ for all $x \in \Delta([m_1])$ and $(x^\star)^\top B y^\star \geq (x^\star)^\top B y$ for all $y \in \Delta([m_2])$.}}
    \end{center}
    The two-player Nash equilibrium problem is known to be PPAD-hard \citep{CD06:ppad}. Moreover, \cite{MS64:nash} have shown that the Nash equilibria in this problem coincide precisely with the optimal solutions $(x^\star, y^\star, \alpha^\star, \beta^\star)$ to the bilinear program
    \begin{equation}\label{eq:ppad_hard1}
        \begin{array}{l@{\quad}l}
            \displaystyle \mathop{\text{maximize}}_{x, \, y, \, \alpha, \, \beta} & \displaystyle x^\top (A + B) y - \alpha - \beta \\
            \displaystyle \text{subject to} & \displaystyle \alpha \cdot \mathbf{1} \geq Ay, \;\; \beta \cdot \mathbf{1} \geq B^\top x \\
            & \displaystyle x \in \Delta([m_1]), \;\; y \in \Delta([m_2]), \;\; \alpha, \beta \in \mathbb{R},
        \end{array}
    \end{equation}
    and thus solving problem~\eqref{eq:ppad_hard1} is also PPAD-hard. Note that any optimal solution $(x^\star, y^\star, \alpha^\star, \beta^\star)$ to problem~\eqref{eq:ppad_hard1} satisfies $\alpha^\star = a_{i_1}{}^\top y^\star = \max_{i \in [m_1]} \, a_i{}^\top y^\star$ and $\beta^\star = b_{i_2}{}^\top x^\star = \max_{i \in [m_2]} \, b_i{}^\top x^\star$ for some $(i_1, i_2) \in [m_1] \times [m_2]$, where $a_i{}^\top$ and $b_j$ refer to the $i$-th row of matrix $A$ and the $j$-th column of matrix $B$, respectively. We can thus optimally solve problem~\eqref{eq:ppad_hard1} by solving the $m_1 \cdot m_2$ auxiliary bilinear programs
    \begin{equation}\label{eq:ppad_hard2}
        \begin{array}{l@{\quad}l@{\qquad}l}
            \displaystyle \mathop{\text{maximize}}_{x, \, y} & \multicolumn{2}{l}{\mspace{-8mu} \displaystyle x^\top (A + B) y - a_{i_1}{}^\top y - b_{i_2}{}^\top x} \\
            \displaystyle \text{subject to} & \displaystyle a_{i_1}{}^\top y \geq a_i{}^\top y & \displaystyle \forall i \in [m_1] \\
            & \displaystyle b_{i_2}{}^\top x \geq b_i{}^\top x & \displaystyle \forall i \in [m_2] \\
            & \multicolumn{2}{l}{\mspace{-8mu} \displaystyle x \in \Delta([m_1]), \;\; y \in \Delta([m_2])}
        \end{array}
    \end{equation}
    parameterized by $(i_1, i_2) \in [m_1] \times [m_2]$, where we have substituted out the decision variables $\alpha$ and $\beta$. In particular, an optimal solution $(x^\star, y^\star)$ to an instance $(i_1^\star, i_2^\star)$ of problem~\eqref{eq:ppad_hard2} achieving the largest objective value among all instances $(i_1, i_2) \in [m_1] \times [m_2]$ gives rise to an optimal solution $(x^\star, y^\star, \alpha^\star, \beta^\star)$ via $\alpha^\star = a_{i_1^\star}{}^\top y^\star$ and $\beta^\star = b_{i_2^\star}{}^\top x^\star$.
    Due to their ability to jointly solve the PPAD-hard problem~\eqref{eq:ppad_hard1}, solving the bilinear programs~\eqref{eq:ppad_hard2} is PPAD-hard as well. In the remainder, we show that each bilinear program~\eqref{eq:ppad_hard2} can be cast as an instance of our worst-case expectation problem from the statement of the proposition. This will imply that our worst-case expectation problem is PPAD-hard as well, as desired.

    Fix an instance of problem~\eqref{eq:ppad_hard2}, which is described by the matrices $A, B \in \mathbb{R}^{m_1 \times m_2}$ as well as the indices $i_1 \in [m_1]$ and $i_2 \in [m_2]$. We construct a factored MDP with $N = 2$ sub-state spaces $\mathcal{S}_n = \{ s_n \in \mathbb{B}^{m_n} \, : \sum_{i \in [m_n]} s_{n,i} = 1 \}$, $n \in [N]$, and the ambiguity sets $\mathcal{P}^1 = \{ p \in \Delta ([m_1]) \, : \, b_{i_2}{}^\top p \geq b_i{}^\top p \;\; \forall i \in [m_2] \}$ as well as $\mathcal{P}^2 = \{ p \in \Delta ([m_2]) \, : \, a_{i_1}{}^\top p \geq a_i{}^\top p \;\; \forall i \in [m_1] \}$. According to our definition of the factored ambiguity sets, the transition probability $p (s')$ for $s' = (s'_1, s'_2) \in \mathcal{S}$ satisfies $p (s') = p_1 (s'_1) \cdot p_2 (s'_2)$ for some $p_1 \in \mathcal{P}^1$ and $p_2 \in \mathcal{P}^2$. In the following, $p_1$ and $p_2$ will play the roles of $x$ and $y$ in problem~\eqref{eq:ppad_hard2}, respectively. Setting $v (s') = A_{j_1, j_2} + B_{j_1, j_2} - a_{i_1, j_2} - b_{i_2, j_1}$ for all such $s'$ then results in the expected value
    \begin{align*}
        \sum_{s' \in \mathcal{S}} p (s') \cdot v (s')
        \;\; &= \;\;
        \sum_{j_1 \in [m_1]} \sum_{j_2 \in [m_2]} (A_{j_1, j_2} + B_{j_1, j_2} - a_{i_1, j_2} - b_{i_2, j_1}) \cdot p_1 (s'_{1, j_1}) \cdot p_2 (s'_{2, j_2}) \\
        &= \;\;
        p_1{}^\top (A + B) p_2 - a_{i_1}{}^\top p_2 - b_{i_2}{}^\top p_1.
    \end{align*}
    The statement now follows from the fact that the constraint $p \in \mathcal{P} (\phi (s, a))$ in the worst-case expectation problem coincides precisely with the constraints of problem~\eqref{eq:ppad_hard2}.
\qed

~\\[-8mm]

\noindent \textbf{Proof of Observation~\ref{obs:lifted_policies}.}
In view of the first claim, we first show via contraposition that $\hat{\Pi} (\pi) \cap \hat{\Pi} (\pi') = \emptyset$ whenever $\pi \neq \pi'$. To this end, fix any two non-lifted policies $\pi, \pi' \in \Pi$ that satisfy $\hat{\Pi} (\pi) \cap \hat{\Pi} (\pi') \neq \emptyset$ as well as any $\hat{\pi} \in \hat{\Pi} (\pi) \cap \hat{\Pi} (\pi')$. We then have
\begin{equation*}
    \hat{\pi} (\hat{s})
    \; = \;
    \pi (u)
    \; = \;
    \pi' (u)
    \qquad \forall \hat{s} = (u, s^0, a^0, 1) \in \hat{\mathcal{S}}.
\end{equation*}
Since $\{ u : (u, s^0, a^0, 1) \in \hat{\mathcal{S}} \} = \mathcal{S}$, we have $\pi (u) = \pi' (u)$ for all $u \in \mathcal{S}$, that is, we obtain the desired contradiction that $\pi = \pi'$. The first claim now follows if we can show that $\bigcup_{\pi \in \Pi} \hat{\Pi} (\pi) = \hat{\Pi}$. Indeed, fix any $\hat{\pi} \in \hat{\Pi}$ and construct the policy $\pi$ satisfying $\pi (u) = \hat{\pi} ((u, s^0, a^0, 1))$ for all $u \in \mathcal{S}$. We then have $\hat{\pi} \in \hat{\Pi} (\pi)$ as desired.

As for the second claim, fix any $\hat{\pi}, \hat{\pi}' \in \hat{\Pi}$ satisfying $\hat{\pi}, \hat{\pi}' \in \hat{\Pi} (\pi)$ for some $\pi \in \Pi$. This implies that $\hat{\pi} (\hat{s}) = \hat{\pi}' (\hat{s})$ for all states $\hat{s} = (u, m, a, n) \in \hat{\mathcal{S}}$ with $(m, a, n) = (s^0, a^0, 1)$. Fix any transition kernel $\hat{p} \in \hat{\mathcal{P}}$
as well as any finite state trajectory $\hat{s}^1, \ldots, \hat{s}^{T + 1} \in \hat{\mathcal{S}}$, $T \in \mathbb{N}$, and set $M = \lfloor T / N \rfloor$. The probability of observing this trajectory under $\hat{p}$ coincides for both policies $\hat{\pi}$ and $\hat{\pi}'$ since
\begin{align*}
    & q (\hat{s}^1) \cdot \prod_{t \in [T]} \hat{p} (\hat{s}^{t+1} \, | \, \hat{s}^t, \hat{\pi} (\hat{s}^t)) \\
    = \; &
    q (\hat{s}^1) \cdot \prod_{\tau = 0}^M \prod_{t \in [N]} \hat{p} (\hat{s}^{\tau N + t + 1} \, | \, \hat{s}^{\tau N + t}, \hat{\pi} (\hat{s}^{\tau N + t}))
    \; \cdot \;
    \prod_{t \in [T - MN]} \hat{p} (\hat{s}^{MN + t + 1} \, | \, \hat{s}^{MN + t}, \hat{\pi} (\hat{s}^{MN + t})) \\
    = \; &
    q (\hat{s}^1) \cdot \prod_{\tau = 0}^M \prod_{t \in [N]} \hat{p} (\hat{s}^{\tau N + t + 1} \, | \, \hat{s}^{\tau N + t}, \hat{\pi}' (\hat{s}^{\tau N + t}))
    \; \cdot \;
    \prod_{t \in [T - MN]} \hat{p} (\hat{s}^{MN + t + 1} \, | \, \hat{s}^{MN + t}, \hat{\pi}' (\hat{s}^{MN + t})) \\
    = \; &
    q (\hat{s}^1) \cdot \prod_{t \in [T]} \hat{p} (\hat{s}^{t+1} \, | \, \hat{s}^t, \hat{\pi}' (\hat{s}^t)).
\end{align*}
Here, the first and last identity group the transitions into $M$ sub-sequences of $N$ transitions each, plus one final sub-sequence collecting the last $T - MN < N$ transitions (if present). To see that the second equality holds as well, observe first that we either have $\hat{s}^{\tau N + 1} = (u, s^0, a^0, 1)$ for some $u \in \mathcal{S}$ or the probability of observing the trajectory is zero (in which case the entire chain of equalities trivially holds). A similar argument shows that we either have $\hat{s}^{\tau N + t} = (u, m, a, t)$ for some $u, m \in \mathcal{S}$ and $a \in \mathcal{A}$ or the probability of observing the trajectory is zero. The second equality then follows from a case distinction on $t$. Whenever $t = 1$, we can assume that $\hat{s}^{\tau N + t} = (u, s^0, a^0, 1)$ for some $u \in \mathcal{S}$, and we observed before that $\hat{\pi} (\hat{s}) = \hat{\pi}' (\hat{s})$ for all such states, that is, the probabilities on both sides of the second equality must coincide for $t = 1$. Whenever $t > 1$, on the other hand, the construction of our lift implies that the transition probabilities are unaffected by the policy, and thus the probabilities on both sides of the second equality trivially coincide for $t > 1$.

To conclude the second claim, we show that the two policies $\hat{\pi}$ and $\hat{\pi}'$ generate the same rewards on the state trajectory $\hat{s}^1, \ldots, \hat{s}^{T + 1}$. To this end, we can again assume that $\hat{s}^{\tau N + 1} = (u, s^0, a^0, 1)$ for some $u \in \mathcal{S}$ and $\hat{s}^{\tau N + t} = (u, m, a, t)$ for some $u, m \in \mathcal{S}$ and $a \in \mathcal{A}$. We then observe that
\begin{align*}
    & \sum_{t \in [T + 1]} \hat{\gamma}^t \cdot \hat{r} (\hat{\phi} (\hat{s}^t, \hat{\pi} (\hat{s}^t))) \\
    = \; &
    \sum_{\tau = 0}^M \sum_{t \in [N]} \hat{\gamma}^{\tau N + t} \cdot \hat{r} (\hat{\phi} (\hat{s}^{\tau N + t}, \hat{\pi} (\hat{s}^{\tau N + t})))
    \; + \;
    \sum_{t \in [T + 1 - MN]} \hat{\gamma}^{MN + t} \cdot \hat{r} (\hat{\phi} (\hat{s}^{MN + t}, \hat{\pi} (\hat{s}^{MN + t}))) \\
    = \; &
    \sum_{\tau = 0}^M \hat{\gamma}^{\tau N} \cdot \hat{r} (\hat{\phi} (\hat{s}^{\tau N}, \hat{\pi} (\hat{s}^{\tau N})))
    \;\; = \;\;
    \sum_{\tau = 0}^M \hat{\gamma}^{\tau N} \cdot \hat{r} (\hat{\phi} (\hat{s}^{\tau N}, \hat{\pi}' (\hat{s}^{\tau N}))) \\
    = \; &
    \sum_{\tau = 0}^M \sum_{t \in [N]} \hat{\gamma}^{\tau N + t} \cdot \hat{r} (\hat{\phi} (\hat{s}^{\tau N + t}, \hat{\pi}' (\hat{s}^{\tau N + t})))
    \; + \;
    \sum_{t \in [T + 1 - MN]} \hat{\gamma}^{MN + t} \cdot \hat{r} (\hat{\phi} (\hat{s}^{MN + t}, \hat{\pi}' (\hat{s}^{MN + t}))) \\
    = \; &
    \sum_{t \in [T + 1]} \hat{\gamma}^t \cdot \hat{r} (\hat{\phi} (\hat{s}^t, \hat{\pi}' (\hat{s}^t))).
\end{align*}
Here, the first and last identity group the transitions into sub-sequences as before. The second and penultimate identity exploit that the rewards vanish at all states $\hat{s} = (u, m, a, n)$ with $(m, a, n) \neq (s^0, a^0, 1)$. The third identity, finally, follows from the fact that $\hat{\pi} (\hat{s}^{\tau N}) = \pi (u) = \hat{\pi}' (\hat{s}^{\tau N})$ for $\hat{s}^{\tau N} = (u, s^0, a^0, 1)$ with $u \in \mathcal{S}$ since $\hat{\pi}, \hat{\pi}' \in \hat{\Pi} (\pi)$. We thus conclude that $\hat{\pi}$ and $\hat{\pi}'$ generate the same transitions and rewards across any finite state trajectory, which implies the second claim.
\qed

~\\[-8mm]

\noindent \textbf{Proof of Theorem~\ref{thm:rel_ambiguity_sets}.}
    In view of the first statement, we first show that $\mathrm{F} (\pi; \mathcal{R}) \geq \mathrm{L} (\hat{\pi}; \mathcal{R})$ for all instances $\mathcal{R}$, all non-lifted policies $\pi \in \Pi$ and all lifted policies $\hat{\pi} \in \hat{\Pi} (\pi)$. To this end, fix any $\mathcal{R}$, $\pi \in \Pi$, $\hat{\pi} \in \hat{\Pi} (\pi)$ and $p \in \mathcal{P}^{\mathrm{F}}$, and construct the lifted transition kernel $\hat{p} \in \hat{\mathcal{P}}$ associated with $p$ as outlined in the beginning of Section~\ref{sec:robust_FMDP:factored}. Our claims then follows if we show that
    \begin{equation*}
        \mathbb{E}_p \Big[ \sum_{t \in \mathbb{N}_0} \gamma^t \cdot r (\tilde{s}^t, \pi (\tilde{s}^t)) \, | \, \tilde{s}^0 \sim q \Big]
        \;\; = \;\;
        \mathbb{E}_{\hat{p}} \Big[ \sum_{t \in \mathbb{N}_0} \hat{\gamma}^t \cdot \hat{r} (\tilde{\hat{s}}^t, \hat{\pi} (\tilde{\hat{s}}^t)) \, | \, \tilde{\hat{s}}^0 \sim \hat{q} \Big],
    \end{equation*}
    where we notationally emphasize the transition kernels governing the respective stochastic processes. To this end, we proceed along the lines of the proof of Observation~\ref{obs:lifted_policies} and notice that there is a one-to-one correspondence between the state trajectories $s^1, \ldots, s^T \in \mathcal{S}$, $T \in \mathbb{N}$, in the factored RFMDP and the state trajectories $\hat{s}^1, \ldots, \hat{s}^{NT} \in \hat{S}$ in the lifted RFMDP that satisfy
    \begin{equation*}
        \hat{s}^{N(t-1) + n} \; = \;
        \left(
            \left[ \begin{matrix}
                (s^{t + 1}_1, \ldots, s^{t + 1}_n)^\top \\
                (s^t_{n+1}, \ldots, s^t_N)^\top
            \end{matrix} \right], \;
            s^t(n), \;
            a^t(n), \;
            n
        \right)
        \qquad \forall t \in [T], \; \forall n \in [N],
    \end{equation*}
    where $(s^t (1), a^t (1)) = (s^0, a^0)$ and $(s^t (n), a^t (n)) = (s^t, \pi(s^t))$ otherwise, in the sense that both trajectories share the same occurrence probabilities and rewards, and no other trajectories are observed with positive probability in either RFMDP.
    
    We next show that $\mathrm{L} (\hat{\pi}; \mathcal{R}) \geq \mathrm{NF} (\pi; \mathcal{R})$ for all instances $\mathcal{R}$, all non-lifted policies $\pi \in \Pi$ and all lifted policies $\hat{\pi} \in \hat{\Pi} (\pi)$. To this end, fix any $\mathcal{R}$, $\pi \in \Pi$, $\hat{\pi} \in \hat{\Pi} (\pi)$ and $\hat{p} \in \hat{\mathcal{P}}$, and consider the non-lifted and non-factored transition kernel $p$ satisfying
    \begin{equation*}
        \mspace{-25mu}
        p_{sa} (s') \; = \;
        \prod_{n \in [N]}
        \hat{p} \left(
        \left(
            \left[ \begin{matrix}
                (s'_1, \ldots, s'_n) \\
                (s_{n+1}, \ldots, s_N)^\top
            \end{matrix} \right],
            s(n \oplus 1), a(n \oplus 1), n \oplus 1
        \right) \; \Bigg| \;
        \left(
            \left[ \begin{matrix}
                (s'_1, \ldots, s'_{n-1}) \\
                (s_n, \ldots, s_N)^\top
            \end{matrix} \right],
            s(n), a(n), n
        \right) \right)
    \end{equation*}
    for all $s, s' \in \mathcal{S}$ and $a \in \mathcal{A}$, where $N \oplus 1 = 1$ and $n \oplus 1 = n + 1$ otherwise, $(s(1), a(1)) = (s^0, a^0)$ and $(s(n), a(n)) = (s, a)$ otherwise. The assumed convexity of the marginal ambiguity sets ensures that $p \in \mathcal{P}^\mathrm{NF}$. A similar argument as in the preceding paragraph shows that
    \begin{equation*}
        \mathbb{E}_{\hat{p}} \Big[ \sum_{t \in \mathbb{N}_0} \hat{\gamma}^t \cdot \hat{r} (\tilde{\hat{s}}^t, \hat{\pi} (\tilde{\hat{s}}^t)) \, | \, \tilde{\hat{s}}^0 \sim \hat{q} \Big]
        \;\; = \;\;
        \mathbb{E}_p \Big[ \sum_{t \in \mathbb{N}_0} \gamma^t \cdot r (\tilde{s}^t, \pi (\tilde{s}^t)) \, | \, \tilde{s}^0 \sim q \Big],
    \end{equation*}
    which concludes that indeed $\mathrm{NF} (\pi; \mathcal{R}) \leq \mathrm{L} (\hat{\pi}; \mathcal{R}) \leq \mathrm{F} (\pi; \mathcal{R})$ for all $\mathcal{R}$, $\pi \in \Pi$ and $\hat{\pi} \in \hat{\Pi} (\pi)$.

    The case where $\mathrm{NF} (\pi; \mathcal{R}) = \mathrm{L} (\hat{\pi}; \mathcal{R}) = \mathrm{F} (\pi; \mathcal{R})$ arises, for example, when all marginal ambiguity sets are singleton sets that contain Dirac distributions, that is, when we consider a non-robust deterministic FMDP. The examples provided below to prove the second and third statement of the theorem constitute instances where $\mathrm{L} (\hat{\pi}; \mathcal{R})$ coincides with either $\mathrm{NF} (\pi; \mathcal{R})$ or $\mathrm{F} (\pi; \mathcal{R})$ but not both. Example~\ref{ex:factored_nonfactored}, finally, provides an instance where all three formulations attain different values. This proves the first statement of the theorem.

    As for the second statement, consider the family of RFMDP instances with state space $\mathcal{S} = \prod_{n \in [N]} \mathcal{S}_n$ and $\mathcal{S}_n = \mathbb{B}$, $n \in [N]$, action space $\mathcal{A} = \{ 1 \}$, $\phi (s, a) = (s, a)$, a deterministic initial state $s^0$ satisfying $s^0_n = \mathbf{1} [n \text{ is even}]$, $n \in [N]$, marginal ambiguity sets $\mathcal{P}^n$ under which $p_n (s'_n = 0 | s, a) = p_n (s'_n = 1 | s, a) = 1/2$ for all even sub-states $n \in [N]$ as well as $\mathcal{P}^n = \Delta (\mathcal{S}_n)$ for all odd sub-states $n \in [N]$, and unit rewards whenever there is $\eta \in [ \lfloor N / 2 \rfloor ]$ such that $s_{2 \eta} \neq s_{2 \eta + 1}$. In this case, $\mathrm{NF}^\star (\mathcal{R}) = \mathrm{L}^\star (\mathcal{R}) = 0$ since the worst-case distribution can always ensure that $s_{2 \eta} = s_{2 \eta + 1}$ for all $\eta \in [ \lfloor N / 2 \rfloor ]$ with probability $1$. We have $\mathrm{F}^\star (\mathcal{R}) = \gamma \cdot (1 - \gamma)^{-1} \cdot (1 - 2^{-\lfloor N / 2 \rfloor})$, on the other hand, since no matter which distribution from the ambiguity set is selected, the probability that $s_{2 \eta} = s_{2 \eta + 1}$ for all $\eta \in [ \lfloor N / 2 \rfloor ]$ is $1 / 2^{\lfloor N / 2 \rfloor}$ because all even sub-states transition independently.

    In view of the third statement, finally, consider the family of RFMDP instances with state space $\mathcal{S} = \prod_{n \in [N]} \mathcal{S}_n$ and $\mathcal{S}_n = \mathbb{B}$, $n \in [N]$, action space $\mathcal{A} = \{ 1 \}$, $\phi (s, a) = (s, a)$, a deterministic initial state $s^0$ satisfying $s^0_n = \mathbf{1} [n \text{ is even}]$, $n \in [N]$, marginal ambiguity sets $\mathcal{P}^n$ under which $p_n (s'_n = 0 | s, a) = p_n (s'_n = 1 | s, a) = 1/2$ for all $n \in [N]$ and unit rewards whenever $s \neq (0, \ldots, 0)$ and $s \neq (1, \ldots, 1)$. In this case, $\mathrm{NF}^\star (\mathcal{R}) = 0$ since the worst-case distribution transitions into the two states $s = (0, \ldots, 0)$ and $s = (1, \ldots, 1)$ with probability $1/2$ each. We have $\mathrm{L}^\star (\mathcal{R}) = \mathrm{F}^\star (\mathcal{R}) = \gamma \cdot (1 - \gamma)^{-1} \cdot (1 - 2^{-(N-1)})$, on the other hand, since at any time $t = 2, \ldots$, the new state is $s = (0, \ldots, 0)$ or $s = (1, 
    \ldots, 1)$ with probability $1 / 2^N$ each because all sub-states transition independently.
\qed

\end{document}